\documentclass[11pt,reqno]{amsart}
 \usepackage{amsmath,amssymb,amsthm}
\usepackage[margin=1.1in]{geometry}
\usepackage{booktabs,array,dcolumn}
\usepackage{tikz}
\usepackage{parskip}
\usepackage[colorlinks=true,linkcolor=blue,citecolor=blue,urlcolor=blue]{hyperref}
\usepackage{tcolorbox}
\tcbuselibrary{skins}
\newtcolorbox{mainresult}{enhanced,frame hidden,colback=white,boxrule=0.5pt,
  borderline={0.5pt}{0.25pt}{black},arc=0pt,outer arc=0pt,left=7pt,right=7pt,top=5pt,bottom=5pt,
  before skip=12pt plus 2pt,after skip=12pt plus 2pt}

\newtheorem{theorem}{Theorem}[section]
\newtheorem{lemma}[theorem]{Lemma}
\newtheorem{proposition}[theorem]{Proposition}
\newtheorem{corollary}[theorem]{Corollary}

\theoremstyle{definition}
\newtheorem{definition}[theorem]{Definition}
\newtheorem{example}[theorem]{Example}
\theoremstyle{remark}
\newtheorem{remark}[theorem]{Remark}
\theoremstyle{plain}
\newtheorem{thmE}{Theorem}

\newtheorem{lemE}{Lemma}

\numberwithin{equation}{section}

\newcommand{\Z}{\mathbb{Z}}
\newcommand{\Q}{\mathbb{Q}}
\newcommand{\C}{\mathbb{C}}
\newcommand{\R}{\mathbb{R}}
\newcommand{\eps}{\varepsilon}
\newcommand{\vk}{\varkappa}

\newcommand{\lcm}{\operatorname{lcm}}
\newcommand{\Log}{\operatorname{Log}}
\newcommand{\im}{\operatorname{Im}}
\newcommand{\re}{\operatorname{Re}}
\newcommand{\fl}[1]{\lfloor #1\rfloor}
\newcommand{\fr}[1]{\{#1\}}
\newcommand{\rt}{\sqrt{2}}
\newcommand{\dd}{\mathrm{d}}
\newcommand{\Dom}{\Omega}
\newcommand{\cA}{\mathcal{A}}
\newcommand{\lev}[1]{\mathfrak{L}_{#1}}
\newcommand{\Zp}{\Z_{(p)}}
\newcommand{\XX}[2]{X_{#1}^{(#2)}}
\newcommand{\dg}{\Psi}

\makeatletter
\def\@setauthors{\begingroup
  \parskip\z@
  \trivlist
  \centering\footnotesize \@topsep30\p@\relax
  \advance\@topsep by -\baselineskip
  \item\relax
  \author@andify\authors
  {\normalsize\bfseries\authors\par}\vspace{3pt}{\footnotesize\textsuperscript{1}ETH Zurich, Switzerland\par}{\footnotesize\textsuperscript{2}University of Fribourg, Switzerland\par}\vspace{6pt}{\footnotesize\itshape First released October 3, 2026$^\dagger$\enspace\textperiodcentered\enspace This version: \today\ (version~5)\par}\endtrivlist
  \endgroup}
\def\@evenhead{\normalfont\scriptsize\rlap{\thepage}\hfil}
\def\@settitle{\begin{center}\Large\bfseries\boldmath
  \@title
  \end{center}}
\def\section{\@startsection{section}{1}\z@{.8\linespacing\@plus\linespacing}{.5\linespacing}{\normalfont\large\bfseries\boldmath\raggedright\def\@secnumfont{\bfseries}}}
\makeatother
\begin{document}

\title[An integral for $\pi$ and the bound $\mu(\pi)\le 6.0446$]{A family of integrals for $\pi$ with the symmetry $z\mapsto 2/z$,\\ and the bound $\mu(\pi)\le 6.0446$}
\author[Michael Moor]{Michael Moor\textsuperscript{1}}
\author[Bastian Rieck]{Bastian Rieck\textsuperscript{2}}
\thanks{$^\dagger$Public release on GitHub: \url{https://github.com/mi92/mu-pi-pub}}
\hypersetup{pdfauthor={Michael Moor, Bastian Rieck},pdftitle={A family of integrals for pi with the symmetry z -> 2/z, and the bound mu(pi) <= 6.0446}}
\date{}

\begin{abstract}
We prove that the irrationality measure of $\pi$ satisfies $\mu(\pi)\le 6.0446$,
improving the bound $7.1032\ldots$ of Zeilberger and Zudilin~(2020). Our proof uses the integrals
\[
J_n=-i\int_{1-i}^{1+i}\biggl(\frac{\bigl((z-1)(z-2)\bigr)^{5}\,(z^2-2z+2)^{3}}{z^{4}\,(z^2-2)^{4}}\biggr)^{\!n}\,\frac{\dd z}{z}=u_n+v_n\pi,
\]
where the rational function in brackets is invariant under $z\mapsto 2/z$. We bound the
denominators of the rational numbers $u_n,v_n$ using primes that can be removed from
$\lcm(1,\dots,4n)$, determine the exponential growth of $v_n$, and bound the decay of $J_n$. For $n=1$, this integral recovers Zu Chongzhi's
approximation $355/113$. We also give an effective version: $|\pi-a/b|>b^{-6.0455}$ for all
integers $a$ and $b\ge2$, improving Mignotte's exponent $20.6$ (1974), and an explicit $b_0(\mu)$
for every $\mu>6.04459707$. The bound
$6.0446$ is formalised in Lean~4 with Mathlib, whereas the formula
for $b_0(\mu)$ is formalised assuming four published estimates for Chebyshev's functions.
The discovery of the integrand was largely carried out by an AI agent under human supervision.
\end{abstract}

\maketitle

\section{Introduction}\label{sec:intro}

The following quantity measures the rational approximability of irrational numbers.
\begin{definition}[Irrationality measure]\label{def:mu}
For $\theta\in\R\setminus\Q$, the \emph{irrationality measure} (or \emph{irrationality exponent}) of $\theta$ is
\[
\mu(\theta)=\inf\Bigl\{\nu\in\R:\ \text{there is an integer $b_0\ge1$ such that }\Bigl|\theta-\frac ab\Bigr|>\frac1{b^{\nu}}\ \text{for all $a\in\Z$, $b\in\Z_{\ge b_0}$}\Bigr\}.
\]
\end{definition}

By Dirichlet's approximation theorem, $\mu(\theta)\ge2$ for every irrational
number. While equality holds for Lebesgue-almost every $\theta$, the exact value
of $\mu(\pi)$ has been an open problem for decades, with Mahler's first upper
bound of $\mu(\pi) \leq 42$ published in 1953~\cite{Mah53}, and
Table~\ref{tab:history} listing the main improvements since. Starting from Hata's
bound $8.0160\dots$~\cite{Hata}, every new bound was obtained from complex
integrals of rational functions whose values are linear forms in $1$ and $\pi$,
respectively: Salikhov, for instance, introduced a symmetrised integral, while
Zeilberger and Zudilin found a better member of Salikhov's family. Bai's recent
preprint~\cite{Bai} lowers the bound to $7.1018\dots$ within said family, and a further
preprint~\cite{NG} conjectures that this is the best bound obtainable from it; we leave
this family and instead use integrands with a different symmetry, namely $z\mapsto 2/z$.

\begin{table}[tbp]
\small
\begin{tabular}{llD{.}{.}{2.7}l}
\toprule
year & author & \multicolumn{1}{c}{$\mu(\pi)\le$} & method\\
\midrule
1953 & Mahler \cite{Mah53} & 42 & Hermite's approximations to $e^z$\\
1974 & Mignotte \cite{Mig74} & 20 & refinement of Mahler's method\\
1982 & G.~Chudnovsky \cite{Chud82} & 19.8899\ldots & Hermite--Pad\'e approximations\\
1993 & Hata \cite{HataJNT} & 13.398 & Hermite-type approximations to $\log^j(1-x)$\\
1993 & Hata \cite{Hata} & 8.0160\ldots & the integrals $\int\bigl((z-1)(z-2)(z-1-i)\bigr)^{2n}z^{-3n-1}\dd z$\\
2008 & Salikhov \cite{Sal08,Sal10} & 7.6063\ldots & symmetrised complex integral\\
2020 & Zeilberger--Zudilin \cite{ZZ} & 7.1032\ldots & another member of Salikhov's family\\
2026 & Bai \cite{Bai} (preprint) & 7.1018\ldots & deformation inside the same family\\
2026 & this paper (Oct 3, 2026) & 6.0446 & integrals with the symmetry $z\mapsto2/z$\\
\bottomrule
\end{tabular}
\medskip
\caption{Upper bounds for the irrationality measure of $\pi$, truncated to four decimal places. ``$\ldots$'' marks a bound with further
digits, e.g., $19.8899944\ldots$ (Chudnovsky). The values for Mahler and Mignotte are those of~\cite{Mah53,Mig74}, the value for Chudnovsky is $1+\chi_2$ with the $\chi_2=18.88999444\ldots$ of~\cite[Corollary~3.6]{Chud82}, and the value $13.398$ is the one stated in the published paper~\cite{HataJNT}. Mahler proved $|\pi-a/b|>b^{-42}$ for all $b\ge2$ and remarked that $42$ can be replaced by $30$ for large $b$; Mignotte proved $|\pi-a/b|>b^{-20.6}$ for all $b\ge2$, and his exponent $20$ holds for $b\ge\exp(e^{110})$; Chudnovsky stated the exponent $19.9$ for all $b\ge2$, to the best of our knowledge without proof (Section~\ref{ss:eff-results}). A recent preprint \cite{NG} conjectured that $7.1018\ldots$ is the best bound obtainable from the family of \cite{Sal10,ZZ,Bai}.}\label{tab:history}
\end{table}

\subsection{Results}

Our main result, whose proof is sketched in Section~\ref{sec:roadmap},  can be stated as follows:
\begin{mainresult}
\begin{theorem}\label{thm:main}
The irrationality measure of $\pi$ satisfies $\mu(\pi)\le6.0446$.
\end{theorem}
\end{mainresult}
Previous upper bounds comprise $7.1018\dots$~\cite{Bai} and $7.1032\dots$~\cite{ZZ}.
Section~\ref{sec:effective}, unlike some prior work, provides a value for the
threshold $b_0$ from Definition~\ref{def:mu}.
However, a simplified version of the statement~(Theorem~\ref{thm:E1}) holds without any thresholds, namely
\[
\Bigl|\pi-\frac ab\Bigr|>b^{-6.0455}
\]
for all integers $a, b$ with $b \geq 2$, while Theorem~\ref{thm:E2} provides a closed formula for $b_0(\mu)$ for every $\mu > 6.04459707$.
To the best of our knowledge, the last exponent that was proved for all $b\ge2$ is the $20.6$ of
Mignotte~\cite{Mig74} from 1974, which improved Mahler's $42$~\cite{Mah53}, and the exponent $19.9$ that
Chudnovsky~\cite{Chud82} stated for all $b\ge2$ in 1982 was given without proof; the bounds of
Table~\ref{tab:history} from 1993 on are stated for $b\ge b_0$ without a value of $b_0$
(Section~\ref{ss:eff-results}).

\subsection{Linear forms and the integrals}\label{ss:linforms}

Several recent entries of Table~\ref{tab:history} are proved in a similar way. One constructs integers $q_n$ and $p_n$ ($n=1,2,\dots$) such that
the \emph{linear form} $q_n\pi-p_n$ is small while $q_n$ is not too large. The following lemma converts the two exponential rates into a bound for the irrationality measure.

\begin{lemma}[Hata {\cite[Remark~2.1]{Hata}}]\label{lem:hata}
Let $\theta$ be a real irrational number and let $q_n$, $p_n$ $(n\ge1)$ be integers such that
\[
\lim_{n\to\infty}\frac1n\log|q_n|=\sigma\qquad\text{and}\qquad\limsup_{n\to\infty}\frac1n\log|q_n\theta-p_n|\le-\tau
\]
for some real numbers $\sigma>0$ and $\tau>0$. Then $\mu(\theta)\le1+\sigma/\tau$.\footnote{Hata states the lemma as a remark to his Lemma~2.1, which concerns two linear forms and is proved in full there. In his terminology, the conclusion reads ``$\theta$ has an irrationality measure $1+\sigma/\tau$,'' which means $\mu(\theta)\le1+\sigma/\tau$ in the sense of Definition~\ref{def:mu}. A proof in the form used by the formal development is part of the Lean project (Section~\ref{sec:formal}).
}
\end{lemma}

For the bound of Zeilberger and Zudilin, $\sigma=11.6138\ldots$ and
$\tau=1.9029\ldots$, yielding $\mu(\pi)\le7.1032\ldots$. By comparison, in this paper, $\sigma\le14.2464$ and $\tau\ge2.824091$, giving $\mu(\pi)\le6.0446$.
Our linear forms come from the following rational function
\begin{equation}\label{eq:f}
f(z)=\frac{\bigl((z-1)(z-2)\bigr)^5\,(z^2-2z+2)^3}{z^4\,(z^2-2)^4},
\end{equation}
where  $z \in \Dom:=\C\setminus\{0,\rt,-\rt\}$.
We write $N(z)=\bigl((z-1)(z-2)\bigr)^5(z^2-2z+2)^3\in\Z[z]$ for its numerator. For every integer $n\ge1$, let
\begin{equation}\label{eq:J}
J_n=-i\int_{1-i}^{1+i}f(z)^n\,\frac{\dd z}{z},
\end{equation}
the integral being taken along the straight segment from $1-i$ to $1+i$.
The numerator $N(z)$ vanishes at $1$, $2$ and $1\pm i$, the function $f$ has poles at $0$ and $\pm\rt$, and $f(2/z)=f(z)$ (Lemma~\ref{lem:sym}). The factor $-i$ ensures that $J_n\in\R$ (Section~\ref{ss:real}).

\begin{definition}\label{def:removable}
Let $n\ge1$. A prime number $p$ is \emph{removable for $n$} if
\[
p^2>4n,\qquad \fl{4n/p}\ \text{is even},\qquad (5n\bmod p)+(3n\bmod p)<(4n\bmod p).
\]
We denote by $\Phi_n$ the product of all removable primes $p\le4n$, and we put
\[
D_n=\frac{\lcm(1,2,\dots,4n)}{\Phi_n}.
\]
\end{definition}

Since a removable prime $p\le4n$ satisfies $p^2>4n$, it divides $\lcm(1,\dots,4n)$ exactly once. Thus, $D_n$ is an integer. We shall see
(Lemma~\ref{lem:windows}) that a prime $p>\sqrt{4n}$ is removable if and only if $\fr{n/p}\in[\frac15,\frac14)\cup[\frac23,\frac34)$, where
$\fr{x}=x-\fl{x}$ denotes the fractional part of a real number $x$. For example $\Phi_1=1$, $\Phi_2=3$ (because $\fr{2/3}=\frac23$) and $\Phi_{21}=17\cdot29\cdot31$.

Removing primes by a condition on $\fr{n/p}$ is a standard device: Hata~\cite[Lemma~2.2]{Hata} used $\frac12\le\fr{n/p}<\frac23$ for his
approximations to $\pi$, Rhin and Viola~\cite[Section~4]{RV96} used conditions of the same kind for $\zeta(2)$, and Zeilberger and
Zudilin~\cite[Lemmas~2 and~3]{ZZ} divided a product of primes from Hata's interval out of $\lcm(1,\dots,4n)$ in their construction for $\pi$.

The last two conditions of Definition~\ref{def:removable} have the following meaning (the proof is in Section~\ref{ss:removable}).
Write $5n$, $3n$ and $4n$ as multiples of $p$ plus remainders $\alpha_1$, $\alpha_2$, $\delta$ (each less than $p$). Since the exponents of the
integrand are $5n$, $3n$ and $4n$, the integrand $h_n:=f^n/z$ then splits as a $p$-th power times a small factor, $h_n=g^p\cdot\omega$, where $\omega$ contains
$(z-1)(z-2)$, $z^2-2z+2$ and $z^2-2$ only with the small exponents $\alpha_1$, $\alpha_2$, $\delta$. Because all these exponents are less than $p$ (for
$\alpha_1+\alpha_2$ this follows from $\alpha_1+\alpha_2<\delta$), $\omega$ can be integrated without dividing by $p$. The parity of $\fl{4n/p}$
makes $\delta$ even, which this integration needs. The inequality $\alpha_1+\alpha_2<\delta$ also means that $\omega$ has no pole at $0$,
so that its primitive contains no logarithm. Finally, the derivative
of a $p$-th power is divisible by $p$, since $(g^p)'=p\,g^{p-1}g'$ (the same reason why $(x^p)'=px^{p-1}$ vanishes modulo~$p$). Integrating by parts,
this factor $p$ cancels the division by $p$.
Removing these primes is worth a lot. With the multiplier $\lcm(1,\dots,4n)$ the method would only give $\mu(\pi)\le6.747$, while with $D_n$
it gives $6.0446$.

\begin{theorem}[the linear forms]\label{thm:forms}
For every integer $n\ge1$:
\begin{enumerate}
\item[(a)] $J_n=u_n+v_n\pi$ with rational numbers $u_n$ and $v_n$. More precisely $v_n=\frac12c_n$, where $c_n$ is the coefficient of
$z^{4n}$ in the power series expansion of $N(z)^n(z^2-2)^{-4n}$ at $z=0$, and $c_n$ is a positive integer.
\item[(b)] $D_nu_n$ and $D_nv_n$ are integers.
\item[(c)] The limit $r=\lim_{n\to\infty}\frac1n\log v_n$ exists, and $8.9<r\le10.540357$.
\item[(d)] $|J_n|\le4\rho^n$ with $\rho=0.001458809<e^{-6.530134}$.
\item[(e)] $\lim_{n\to\infty}\frac1n\log D_n=4-\varpi$, where
\[
\varpi=\sum_{k=0}^\infty\frac{1}{(3k+2)(4k+3)}+\sum_{k=1}^\infty\frac{1}{(4k+1)(5k+1)}=0.2939598992\ldots
\]
\end{enumerate}
\end{theorem}

Theorem~\ref{thm:main} follows from Theorem~\ref{thm:forms} and Lemma~\ref{lem:hata} applied to the integers $q_n=D_nv_n$ and $p_n=-D_nu_n$. One gets $\sigma=r+(4-\varpi)\le14.2464$ and $\tau\ge6.530134-(4-\varpi)\ge2.824091$, and
$1+14.2464/2.824091<6.0446$. The prime number theorem is used only to determine the size of the multiplier (part~(e), and Proposition~\ref{prop:mult} below).
The details are in Section~\ref{sec:conclusion}. There, as in the formal proof, we work with a multiplier $D_n^{(K)}$ that removes from
$\lcm(1,\dots,4n)$ only the primes of finitely many ``windows'' (Proposition~\ref{prop:mult}(a)), so that part~(e) is needed only in a
truncated form. It gives the same bound.
The exact values of the three rates are discussed in Remark~\ref{rem:exact}. They would give $6.04458594\ldots$, so nothing is lost
in the rounding.

\begin{example}[$n=1$]\label{ex:n1}
For $n=1$ a direct computation (formula \eqref{eq:cformula} for $c_1$ and the recursion of Section~\ref{ss:rec} for $u_1$) gives
\[
J_1=\frac{7345\,\pi-23075}{2}=\frac{7345}{2}\Bigl(\pi-\frac{355}{113}\Bigr)=-0.000979691\ldots ,
\]
because $7345=65\cdot113$ and $23075=65\cdot355$. So the first integral reproduces the approximation $355/113$. Here $D_1=\lcm(1,2,3,4)=12$,
and indeed $12u_1=-138450$ and $12v_1=44070$ are integers.
\end{example}

\subsection{What is new}\label{ss:new}
The integrals \eqref{eq:J} give $\mu(\pi)\le6.0446$ (Theorem~\ref{thm:main}). The proof has been formalised in Lean, with the prime number
theorem as its only hypothesis (Section~\ref{sec:formal}). Section~\ref{sec:effective} also proves $|\pi-a/b|>b^{-6.0455}$ for every integer
$a$ and every integer $b\ge2$.

Hata's integrand of 1993 has the zeros $1$, $2$ and $1+i$ and a pole at $0$
only. The integrand used here has in addition the zero $1-i$ and the two poles $\pm\rt$, which are fixed by the symmetry $z\mapsto2/z$ and contribute no
logarithm to the integral (Remark~\ref{rem:classical}). Three facts make the construction work. First, the only ``bad'' prime is $2$, and its contribution to the denominators is only $2^{O(\log n)}$
(Section~\ref{ss:two}). Second, many primes up to $4n$ can be removed from the denominators, their product being $e^{0.2939\ldots n}$
(Sections~\ref{ss:removable}, \ref{sec:mult}). Third, the integral is exponentially smaller than the straight path of integration suggests
(Section~\ref{sec:decay}).

The integrand $f$ is one member of a three-parameter family of integrands with the same symmetry, which also contains a member apparently equivalent
to the integrals of Zeilberger and Zudilin. Section~\ref{sec:family} discusses this family. Apart from one identity, nothing in it is proved.

To the best of our knowledge, these integrals have not appeared in the literature before.

\subsection{Status of the result}\label{ss:status}
To our knowledge the result is new. It improves the published record, which has stood since 2020, by a wide margin. The work was largely executed by an AI agent under human supervision. The proof has been formalised in Lean~4 with Mathlib. The formal statement is Theorem~\ref{thm:main} under the hypothesis of the classical result
$\theta(x)\sim x$ (the prime number theorem, which is not yet part of Mathlib). Every other step of the proof of Theorem~\ref{thm:main},
including all numerical estimates that it uses, is machine-checked (in a few places by a slightly different route, Section~\ref{sec:formal}). Several statements of this paper that are not needed for Theorem~\ref{thm:main} are
\emph{not} formalised, for instance Theorem~\ref{thm:forms}(e), the lower bound $8.9<r$ in
Theorem~\ref{thm:forms}(c) and the exact values of the constants. Section~\ref{ss:notcovered}
lists the main ones. Section~\ref{sec:formal} states the verified theorem exactly and describes the formal development. The proofs in this paper follow the formalised
route. Where it differs from the classical one we say so in a remark.
Section~\ref{sec:effective} makes the bound effective (Theorems~\ref{thm:E1} and~\ref{thm:E2}). Theorem~\ref{thm:E2}, and
Theorem~\ref{thm:E1} for $b\ge10^{122\,680\,112\,076\,457}$, are machine-checked under four published estimates for Chebyshev's functions, which
enter as hypotheses; for the smaller $b$ Theorem~\ref{thm:E1} rests on computations with recorded logs, certified in ball arithmetic
(Table~\ref{tab:trust}).

\subsection{Organisation}
Section~\ref{sec:roadmap} gives the proof in one page. Section~\ref{sec:integrand} describes the integrand.
Sections~\ref{sec:form}--\ref{sec:mult} prove the five parts of Theorem~\ref{thm:forms}, and Section~\ref{sec:conclusion} proves
Theorem~\ref{thm:main}. Section~\ref{sec:effective} makes the bound effective, with an explicit $b_0(\mu)$ (details in Appendix~\ref{app:effective}). Section~\ref{sec:formal} describes the formal verification and the numerical checks, and Section~\ref{sec:family}
places the integrand in a three-parameter family, without proofs. Appendix~\ref{app:details} contains the proofs of several lemmas and the
details of the estimates of Sections~\ref{sec:growth}, \ref{sec:decay} and~\ref{sec:mult}.

\subsection{Accompanying material}\label{ss:material}
The Lean project (directory \texttt{lean/}, Section~\ref{sec:formal}) and the programs that recompute the numbers printed in this paper
(directory \texttt{scripts/}, Section~\ref{sec:numerics}) are publicly available at \url{https://github.com/mi92/mu-pi-pub}.
Everything that the proof of Theorem~\ref{thm:main} given in Sections~\ref{sec:integrand}--\ref{sec:conclusion} and Appendix~\ref{app:details}
uses is printed in this paper, with one exception. The $99$ rational numbers of the certificate in Proposition~\ref{prop:contour} are too long to
print. They are contained in the file \texttt{lean/MuPi/K2H/ContourCert.lean} of that repository, and Appendix~\ref{app:decay} describes how to
recompute them. The programs and recorded logs of the effective version (Section~\ref{sec:effective}) are in the directory \texttt{effective/} of
the same repository (Section~\ref{ss:eff-repro}). It contains the certificate of Theorem~\ref{thm:E3}(i), which is too long to print, and the
programs that recompute it, the table of block bounds of Theorem~\ref{thm:E3}(ii) and the partial quotients of Theorem~\ref{thm:E3}(iii). The
table of block bounds itself ($650$\,MB) is not stored there.

\section{The proof in one page}\label{sec:roadmap}

Write $h_n(z)=f(z)^n/z$ ($z\in\Dom$) for the integrand of \eqref{eq:J}.

\begin{enumerate}
\item[\textbf{1.}] \textbf{Primitive pairs.} A pair $(G,c)$, with $G$ a holomorphic function on $\Dom$ and $c$ a complex number, is called a
\emph{primitive pair} if $G'(z)=h_n(z)-c/z$ on $\Dom$. Any two primitive pairs have the same $c$ and the same difference
$G(1+i)-G(1-i)$, and the fundamental theorem of calculus gives
\[
J_n=-i\bigl(G(1+i)-G(1-i)\bigr)+\frac{c\,\pi}{2},
\]
because $\log z$ increases by $i\pi/2$ along the segment. So $u_n=-i(G(1+i)-G(1-i))$ and $v_n=c/2$ can be computed from \emph{any} primitive pair
(Section~\ref{ss:pairs}). The whole proof consists in exhibiting several primitive pairs, each adapted to one purpose.
\item[\textbf{2.}] \textbf{Rationality.} In the variables $W=z+2/z$ and $Y=z-2/z$ one has $f=(W-3)^5(W-2)^3/Y^4$, and a two-term recursion
integrates every function $W^jY^{-2m}/z$. This gives a primitive pair with rational data, hence $u_n,v_n\in\Q$ (Section~\ref{ss:rec}).
\item[\textbf{3.}] \textbf{Denominators.} For every prime $p$ we exhibit a primitive pair that is ``$p$-integral'' up to an explicit power
of $p$.
For odd $p$, the partial fraction decomposition of $h_n$, integrated term by term, has denominators $\le4n$ (Section~\ref{ss:odd}).
For a removable prime, an identity $h_n=(g^pG_{\mathrm{red}})'-p\,g^{p-1}g'G_{\mathrm{red}}$ with suitable rational functions $g$ and
$G_{\mathrm{red}}$ shows that the factor $p$ in those denominators cancels (Section~\ref{ss:removable}).
For $p=2$, a M\"obius change of variable $z\mapsto w$ turns $h_n$ into a combination of elementary integrands whose primitives are divisible
by explicit powers of $2$ (Section~\ref{ss:two}). Together, $D_nu_n,D_nv_n\in\Z$ (Section~\ref{ss:integrality}).
\item[\textbf{4.}] \textbf{Growth.} After the substitution $z\mapsto-z$ the number $c_n$ is the coefficient of $z^{4n}$ in the $n$-th power of a
power series with non-negative coefficients. Hence $c_{m+n}\ge c_mc_n$, the limit $r$ exists (Fekete's lemma), and $c_n\le f(-x)^n$ for
$0<x<\rt$. With $x=0.3913$ this gives $r\le10.540357$ (Section~\ref{sec:growth}).
\item[\textbf{5.}] \textbf{Decay.} The primitive $G+c\Log z$ (with the principal logarithm $\Log$) also computes the integral of $h_n$ along the polygon
$1-i\to\bar z_s\to1\to z_s\to1+i$ through a point $z_s$ next to a saddle point of $|f|$. On this polygon $|f|\le\rho$, which is a polynomial
inequality in one real variable and is certified by a Bernstein expansion (Appendix~\ref{app:decay}). Hence $|J_n|\le4\rho^n$ (Section~\ref{sec:decay}).
\item[\textbf{6.}] \textbf{The multiplier.} By the prime number theorem $\log\lcm(1,\dots,4n)\sim4n$, and the removable primes are the primes $p>\sqrt{4n}$ in
the ``windows'' $\fr{n/p}\in[\frac15,\frac14)\cup[\frac23,\frac34)$, so $\frac1n\log D_n\to4-\varpi$ (Section~\ref{sec:mult}).
\item[\textbf{7.}] \textbf{Conclusion.} Hata's lemma with $q_n=D_nv_n$, $p_n=-D_nu_n$. In Section~\ref{sec:conclusion}, as in the formal proof, $D_n$ is
replaced by $D_n^{(K)}$, which is easier to estimate.
\end{enumerate}

\noindent The numbers fit together as shown in Table~\ref{tab:budget}.

\begin{table}[ht]
\small
\begin{tabular}{llll}
\toprule
rate & meaning & bound used & exact value\\
\midrule
$r$ & growth of $v_n$ & $r\le10.540357$ & $10.540355111696\ldots$\\
$s$ & decay of $J_n$: $|J_n|\le4e^{-sn}$ & $s=6.530134$ & $6.530136147360\ldots$ (saddle point)\\
$C$ & size of the multiplier: $e^{Cn+o(n)}$ & $C\le3.706043$ & $4-\varpi=3.706040100783\ldots$\\
\midrule
$\sigma=r+C$ & growth of $q_n$ & $\sigma\le14.2464$ & $14.2463952124\ldots$\\
$\tau=s-C$ & decay of $q_n\pi-p_n$ & $\tau\ge2.824091$ & $2.8240960465\ldots$\\
$1+\sigma/\tau$ & bound for $\mu(\pi)$ & $<6.0446$ & $6.04458594\ldots$\\
\bottomrule
\end{tabular}
\medskip
\caption{The three rates and how they combine. The column ``bound used'' contains what is proved, and formally verified, in this paper. There
the multiplier is $D_n^{(K)}$ with $K=10^5$ (Section~\ref{sec:mult}), for which $C=4-\varpi_K\le3.706043$ (this exceeds $4-\varpi$ by less than $3\cdot10^{-6}$, the price of truncating the series for $\varpi$ and of rounding the partial sum $\varpi_K\ge0.293957$ down). The last column contains the exact values of the three rates, which are discussed in Remark~\ref{rem:exact} and are not needed for the proof. The value of $C$ is proved in Proposition~\ref{prop:mult}(b) (assuming the prime number theorem, for the full multiplier $D_n$), and $r=\log\lambda_3$, with $\lambda_3$ the minimum of $f(-x)$ (Remark~\ref{rem:r}), is proved in Theorem~\ref{thm:E3}(i). The value of $s$ is computed numerically, from the saddle point of $|f|$ (Remark~\ref{rem:saddle}), and the paper proves nothing more about it than the bound $|J_n|\le4e^{-6.530134\,n}$ of the middle column.}\label{tab:budget}
\end{table}

\section{The integrand}\label{sec:integrand}

\subsection{Zeros, poles and the symmetry}
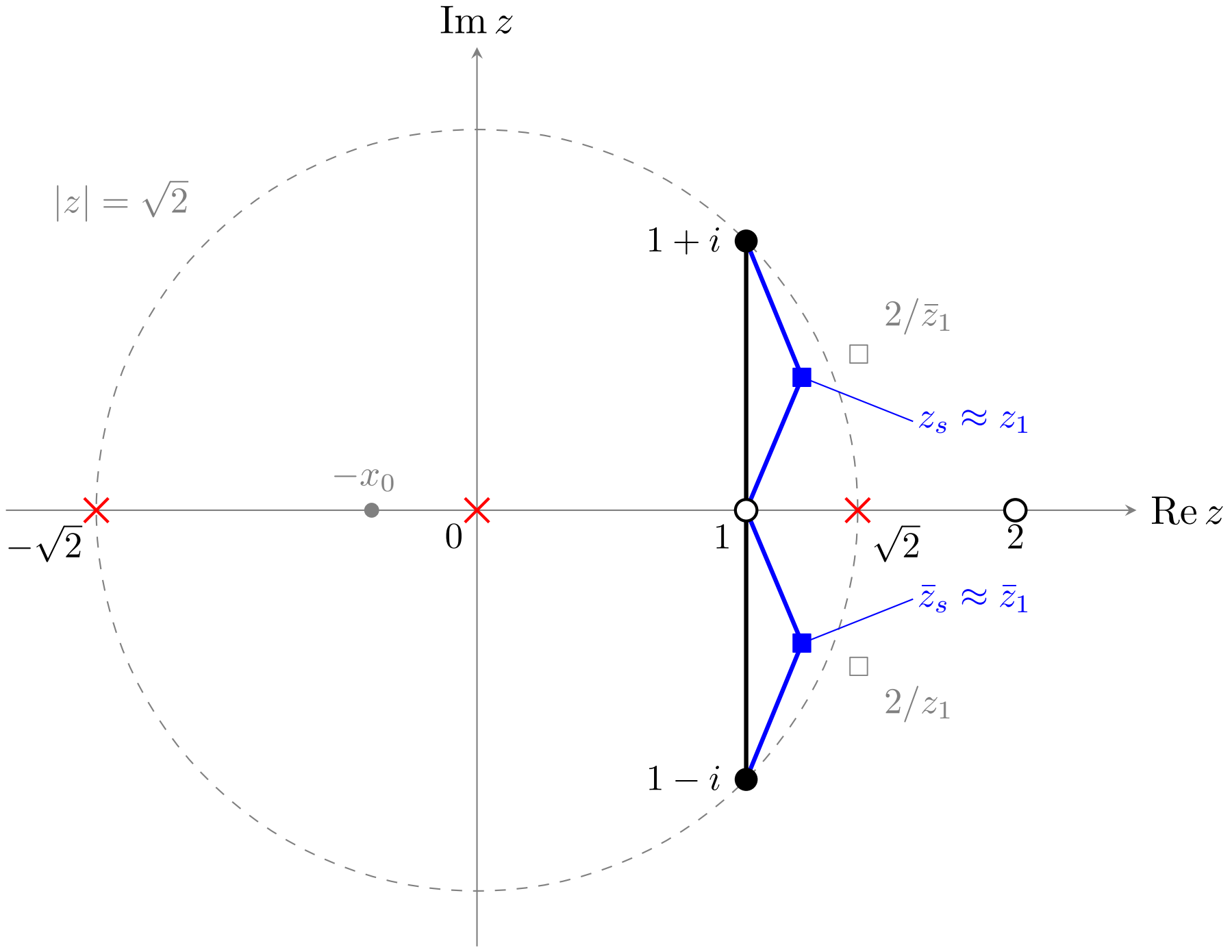
\begin{figure}[ht]
\begin{tikzpicture}[scale=2.6,>=stealth]
  \draw[->,gray] (-1.75,0)--(2.45,0) node[right,black] {$\re z$};
  \draw[->,gray] (0,-1.62)--(0,1.72) node[above,black] {$\im z$};
  \draw[dashed,gray] (0,0) circle (1.41421);
  \node[gray,above left] at (-1.02,1.02) {\small $|z|=\rt$};
\draw[very thick] (1,-1)--(1,1);
\draw[very thick,blue] (1,-1)--(1.2077,-0.4939)--(1,0)--(1.2077,0.4939)--(1,1);
\foreach \x/\lab/\pos in {0/0/below left,1.41421/\rt/below right,-1.41421/-\rt/below left}
    {\draw[thick,red] (\x-0.045,-0.045)--(\x+0.045,0.045) (\x-0.045,0.045)--(\x+0.045,-0.045);
     \node[\pos] at (\x,0) {\small $\lab$};}
\draw[fill=white,thick] (1,0) circle (0.04) node[below left] {\small $1$};
  \draw[fill=white,thick] (2,0) circle (0.04) node[below] {\small $2$};
  \draw[fill=black] (1,1) circle (0.04) node[left=3pt] {\small $1+i$};
  \draw[fill=black] (1,-1) circle (0.04) node[left=3pt] {\small $1-i$};
\node[blue] at (1.2077,0.4939) {\tiny$\blacksquare$};
  \node[blue] at (1.2077,-0.4939) {\tiny$\blacksquare$};
  \draw[blue,thin] (1.2077,0.4939)--(1.62,0.33);
  \draw[blue,thin] (1.2077,-0.4939)--(1.62,-0.33);
  \node[blue,right,inner sep=1pt] at (1.62,0.33) {\small $z_s\approx z_1$};
  \node[blue,right,inner sep=1pt] at (1.62,-0.33) {\small $\bar z_s\approx\bar z_1$};
  \node[gray] at (1.4187,0.5802) {\tiny$\square$};
  \node[gray] at (1.4187,-0.5802) {\tiny$\square$};
  \node[gray,right] at (1.46,0.72) {\small $2/\bar z_1$};
  \node[gray,right] at (1.46,-0.72) {\small $2/z_1$};
  \draw[fill=gray,gray] (-0.3913,0) circle (0.025);
  \node[gray,above] at (-0.42,0.02) {\small $-x_0$};
\end{tikzpicture}
\caption{The $z$-plane. Crosses mark the poles $0,\pm\rt$ of $f$, open circles the zeros $1$, $2$ (order~$5$) and black dots the zeros $1\pm i$
(order~$3$), which are the end points of the path. The black segment is the path of integration in \eqref{eq:J}, and the blue line is the polygon $1-i\to\bar z_s\to1\to z_s\to1+i$ of Section~\ref{sec:decay}. Its vertex
$z_s$ is a rational point at distance less than $10^{-5}$ from the saddle point $z_1=1.2077099\ldots+0.4939472\ldots i$ of $|f|$, so the two are drawn
at the same place. The symmetry $z\mapsto2/z$ maps the dashed circle to itself. The points $z_1$ and $\bar z_1$ lie inside it, and their images, the saddle points
$2/z_1$ and $2/\bar z_1$ (grey squares), lie outside it. The grey dot is the point $-x_0=-0.39126\ldots$ where $|f|$ is minimal on the interval $(-\rt,0)$ (Section~\ref{sec:growth}).}\label{fig:plane}
\end{figure}

The numerator $N(z)=(z-1)^5(z-2)^5(z-1-i)^3(z-1+i)^3$ has degree $16$, so
$f=N/(z^4(z-\rt)^4(z+\rt)^4)$ has zeros of order $5$ at $1$ and $2$, zeros of order $3$ at the end points $1\pm i$ of the path, and poles of order $4$
at $0$, $\rt$, $-\rt$ (and at $\infty$). Figure~\ref{fig:plane} shows these points.

\begin{lemma}[symmetry]\label{lem:sym}
$f(2/z)=f(z)$ for all $z\in\Dom$.
\end{lemma}
\begin{proof}
Replace $z$ by $2/z$ in each factor:
\[
\Bigl(\frac2z-1\Bigr)\Bigl(\frac2z-2\Bigr)=\frac{2(z-1)(z-2)}{z^2},\quad
\Bigl(\frac2z\Bigr)^2-2\cdot\frac2z+2=\frac{2(z^2-2z+2)}{z^2},\quad
\Bigl(\frac2z\Bigr)^2-2=-\frac{2(z^2-2)}{z^2}.
\]
Hence the numerator of $f(2/z)$ is $2^{5+3}N(z)/z^{16}$ and the denominator is $(2/z)^4\cdot2^4(z^2-2)^4/z^8=2^8z^{-12}(z^2-2)^4$. The quotient is $f(z)$.
\end{proof}

The map $z\mapsto2/z$ exchanges $1$ and $2$, exchanges $1+i$ and $1-i$, exchanges $0$ and $\infty$, and fixes $\pm\rt$. It maps the circle
$|z|=\rt$, which passes through $1\pm i$ and $\pm\rt$, to itself. The integrand of \eqref{eq:J} is built to respect this symmetry. This is the
reason why the two fixed points $\pm\rt$ will turn out to be harmless poles (they contribute no logarithm, Remark~\ref{rem:classical}).

\subsection{The invariant coordinates \texorpdfstring{$W$ and $Y$}{W and Y}}\label{ss:WY}
The symmetry suggests working with a function that is invariant under $z\mapsto 2/z$ and one that changes sign:
\begin{equation}\label{eq:WY}
W=z+\frac2z,\qquad Y=z-\frac2z .
\end{equation}

\begin{lemma}\label{lem:WY}
On $\C\setminus\{0\}$:
\begin{enumerate}
\item[(a)] $W^2-Y^2=8$;
\item[(b)] $z\,\dfrac{\dd W}{\dd z}=Y$ and $z\,\dfrac{\dd Y}{\dd z}=W$;
\item[(c)] $W-3=\dfrac{(z-1)(z-2)}{z}$, \ $W-2=\dfrac{z^2-2z+2}{z}$, \ $Y=\dfrac{z^2-2}{z}$, and therefore
\begin{equation}\label{eq:fWY}
f=\frac{(W-3)^5\,(W-2)^3}{Y^4}\qquad\text{on }\Dom;
\end{equation}
\item[(d)] at the end points, $W(1\pm i)=2$ and $Y(1\pm i)=\pm2i$.
\end{enumerate}
\end{lemma}
The proof is a direct computation. Formula \eqref{eq:fWY} explains the shape of $f$. In the variable $W$ the zeros are at
$W=3$ (that is $z=1,2$) and $W=2$ (that is $z=1\pm i$), and $Y^2=W^2-8$ vanishes at $W=\pm2\rt$ (that is $z=\pm\rt$).

\subsection{The integral is real}\label{ss:real}
Parametrise the segment by $z=1+it$, $-1\le t\le1$, so that $\dd z=i\,\dd t$ and
$J_n=\int_{-1}^1h_n(1+it)\,\dd t$. Since $f$ has real coefficients, $h_n(\bar z)=\overline{h_n(z)}$, so the values of the integrand at $t$ and
$-t$ are complex conjugates and $J_n$ is real.

\section{The linear form \texorpdfstring{$J_n=u_n+v_n\pi$}{Jn = un + vn pi}}\label{sec:form}

\subsection{Primitive pairs}\label{ss:pairs}
\begin{definition}\label{def:pair}
Let $n\ge1$. A \emph{primitive pair} (for $n$) is a pair $(G,c)$ consisting of a function $G\colon\Dom\to\C$ and a number $c\in\C$ such that
$G$ is complex differentiable on $\Dom$ and
\[
G'(z)=\frac{f(z)^n}{z}-\frac cz\qquad\text{for all }z\in\Dom .
\]
\end{definition}

\begin{lemma}[uniqueness]\label{lem:unique}
If $(G_1,c_1)$ and $(G_2,c_2)$ are primitive pairs for the same $n$, then
\[
c_1=c_2\qquad\text{and}\qquad G_1(1+i)-G_1(1-i)=G_2(1+i)-G_2(1-i).
\]
\end{lemma}
\begin{proof}
The function $\Delta=G_1-G_2$ satisfies $\Delta'(z)=(c_2-c_1)/z$ on $\Dom$. Integrate $\Delta'$ over the circle $|z|=\frac12$, which lies in $\Dom$. The
integral of a derivative over a closed curve is $0$, while the integral of $1/z$ is $2\pi i$. Hence $c_1=c_2$, and then $\Delta'=0$ on $\Dom$.
The strip $\{z:\frac12<\re z<\frac54\}$ is convex, is contained in $\Dom$ (because $\rt>\frac54$) and contains both $1+i$ and $1-i$, and a function
with vanishing derivative on a convex open set is constant, so $\Delta(1+i)=\Delta(1-i)$.
\end{proof}

\begin{lemma}[the fundamental theorem of calculus]\label{lem:ftc}
If $(G,c)$ is a primitive pair for $n$, then
\[
J_n=-i\bigl(G(1+i)-G(1-i)\bigr)+c\cdot\frac{\pi}{2}.
\]
\end{lemma}
\begin{proof}
Let $\Log$ be the principal branch of the logarithm, which is holomorphic on the half-plane $\re z>0$ with derivative $1/z$. On the set
$\{\re z>0\}\setminus\{\rt\}$ the function $G+c\Log$ is a primitive of $h_n$, and the segment from $1-i$ to $1+i$ lies in this set. Therefore
\[
\int_{1-i}^{1+i}h_n(z)\,\dd z=G(1+i)-G(1-i)+c\bigl(\Log(1+i)-\Log(1-i)\bigr),
\]
and $\Log(1\pm i)=\frac12\log2\pm i\frac\pi4$, so the last bracket equals $i\pi/2$. Multiply by $-i$.
\end{proof}

Consequently, \emph{as soon as one primitive pair exists}, the numbers
\begin{equation}\label{eq:uv}
u_n:=-i\bigl(G(1+i)-G(1-i)\bigr),\qquad c_n:=c,\qquad v_n:=\tfrac12c_n
\end{equation}
do not depend on the pair and satisfy $J_n=u_n+v_n\pi$. This is how $u_n$ and $v_n$ are \emph{defined} in this paper. We shall use
three primitive pairs, namely one with rational data (Section~\ref{ss:rec}), one from the partial fractions of $h_n$, which serves all odd primes (Sections~\ref{ss:odd},
\ref{ss:removable}), and one for the prime $2$ (Section~\ref{ss:two}). By Lemma~\ref{lem:unique} each of them gives information about the
same numbers $u_n$ and $c_n$.

\begin{remark}[the classical formulation]\label{rem:classical}
In the language of residues, a primitive pair exists if and only if the residues of $h_n$ at $\pm\rt$ vanish, and then $c$ is the residue at $0$. Classically one proves the
vanishing from the symmetry. The form $h_n\,\dd z$ changes sign under $z\mapsto 2/z$, and $\pm\rt$ are fixed points, so the residue there
equals its own negative. In the route followed here the vanishing is a by-product of the explicit construction in Section~\ref{ss:rec}.
No arithmetic property of $\pi$ is used in the proof of Theorem~\ref{thm:forms}.
The irrationality of $\pi$ is used only in Hata's lemma.
\end{remark}

\subsection{A primitive pair with rational data}\label{ss:rec}
By \eqref{eq:fWY},
\begin{equation}\label{eq:hWY}
h_n=\frac{(W-3)^{5n}(W-2)^{3n}}{Y^{4n}}\cdot\frac1z ,
\end{equation}
so $h_n$ is a combination, with integer coefficients, of the functions
\[
\XX jm:=\frac{W^j}{Y^{2m}}\cdot\frac1z\qquad(j\ge0,\ m\ge0)
\]
with $m=2n$ and $0\le j\le8n$. These can be integrated by a recursion on $m$.

\begin{lemma}\label{lem:rec}
For all integers $j\ge0$ and $m\ge0$, on $\Dom$,
\begin{equation}\label{eq:rec}
\frac{\dd}{\dd z}\Bigl(\frac{W^{j+1}}{Y^{2m+1}}\Bigr)=(j-2m)\,\XX jm-8(2m+1)\,\XX j{m+1},
\end{equation}
and
\begin{equation}\label{eq:rec0}
\XX j0=\frac{W^j}{z}=\sum_{k=0}^{j}\binom jk2^k\,z^{\,j-2k-1}.
\end{equation}
\end{lemma}
\begin{proof}
By the product rule and Lemma~\ref{lem:WY}(b), the derivative of $W^{j+1}Y^{-2m-1}$ is
$(j+1)W^jY^{-2m}/z-(2m+1)W^{j+2}Y^{-2m-2}/z$. Now write $W^{j+2}=W^j(Y^2+8)$ (Lemma~\ref{lem:WY}(a)). The second formula is the binomial theorem
for $(z+2/z)^j$.
\end{proof}

Read \eqref{eq:rec} as a formula for $\XX j{m+1}$:
\begin{equation}\label{eq:recsolved}
\XX j{m+1}=\frac{1}{8(2m+1)}\Bigl[(j-2m)\,\XX jm-\frac{\dd}{\dd z}\Bigl(\frac{W^{j+1}}{Y^{2m+1}}\Bigr)\Bigr].
\end{equation}
It expresses $\XX j{m+1}$ through $\XX jm$ and an exact derivative, so after $m$ steps everything is reduced to $\XX j0$, which by
\eqref{eq:rec0} is a sum of powers of $z$. Each power $z^{l-1}$ with $l\ne0$ is the derivative of $z^l/l$. The only term that is not a
derivative is the one with $l=0$, a multiple of $1/z$.

So every function $\XX jm$, and hence every combination $P(W)\,Y^{-2m}/z$ of such functions, is a multiple of $1/z$ plus an exact
derivative. This is the shape $c/z+G'$ that a primitive pair requires. For example, take the function $\XX 21$. We have
$\XX 20=W^2/z=z+4/z+4/z^3$, and the step \eqref{eq:recsolved} from $m=0$ to $m=1$ gives
\[
\XX 21=\frac{W^2}{Y^2}\cdot\frac1z=\frac18\Bigl[2\,\XX 20-\frac{\dd}{\dd z}\Bigl(\frac{W^3}{Y}\Bigr)\Bigr]
=\frac1z+\frac{\dd}{\dd z}\Bigl(\frac{z^2}{8}-\frac{1}{2z^2}-\frac{W^3}{8Y}\Bigr).
\]
Here $c=1$, and we divided only by $8$ (in the step of the recursion) and by $\pm2$ (in integrating the powers $z$ and $z^{-3}$).

In general two things have to be recorded. The first is the multiple of $1/z$, which becomes the number $c$. In $\XX j0=W^j/z$ it is the
middle term $\binom j{j/2}2^{j/2}/z$ of the binomial expansion \eqref{eq:rec0}, which is present only for even $j$, and the step
\eqref{eq:recsolved} from $m$ to $m+1$ multiplies it by $(j-2m)/(8(2m+1))$. The second is the list of numbers by which one divides,
because it decides in which ring the coefficients of $G$ lie. These are the integers $l=j-2k\ne0$ with $0<|l|\le j$ (the power $z^{l-1}$ of \eqref{eq:rec0} has the primitive $z^l/l$), and the numbers
$8\cdot1,\ 8\cdot3,\ \dots,\ 8(2m-1)$. The following proposition records both. It is stated for a general ring $R$ of coefficients
because it will be used twice, with $R=\Q$ in Corollary~\ref{cor:rational} and, in Section~\ref{ss:removable}, with the ring $R=R_p$ of
the $p$-integral numbers (Definition~\ref{def:pint}), that is, of the numbers without a given prime $p$ in the denominator.

\begin{proposition}\label{prop:rec}
Let $R$ be a subring of $\C$ containing $\frac12$, let $m\ge0$, and let $P\in R[W]$ be a polynomial of degree at most $d$. Assume that
$1/l\in R$ for all integers $1\le l\le d$ and $1/(2k-1)\in R$ for all integers $1\le k\le m$. Then there are a number $c\in R$ and a function
$G$ on $\Dom$ such that
\[
\frac{P(W)}{Y^{2m}}\cdot\frac1z=\frac cz+G'(z)\qquad\text{on }\Dom,
\]
and $G$ is a linear combination, with coefficients in $R$, of the functions $z^l$ $(l\in\Z,\ 0<|l|\le d)$ and
$W^{j+1}Y^{-(2m'+1)}$ $(0\le j\le d,\ 0\le m'<m)$. Explicitly, if $P=\sum_jp_jW^j$, then
\begin{equation}\label{eq:cformula}
c=\sum_{j\ \mathrm{even}}p_j\,\xi_{j,m}\binom j{j/2}2^{j/2},\qquad \xi_{j,m}:=\prod_{i=0}^{m-1}\frac{j-2i}{8(2i+1)} .
\end{equation}
\end{proposition}
\begin{proof}
By linearity it suffices to treat $P=W^j$, by induction on $m$. For $m=0$ use \eqref{eq:rec0}. The exponent $l=j-2k$ satisfies $|l|\le j$,
and the term with $l=0$ (which occurs for even $j$, $k=j/2$) is $\binom j{j/2}2^{j/2}/z$. The step from $m$ to $m+1$ is
\eqref{eq:recsolved}, where one divides by $8(2m+1)$, and it multiplies the constant $c$ by $(j-2m)/(8(2m+1))$.
\end{proof}

In \eqref{eq:cformula} one of the $m$ factors of $\xi_{j,m}$ is $0$ for even $j<2m$, so only the powers $W^j$ with
$j\ge2m$ contribute to $c$. This agrees with the fact that for $j<2m$ the function
$W^jY^{-2m}/z=(z^2+2)^j\,z^{2m-j-1}/(z^2-2)^{2m}$ has no pole at $0$.

\begin{corollary}\label{cor:rational}
For every $n\ge1$ there is a primitive pair $(G,c)$ with $c\in\Q$ and $G(1+i)-G(1-i)\in i\,\Q$. Consequently the numbers $u_n$ and $v_n$ of
\eqref{eq:uv} are rational, and $J_n=u_n+v_n\pi$. This proves the first assertion of Theorem~\textup{\ref{thm:forms}(a)}.
\end{corollary}
\begin{proof}
Apply Proposition~\ref{prop:rec} with $R=\Q$, $P=(W-3)^{5n}(W-2)^{3n}$ and $m=2n$. By \eqref{eq:hWY} the result is a primitive pair. The function $G$ is a rational
combination of powers $z^l$ and of functions $W^{j+1}Y^{-(2m'+1)}$. For these, by Lemma~\ref{lem:WY}(d),
\begin{align*}
(1+i)^l-(1-i)^l&=2i\,\im\,(1+i)^l\in i\,\Q,\\
\frac{W^{j+1}}{Y^{2m'+1}}\Big|_{z=1+i}-\frac{W^{j+1}}{Y^{2m'+1}}\Big|_{z=1-i}&=\frac{2^{j+1}}{(2i)^{2m'+1}}-\frac{2^{j+1}}{(-2i)^{2m'+1}}=-i\,(-1)^{m'}2^{\,j+1-2m'}\in i\,\Q
\end{align*}
(for negative $l$ note $(1\pm i)^{-1}=(1\mp i)/2$). The reason behind both formulas is that $G$ is a rational function of $z$ with
rational coefficients. So $G(1-i)$ is the complex conjugate of $G(1+i)$, and the difference of a number of $\Q(i)$ and its complex
conjugate is $2i$ times its imaginary part. Hence $u_n=-i\bigl(G(1+i)-G(1-i)\bigr)$ is rational. The rest is Lemma~\ref{lem:ftc}.
\end{proof}

\section{The denominators of \texorpdfstring{$u_n$}{un} and \texorpdfstring{$v_n$}{vn}}\label{sec:arith}

In this section we prove Theorem~\ref{thm:forms}(b). The strategy is to treat one prime $p$ at a time. For each $p$ we find a primitive pair
$(G,c)$ in which $c$ and $G(1+i)-G(1-i)$ have no $p$ in the denominator, up to an explicit power of $p$. Since all primitive pairs give the
same $u_n$ and $c_n$ (Lemma~\ref{lem:unique}), this bounds the power of $p$ in the denominators of $u_n$ and $c_n$.

For an odd prime $p$, the denominators arise only from integration. The integrand $h_n$ is a quotient of polynomials with integer
coefficients, and the coefficients of its partial fractions have no $p$ in their denominators (Lemma~\ref{lem:hpf}). But a primitive of
$z^j$ is $z^{j+1}/(j+1)$ and a primitive of $(z-a)^{-k}$ with $k\ge2$ is $-(z-a)^{-(k-1)}/{(k-1)}$, so integrating term by term divides the
coefficients by integers. For $h_n$ these integers are at most $4n$ (Section~\ref{ss:odd}), which is why $\lcm(1,\dots,4n)$ appears as a
common denominator. The prime $2$ is different. It occurs in the denominators of the partial fractions already, and
Section~\ref{ss:two} shows by another route that the power of $2$ contained in $\lcm(1,\dots,4n)$ suffices.

What has to be shown is the following. Let $p$ be a prime and let $p^{e}$ be the largest power of $p$ not exceeding $4n$ (with $e=0$ if $p>4n$), so that $p^e$ is
the exact power of $p$ dividing $\lcm(1,\dots,4n)$. The integer $D_n=\lcm(1,\dots,4n)/\Phi_n$ is divisible by exactly $p^e$ if $p$ is not
removable, and it is not divisible by $p$ if $p$ is removable. So Theorem~\ref{thm:forms}(b) follows from four statements.
\begin{itemize}
\item For every odd prime $p$, the number $p^eu_n$ has no $p$ in its denominator (Proposition~\ref{prop:odd}). Here it suffices to
integrate the partial fractions of $h_n$ term by term.
\item For a removable prime $p$, the number $u_n$ itself has no $p$ in its denominator (Proposition~\ref{prop:removable}). Here we show
that the coefficients of $h_n$ which integration divides by a multiple of $p$ are themselves divisible by $p$ (they are $p$ times a
number without $p$ in its denominator).
\item For $p=2$, the number $2^eu_n$ has no $2$ in its denominator (Proposition~\ref{prop:two}, where $e$ is called $L$). This needs a different coordinate and
takes the most work.
\item The number $c_n=2v_n$ is an integer (Propositions~\ref{prop:odd} and~\ref{prop:two}). Since $D_n$ is even, $D_nv_n$ is then an
integer as well.
\end{itemize}
Section~\ref{ss:integrality} puts the four statements together.

\subsection{\texorpdfstring{$p$}{p}-integers and partial fractions}\label{ss:pint}
The primitive pairs involve the irrational numbers $\rt$ and $i$, so we need a notion of ``no $p$ in the denominator'' for such numbers.
Recall that an \emph{algebraic integer} is a complex number that is a root of a monic polynomial with integer coefficients. The algebraic
integers form a ring, and a rational number is an algebraic integer only if it is an integer \cite[Chapter~2]{Marcus}.

\begin{definition}\label{def:pint}
Let $p$ be a prime. A complex number $x$ is \emph{$p$-integral} if $Mx$ is an algebraic integer for some positive integer $M$ not divisible
by $p$. The set of $p$-integral numbers is denoted by $R_p$, and $\Zp=\{a/b:\ a,b\in\Z,\ p\nmid b\}$.
\end{definition}

For a rational number this is the usual notion: it is $p$-integral if and only if $p$ does not divide its denominator
(Lemma~\ref{lem:pint}(b) below). For instance $\frac34$ is $p$-integral for every odd prime $p$, but not for $p=2$. The definition
extends the notion to the irrational numbers that occur in our formulas. The numbers $\rt$ and $i$ are $p$-integral for every $p$, while
$1/\rt=\rt/2$ and $1/(1+i)=(1-i)/2$ are $p$-integral for every odd $p$ but not for $p=2$ (Appendix~\ref{app:pf}).

\begin{lemma}\label{lem:pint}
\textup{(a)} $R_p$ is a subring of $\C$. It contains $\rt$, $i$, and $1/M$ for every integer $M$ not divisible by $p$.

\textup{(b)} $R_p\cap\Q=\Zp$, that is, a rational number is $p$-integral if and only if $p$ does not divide its denominator.

\textup{(c)} A rational number that is $p$-integral for every prime $p$ is an integer.
\end{lemma}
The proof takes three lines (Appendix~\ref{app:pf}).

All functions that occur are rational functions with poles only at $0$ and $\pm\rt$. For a subring $R$ of $\C$ containing $\rt$ let $\cA(R)$
be the set of functions on $\Dom$ of the form
\[
\varphi(z)=\frac{M(z)}{z^{e_0}(z-\rt)^{e_+}(z+\rt)^{e_-}}\qquad\text{with }M\in R[z]\text{ and integers }e_0,e_+,e_-\ge0 .
\]
We also write $e_a$ for the exponent that belongs to the pole $a\in\{0,\rt,-\rt\}$, so $e_{\rt}=e_+$ and $e_{-\rt}=e_-$. The set $\cA(R)$ is
closed under addition and multiplication (bring the fractions to a common denominator) and under differentiation (quotient rule), and it
contains $W=(z^2+2)/z$ and $1/Y=z/((z-\rt)(z+\rt))$.

A function of $\cA(R)$ is integrated by decomposing it into partial fractions, that is, into powers $z^j$ and negative powers
$(z-a)^{-k}$ at the three poles, and integrating each term. The next lemma says two things. First, if $R$ contains $\rt$ and $\frac12$, the decomposition does not
introduce new denominators: its coefficients lie again in $R$. Second, every term has an explicit primitive except the terms
$(z-a)^{-1}$, whose primitives would be logarithms. A primitive pair allows such a term at $a=0$ (it is the term $c/z$) but not at
$a=\pm\rt$.

\begin{lemma}[partial fractions with coefficients in $R$]\label{lem:pf}
Let $R$ be a subring of $\C$ containing $\rt$ and $\frac12$. Every $\varphi\in\cA(R)$ can be written in exactly one way as
\begin{equation}\label{eq:pf}
\varphi(z)=\sum_{j=0}^{d}t_j\,z^j+\sum_{a\in\{0,\rt,-\rt\}}\ \sum_{k\ge1}\frac{q_{a,k}}{(z-a)^k}\qquad(\text{finite sums}),
\end{equation}
and all coefficients $t_j$, $q_{a,k}$ lie in $R$. The terms $(z-a)^{-k}$ occur only for $k\le e_a$. In particular there are no such terms
if $e_a=0$. Moreover:
\begin{enumerate}
\item[(i)] the decomposition of $\varphi'$ is obtained by differentiating \eqref{eq:pf} term by term. In particular $\varphi'$ has no terms
$(z-a)^{-1}$.
\item[(ii)] if $q_{\rt,1}=q_{-\rt,1}=0$, then $\varphi(z)-q_{0,1}/z=G_\varphi'(z)$ on $\Dom$, where
\begin{equation}\label{eq:Gphi}
G_\varphi(z):=\sum_{j=0}^{d}\frac{t_j}{j+1}\,z^{j+1}-\sum_{a}\sum_{k\ge2}\frac{q_{a,k}}{k-1}\cdot\frac1{(z-a)^{k-1}} .
\end{equation}
\end{enumerate}
\end{lemma}
The point of the hypothesis is that the distances $\rt$ and $2\rt$ between the poles are invertible in $R$, which is true in $R_p$ for odd $p$
and false for $p=2$. The proof is in Appendix~\ref{app:pf}. The numbers by which one divides in \eqref{eq:Gphi}, namely $j+1$ and $k-1$, will be
called the \emph{integration indices}. They are the only source of new denominators: the coefficients $t_j$ and $q_{a,k}$ lie in $R$,
and the coefficients of $G_\varphi$ are these numbers divided by integration indices.

\subsection{Odd primes: the integration indices are at most \texorpdfstring{$4n$}{4n}}\label{ss:odd}

We apply Lemma~\ref{lem:pf} to the integrand $h_n$ itself. The degree of its polynomial part and the orders of its poles determine
which integration indices occur.

\begin{lemma}\label{lem:hpf}
Let $n\ge1$. The decomposition \eqref{eq:pf} of $h_n$ has the form
\begin{equation}\label{eq:hpf}
h_n(z)=\sum_{j=0}^{4n-1}t_jz^j+\sum_{k=1}^{4n+1}\frac{a_k}{z^k}+\sum_{k=1}^{4n}\Bigl(\frac{b_k}{(z-\rt)^k}+\frac{b_k'}{(z+\rt)^k}\Bigr)
\end{equation}
with coefficients $t_j,a_k,b_k,b_k'$ in the ring $\Z[\rt,\frac12]$. In particular all coefficients are $p$-integral for every odd prime $p$.
Moreover $b_1=b_1'=0$ and $a_1=c_n$. Hence $(G_{h_n},c_n)$, with $G_{h_n}$ defined by \eqref{eq:Gphi}, is a primitive pair.
\end{lemma}
\begin{proof}
$h_n=N^n/\bigl(z^{4n+1}(z-\rt)^{4n}(z+\rt)^{4n}\bigr)$ has a numerator with integer coefficients, so Lemma~\ref{lem:pf} applies with
$R=\Z[\rt,\frac12]$, and the polynomial part has degree $16n-(12n+1)=4n-1$ (the degree of the numerator minus that of the denominator).
Let $(G_0,c_0)$ be the primitive pair of Corollary~\ref{cor:rational}, so that $h_n=c_0/z+G_0'$. Its function $G_0$ is a rational combination of
powers of $z$ and of the functions $W^{j+1}Y^{-(2m'+1)}$ (Proposition~\ref{prop:rec}), hence $G_0\in\cA(\Q(\rt))$. By
Lemma~\ref{lem:pf}(i) the function $G_0'$ has no terms $(z-a)^{-1}$, so the uniqueness of the decomposition gives $b_1=b_1'=0$ and $a_1=c_0$, and
$c_0=c_n$ by definition. The last assertion is Lemma~\ref{lem:pf}(ii).
\end{proof}

The statements $b_1=b_1'=0$ and $a_1=c_n$ say that $h_n$ has no term $(z\mp\rt)^{-1}$, which would produce a logarithm, and that its
coefficient of $1/z$ is the number $c_n$ of \eqref{eq:uv}. They were not obtained by computing these coefficients, but by comparing
with the primitive pair of Corollary~\ref{cor:rational}, which we already have.

\begin{proposition}[odd primes]\label{prop:odd}
Let $p$ be an odd prime, $n\ge1$, and let $e=\fl{\log_p(4n)}$ be the largest integer with $p^e\le4n$. Then
\[
c_n\in\Zp\qquad\text{and}\qquad p^{e}\,u_n\in\Zp .
\]
\end{proposition}
\begin{proof}
Use the primitive pair $(G_{h_n},c_n)$ of Lemma~\ref{lem:hpf}. First, $c_n=a_1$ is rational and $p$-integral, so $c_n\in\Zp$ by
Lemma~\ref{lem:pint}(b). Next, in \eqref{eq:Gphi} every coefficient is a $p$-integral number divided by an integration index $m$ with
$1\le m\le4n$. Indeed $j+1\le4n$, and $k-1\le4n$ for the pole $0$ and $k-1\le4n-1$ for the poles $\pm\rt$. For such $m$ the number $p^e/m$ lies
in $\Zp$, because the power of $p$ dividing $m$ is at most $m\le4n$ and hence at most $p^e$. So $p^eG_{h_n}$ is a combination, with $p$-integral coefficients, of powers of $z$, $1/z$, $1/(z-\rt)$ and
$1/(z+\rt)$. At $z=1\pm i$ these four quantities are $p$-integral, because
\begin{equation}\label{eq:units}
\begin{gathered}
\frac1{1\pm i}=\frac{1\mp i}{2},\\
\frac{1}{1\pm i-\rt}=-\frac{(1\pm i+\rt)(1\pm i)}{4},\qquad \frac{1}{1\pm i+\rt}=-\frac{(1\pm i-\rt)(1\pm i)}{4}
\end{gathered}
\end{equation}
(Appendix~\ref{app:pf}); the right-hand sides have no denominators other than $2$ and $4$. Therefore $p^e\bigl(G_{h_n}(1+i)-G_{h_n}(1-i)\bigr)\in R_p$, and $p^eu_n=-i\,p^e\bigl(G_{h_n}(1+i)-G_{h_n}(1-i)\bigr)$ is a rational
number in $R_p$, hence in $\Zp$.
\end{proof}

In other words, the power of an odd prime $p$ in the denominator of $u_n$ is at most $p^e$, the power of $p$ contained in
$\lcm(1,\dots,4n)$. The only loss is the division by the integration indices.

\subsection{Removable primes}\label{ss:removable}
By Proposition~\ref{prop:odd} a prime $p$ with $\sqrt{4n}<p\le4n$ occurs in the denominator of $u_n$ at most to the first power.
For the removable primes of Definition~\ref{def:removable} it does not occur at all.

\begin{proposition}[removable primes]\label{prop:removable}
If the prime $p$ is removable for $n$, then $u_n\in\Zp$.
\end{proposition}

The idea of the proof is as follows. In the primitive $G_{h_n}$ of Lemma~\ref{lem:hpf} the prime $p$, which is odd, can enter the denominators only
through the integration indices that are multiples of $p$. We show that the corresponding coefficients of $h_n$ are divisible by $p$
(that is, they are $p$ times a $p$-integral number), so that nothing is lost. The reason is that $h_n$ is a $p$-th power $g^p$ times a factor $\omega$ with small exponents, and that the
derivative of a $p$-th power is divisible by $p$. Namely, if $\omega=G_{\mathrm{red}}'$, then by the product rule
\[
h_n=g^p\,G_{\mathrm{red}}'=\bigl(g^pG_{\mathrm{red}}\bigr)'-p\,g^{p-1}g'\,G_{\mathrm{red}} .
\]
Now compare the partial fractions of the two sides. The partial fractions of $g^pG_{\mathrm{red}}$ and of $g^{p-1}g'G_{\mathrm{red}}$
have $p$-integral coefficients. For a term of $h_n$ with integration index $m$, the differentiation in the first term on the right
supplies the factor $\pm m$, and the second term carries the factor $p$. So whenever $p$ divides $m$, the coefficient of $h_n$ is $p$
times a $p$-integral number. The three conditions of
Definition~\ref{def:removable} are what this argument needs. The parity of $\fl{4n/p}$ and the inequality
$(5n\bmod p)+(3n\bmod p)<(4n\bmod p)$ guarantee that $\omega$ has a primitive $G_{\mathrm{red}}$ which is a rational function with
$p$-integral coefficients (Step~2 below), and $p^2>4n$ guarantees that no integration index of $h_n$ is divisible by $p^2$, because
these indices are at most $4n$, so that one factor $p$ is enough (Step~4).

For example, let $n=2$ and $p=3$. Then $5n=10=3\cdot3+1$, $3n=6=3\cdot2+0$ and $4n=8=3\cdot2+2$, and the three conditions hold
($9>8$, the quotient $2$ of $8$ by $3$ is even, and $1+0<2$). Here $h_2=g^3\,\omega$ with $g=(W-3)^3(W-2)^2/Y^2$ and
\[
\omega=\frac{W-3}{Y^2}\cdot\frac1z=\frac{(z-1)(z-2)}{(z^2-2)^2}=\frac{\dd}{\dd z}\Bigl(\frac{3W-8}{8Y}\Bigr).
\]
The primitive $G_{\mathrm{red}}=(3W-8)/(8Y)$ has only the denominator $8$. In the decomposition \eqref{eq:hpf} of $h_2$ the integration
indices divisible by $3$ are $3$ and $6$. They belong, for instance, to the terms $-998316\,z^2$ and $873\,z^5$ of the polynomial part.
Both coefficients are divisible by $3$, and integrating these terms gives $-332772\,z^3$ and $\frac{873}{6}\,z^6=\frac{291}{2}\,z^6$,
with no $3$ in the denominators. Taken separately, the two terms $(g^3G_{\mathrm{red}})'$ and $3\,g^2g'G_{\mathrm{red}}$ of the product
rule both contain $\frac{27}{8}\,z^{8}$, whose integration index $9=3^2$ exceeds $4n=8$. They cancel in $h_2$, whose polynomial part has degree $7$. This is
why the proof compares the coefficients of $h_n$ and does not integrate the second term by itself.

\begin{proof}
A removable prime is odd, since $p^2>4n\ge4$. Divide $5n$, $3n$ and $4n$ by $p$ with remainder:
\[
5n=p\,q_1+\alpha_1,\qquad 3n=p\,q_2+\alpha_2,\qquad 4n=p\,q_0+\delta,\qquad 0\le\alpha_1,\alpha_2,\delta\le p-1 .
\]
By definition, $q_0=\fl{4n/p}$ is even and $\alpha_1+\alpha_2<\delta$. As $4n$ and $pq_0$ are even, $\delta$ is even, and $\delta\ge2$
because $\delta>\alpha_1+\alpha_2\ge0$.

\smallskip\emph{Step 1: a $p$-th power splits off.} By \eqref{eq:hWY} the factors $W-3$, $W-2$ and $Y$ occur in $h_n$ with the exponents
$5n$, $3n$ and $-4n$. Writing each of $5n$, $3n$ and $4n$ as a multiple of $p$ plus its remainder gives
\begin{equation}\label{eq:gomega}
h_n=g^p\cdot\omega\qquad\text{with}\qquad g:=\frac{(W-3)^{q_1}(W-2)^{q_2}}{Y^{q_0}},\qquad \omega:=\frac{(W-3)^{\alpha_1}(W-2)^{\alpha_2}}{Y^{\delta}}\cdot\frac1z .
\end{equation}

\emph{Step 2: the small factor $\omega$ is the derivative of a $p$-integral function.} Apply Proposition~\ref{prop:rec} with $R=R_p$,
$P=(W-3)^{\alpha_1}(W-2)^{\alpha_2}$, $d=\alpha_1+\alpha_2$ and $m=\delta/2$. The numbers that must be invertible in $R_p$ are $2$, the integers
$l\le\alpha_1+\alpha_2<\delta<p$ and the odd integers $2k-1\le\delta-1<p$, and none of them is divisible by $p$. So $\omega=c_\omega/z+G_{\mathrm{red}}'$ with
$c_\omega\in R_p$ and $G_{\mathrm{red}}\in\cA(R_p)$ (by Proposition~\ref{prop:rec}, $G_{\mathrm{red}}$ is a combination with coefficients in
$R_p$ of powers of $z$ and of functions $W^{j+1}Y^{-(2m'+1)}$, and these lie in $\cA(R_p)$). Moreover $c_\omega=0$, for the following
reason. In terms of $z$,
\[
\omega=\frac{\bigl((z-1)(z-2)\bigr)^{\alpha_1}(z^2-2z+2)^{\alpha_2}\;z^{\,\delta-\alpha_1-\alpha_2-1}}{(z^2-2)^{\delta}}
\]
has no power of $z$ in its denominator, because $\delta-\alpha_1-\alpha_2-1\ge0$. So it belongs to $\cA(R_p)$ with $e_0=0$, and its
decomposition \eqref{eq:pf} has no term $1/z$ (Lemma~\ref{lem:pf}), and
$G_{\mathrm{red}}'$ has no term $1/z$ either (Lemma~\ref{lem:pf}(i)). Hence $c_\omega/z=\omega-G_{\mathrm{red}}'$ has no term $1/z$, which by the
uniqueness of the decomposition means $c_\omega=0$. Thus $\omega=G_{\mathrm{red}}'$. Here the inequality $\alpha_1+\alpha_2<\delta$ is used once more. The parity of $q_0$ was used to make $\delta$ even, so
that Proposition~\ref{prop:rec} applies.

\smallskip\emph{Step 3: integration by parts.} Put $H:=g^pG_{\mathrm{red}}$ and $\eta:=g^{p-1}g'G_{\mathrm{red}}$. Both belong to $\cA(R_p)$, and by the
product rule $H'=p\,g^{p-1}g'G_{\mathrm{red}}+g^pG_{\mathrm{red}}'$, that is
\begin{equation}\label{eq:cartier}
h_n=H'-p\,\eta .
\end{equation}

\emph{Step 4: compare coefficients.} It remains to read off from \eqref{eq:cartier} the coefficients of the primitive $G_{h_n}$.
Let $t^H_j,q^H_{a,k}$ and $t^\eta_j,q^\eta_{a,k}$ be the coefficients of $H$ and $\eta$ in \eqref{eq:pf}. They
are $p$-integral by Lemma~\ref{lem:pf}. By \eqref{eq:cartier}, Lemma~\ref{lem:pf}(i) and the uniqueness of the decomposition, the coefficients
$t_j$, $q_{a,k}$ of $h_n$ (in \eqref{eq:hpf}: $q_{0,k}=a_k$, $q_{\rt,k}=b_k$, $q_{-\rt,k}=b_k'$) are
\[
t_j=(j+1)\,t^H_{j+1}-p\,t^\eta_j,\qquad q_{a,k}=-(k-1)\,q^H_{a,k-1}-p\,q^\eta_{a,k}\quad(k\ge2).
\]
(The terms with $k=1$ do not enter $G_{h_n}$: by Lemma~\ref{lem:hpf} they are $b_1=b_1'=0$ and $a_1/z=c_n/z$.)
Dividing by the integration indices, the coefficients of the primitive $G_{h_n}$ in \eqref{eq:Gphi} are
\begin{equation}\label{eq:pcancel}
\frac{t_j}{j+1}=t^H_{j+1}-\frac{p}{j+1}\,t^\eta_j,\qquad -\frac{q_{a,k}}{k-1}=q^H_{a,k-1}+\frac{p}{k-1}\,q^\eta_{a,k}.
\end{equation}
The first term on each right-hand side is $p$-integral as it stands. The second term contains the quotient $p/(j+1)$ or $p/(k-1)$, in
which the factor $p$ of \eqref{eq:cartier} meets the integration index.
The indices $j+1$ and $k-1$ are at most $4n<p^2$, so they are divisible by $p$ at most once and $p/(j+1)$, $p/(k-1)$ lie in $\Zp$. Hence
\emph{every} coefficient of $G_{h_n}$ is $p$-integral, and as in the proof of Proposition~\ref{prop:odd} we conclude
$u_n=-i\bigl(G_{h_n}(1+i)-G_{h_n}(1-i)\bigr)\in R_p\cap\Q=\Zp$.
\end{proof}

\subsection{The prime \texorpdfstring{$2$}{2}}\label{ss:two}

\begin{proposition}[the prime $2$]\label{prop:two}
Let $n\ge1$ and $L=\fl{\log_2(4n)}$. Then
\[
c_n\in\Z_{(2)}\qquad\text{and}\qquad 2^{L}\,u_n\in\Z_{(2)} .
\]
\end{proposition}

So the prime $2$ costs no more than what $\lcm(1,\dots,4n)$ already contains. This cannot be proved with the decomposition \eqref{eq:hpf}. The
distances $\rt$ and $2\rt$ between the poles are not invertible in $R_2$, and the individual coefficients have large powers of $2$ in their
denominators (for $n=2$ the denominator of $b_2$ is $1024$, whereas that of $u_2$ is $14$). The powers of $2$ cancel in the sum. To make
the cancellation visible we change the coordinate.

The proof has three parts. First, a M\"obius transformation $z\mapsto w$ brings the integrand into a form in which the symmetry
$z\mapsto2/z$ becomes $w\mapsto-w$ and the end points of the path become $w=\pm1$ (Lemma~\ref{lem:w}). Second, the new integrand is written as
a combination of two families of simple functions, the ``atoms'', whose integrals can be computed explicitly (Lemma~\ref{lem:atoms}).
Third, we count the powers of $2$ in the results and add up (Table~\ref{tab:levels}). No single step is difficult, but there is
some bookkeeping.

\subsubsection*{The coordinate $w$}
Let
\[
\eps:=\rt-1,\qquad \vk:=1+\eps^2=2\rt\,\eps .
\]
Note that $\eps(\rt+1)=1$, so $\eps$ and $1/\eps=\rt+1$ both lie in $\Z[\rt]$ (in particular they are algebraic integers), and that
$\eps^2+2\eps=1$. Put
\begin{equation}\label{eq:w}
w=w(z):=i\eps\,\frac{\rt+z}{\rt-z},\qquad\text{equivalently}\qquad z=\rt\,\frac{w-i\eps}{w+i\eps}.
\end{equation}
This M\"obius map sends $z=-\rt,\ \rt,\ 0,\ \infty$ to $w=0,\ \infty,\ i\eps,\ -i\eps$, the circle $|z|=\rt$ to the real axis, and the end points
$1-i$, $1+i$ to $w=1$, $-1$. One checks that $w(2/z)=-w(z)$, so the symmetry $z\mapsto2/z$ becomes $w\mapsto-w$.

\begin{lemma}\label{lem:w}
For $z\in\Dom$ and $w=w(z)$ one has $w\neq0$, $w^2+\eps^2\neq0$, and
\begin{enumerate}
\item[(a)] $\displaystyle f(z)=\frac{1}{128\,\eps}\cdot\frac{(1+\eps^2w^2)^5\,(1-w^2)^3}{w^4\,(w^2+\eps^2)^4}$.
\item[(b)] $\displaystyle\frac1z=\frac{2i\eps}{w^2+\eps^2}\cdot\frac{\dd w}{\dd z}$.
\item[(c)] $w(1+i)=-1$ and $w(1-i)=1$.
\end{enumerate}
Consequently
\begin{equation}\label{eq:kn}
h_n(z)=\frac{2i\eps}{(128\,\eps)^{n}}\;k_n\bigl(w(z)\bigr)\,\frac{\dd w}{\dd z},\qquad
k_n(w):=\frac{(1+\eps^2w^2)^{5n}\,(1-w^2)^{3n}}{w^{4n}\,(w^2+\eps^2)^{4n+1}} .
\end{equation}
\end{lemma}
The proof is a direct computation with the four factors of $f$: the factors $(z-1)(z-2)$, $z^2-2z+2$, $z^2-2$ and $z$ are constant
multiples of $1+\eps^2w^2$, $1-w^2$, $w$ and $w^2+\eps^2$, each multiplied by $(\rt-z)^2$, and the powers of $\rt-z$ cancel in $f$.
The gain is that $k_n$ is an \emph{even} function of $w$ with
poles only at $0$ and $\pm i\eps$, to be ``integrated from $-1$ to $1$''.

\subsubsection*{Data of a function}
We do the bookkeeping of primitives directly in the variable $w$.

\begin{definition}\label{def:data}
Let $\varphi$ be a function on $\C\setminus\{0,i\eps,-i\eps\}$. We say that $\varphi$ \emph{has data} $(U,V)\in\C^2$ if there is a holomorphic
function $F$ on $\C\setminus\{0,i\eps,-i\eps\}$ such that
\[
F'(w)=\varphi(w)-\frac{V}{w^2+\eps^2}\qquad\text{and}\qquad F(1)-F(-1)=U .
\]
\end{definition}
This is the notion of a primitive pair, in the variable $w$. By Lemma~\ref{lem:w}(b) the function $1/(w^2+\eps^2)$ plays the role
of $1/z$, so, up to constant factors, $V$ corresponds to the number $c$ and $U$ to the difference $G(1+i)-G(1-i)$
(Lemma~\ref{lem:transfer} below gives the exact relations).
Informally, ``$\int_{-1}^{1}\varphi(w)\,\dd w=U+V\vartheta$, where $\vartheta=\int_{-1}^1\dd w/(w^2+\eps^2)$'' (both integrals along the same path
from $-1$ to $1$, which avoids the points $0$ and $\pm i\eps$, where $\varphi$ may have poles), but the logarithmic integral $\vartheta$ is never evaluated. Data
are additive. If $\varphi_1,\varphi_2$ have data $(U_1,V_1)$, $(U_2,V_2)$, then $\varphi_1+\varphi_2$ has data $(U_1+U_2,\ V_1+V_2)$, and
$b\,\varphi_1$ has data $(bU_1,\,bV_1)$ for every complex number $b$.

\begin{lemma}[transfer to $z$]\label{lem:transfer}
If $k_n$ has data $(U,V)$, then
\begin{equation}\label{eq:transfer}
u_n=-\frac{2\eps}{(128\,\eps)^{n}}\,U,\qquad c_n=\frac{1}{(128\,\eps)^{n}}\,V .
\end{equation}
\end{lemma}
\begin{proof}
With $F$ as in Definition~\ref{def:data} put $G(z):=2i\eps(128\eps)^{-n}F(w(z))$ and $c:=(128\eps)^{-n}V$. By the chain rule,
\eqref{eq:kn} and Lemma~\ref{lem:w}(b),
\[
G'(z)=\frac{2i\eps}{(128\eps)^n}\Bigl(k_n(w)-\frac{V}{w^2+\eps^2}\Bigr)\frac{\dd w}{\dd z}=h_n(z)-\frac cz ,
\]
so $(G,c)$ is a primitive pair. By Lemma~\ref{lem:w}(c), $G(1+i)-G(1-i)=2i\eps(128\eps)^{-n}\bigl(F(-1)-F(1)\bigr)$. Multiply by $-i$ and use
Lemma~\ref{lem:unique}.
\end{proof}

Lemma~\ref{lem:transfer} shows what has to be proved. We have $(128\,\eps)^{-n}=2^{-7n}\eps^{-n}$, and the factor $\eps^{-n}$ is
harmless, since $\eps$ and $1/\eps$ are algebraic integers. So $c_n$ is $2$-integral exactly when $V/2^{7n}$ is
$2$-integral, and $2^Lu_n$ is $2$-integral exactly when $U/2^{7n-1-L}$ is $2$-integral. The rest of this subsection shows that the data
of $k_n$ are divisible by these powers of $2$.

\subsubsection*{Measuring powers of $2$}
The data will be numbers of the field $\Q(\rt)$, in which $2=(\rt)^2$. We measure divisibility by powers of $\rt$. This is finer than
counting powers of $2$. For instance, the number $\vk=(\rt)^3\eps$ contains an odd power of $\rt$.

\begin{definition}\label{def:level}
For $k\in\Z$ let $\lev k:=\{(\rt)^k\,y:\ y\in R_2\}$. If $x\in\lev k$ we say that $x$ \emph{has level at least $k$}.
\end{definition}

So $\lev0=R_2$, and level $2k$ for a rational number means ``divisible by $2^k$''. In this language the task is to show that $V$ has
level at least $14n$ and that $U$ has level at least $14n-2-2L$. The facts we need are collected in the following lemma (proofs in Appendix~\ref{app:two}). Here $s_2(a)$ denotes the
number of digits $1$ in the binary expansion of $a$.

\begin{lemma}\label{lem:levels}
\textup{(a)} $\lev j\cdot\lev k\subseteq\lev{j+k}$, \ $\lev k+\lev k\subseteq\lev k$, and $\lev k\subseteq\lev j$ for $j\le k$.

\textup{(b)} $2^k\in\lev{2k}$ for $k\in\Z$. If $m$ is an odd integer, $m$ and $1/m$ lie in $\lev0$. For $k\in\Z$, $\eps^k\in\lev0$ and
$\vk^k\in\lev{3k}$.

\textup{(c)} \textup{(Legendre)} $a!$ equals $2^{\,a-s_2(a)}$ times an odd integer. \textup{(Kummer)} $\binom{2j}{j}$ equals $2^{\,s_2(j)}$ times an odd
integer. Moreover $2s_2(a)\le a+1$ for $a\ge0$, and $s_2(a)\le\fl{\log_2a}+1$ for $a\ge1$.

\textup{(d)} If a rational number $x$ lies in $\lev{-2L}$, then $2^Lx\in\Z_{(2)}$.
\end{lemma}

Part (a) says that a product has at least the sum of the levels of its factors, and that a sum has at least the smaller of the levels
of its terms. Part (c) is where
factorials and central binomial coefficients enter. It gives the exact power of $2$ dividing $a!$, namely $2^{\,a-s_2(a)}$, which is
at least $2^a/(2a)$ for $a\ge1$, and the exact power of $2$ dividing a central binomial coefficient.

\subsubsection*{The building blocks}
Besides the power $w^{-4n}$, the function $k_n$ contains powers of the three polynomials $1-w^2$, $w^2+\eps^2$ and $1+\eps^2w^2$. All
three have degree $1$ in $w^2$ and no two of them are proportional, so each of them is a linear combination of the other two. We normalise the first two so that their sum is $1$.
Introduce
\[
x:=\frac{1-w^2}{\vk},\qquad y:=1-x=\frac{w^2+\eps^2}{\vk}.
\]
The factor $1+\eps^2w^2$ of $k_n$ is also a combination of $x$ and $y$. Since $\eps^2+2\eps=1$,
\begin{equation}\label{eq:third}
1+\eps^2w^2=\vk\,(y+2\eps x).
\end{equation}
Inserting this in \eqref{eq:kn} and expanding $(y+2\eps x)^{5n}$ by the binomial theorem,
\begin{equation}\label{eq:knexp}
k_n(w)=\vk^{\,4n-1}\sum_{i=0}^{5n}\binom{5n}{i}(2\eps)^i\;\frac{x^{\,3n+i}\,y^{\,n-1-i}}{w^{4n}} .
\end{equation}
Here the power of $\vk$ is $5n+3n-(4n+1)$ and the exponent of $y$ is $(5n-i)-(4n+1)$.
(From here to the end of the proof of Proposition~\ref{prop:two} the letter $i$ denotes a summation index, and the imaginary unit does not occur.)
For $i\le n-1$ the exponent of $y$ is non-negative and $y^{n-1-i}=(1-x)^{n-1-i}$ is a polynomial in $x$. For $i\ge n$ it is negative. Then we write
$x^{3n+i}=(1-y)^{3n+i}$ and expand in powers of $y$. Altogether we arrive at two kinds of ``atoms''.
These are the functions $x^a\,w^{-2m}$ and $y^{-k}\,w^{-2m}$ with $a\ge0$ and $k\ge1$. In the first case only atoms of the first kind
occur. In the second case both kinds occur, after the non-negative powers of $y$ have been rewritten as polynomials in $x$ through
$y=1-x$. The relation $x+y=1$ is what makes this possible: a product in which both $x$ and $y$ occur is traded for a sum of terms that
contain only a power of $x$ or only a negative power of $y$.

\begin{lemma}[the atoms]\label{lem:atoms}
Let $m\ge0$ be an integer.
\begin{enumerate}
\item[(a)] For every integer $a\ge0$ the function $x^a\,w^{-2m}$ has data $(B_{m,a},0)$ with
\[
B_{m,a}=\frac{2^{a+1}\,a!}{\vk^{a}\,\prod_{j=0}^{a}(2j+1-2m)}\ \in\ \lev{\,a+2-2s_2(a)} .
\]
\item[(b)] For every integer $k\ge1$ the function $y^{-k}\,w^{-2m}$ has data $(U,V)$ with $U,V\in\lev{4-k}$.
\end{enumerate}
\end{lemma}
For (a) one expands $(1-w^2)^a$ and integrates the even powers of $w$. The alternating sum that arises has a closed form, in which the
factorial $a!$ supplies many factors $2$. For (b) one reduces to $m=0$ by partial fractions in $w^2$ and uses the classical recursion for
$\int(w^2+\eps^2)^{-j}\dd w$, whose coefficients are central binomial coefficients. All details are in Appendix~\ref{app:two}.

To see where the level in (a) comes from, take $m=0$. Then
\[
\vk^a\,B_{0,a}=\int_{-1}^{1}(1-w^2)^a\,\dd w=\frac{2^{a+1}\,a!}{1\cdot3\cdot5\cdots(2a+1)} .
\]
By Lemma~\ref{lem:levels}(c) the numerator $2^{a+1}a!$ is $2^{\,2a+1-s_2(a)}$ times an odd integer, so it has level $4a+2-2s_2(a)$. The odd
denominator does not change the level, and the division by $\vk^a$ costs $3a$ levels. What remains is $a+2-2s_2(a)$, which is large when $a$ is large. The second entry of the data
in (a) is $0$ because $x^aw^{-2m}$ is a combination of even powers of $w$, and these have the primitives $w^{2q+1}/(2q+1)$, without a
logarithm. In (b) the function has poles at the two zeros of $w^2+\eps^2$, and a multiple $V/(w^2+\eps^2)$ has to be subtracted before a primitive exists.
For $m=0$ the count in (b) is this: the data of $(w^2+\eps^2)^{-k}$ have level at least $-4(k-1)$ by \eqref{eq:m0}, and
$y^{-k}=\vk^k(w^2+\eps^2)^{-k}$, where $\vk^k$ has level $3k$. Together this is $4-k$.

\begin{proof}[Proof of Proposition~\ref{prop:two}]
We show that $k_n$ has data $(U,V)$ with
\begin{equation}\label{eq:UVlevels}
U\in\lev{\,14n-2-2L},\qquad V\in\lev{\,14n}.
\end{equation}
Then \eqref{eq:transfer} gives the claim, since $(128\eps)^{-n}=2^{-7n}\eps^{-n}\in\lev{-14n}$ and $2\eps\in\lev2$, so $c_n\in\lev0=R_2$ and
$u_n\in\lev{-2L}$ (the levels add up to $2-14n+(14n-2-2L)=-2L$). As both are rational, $c_n\in\Z_{(2)}$ and $2^Lu_n\in\Z_{(2)}$ by Lemmas~\ref{lem:pint}(b) and~\ref{lem:levels}(d).

Since data are additive, it suffices to find data for each term of \eqref{eq:knexp} separately.
Consider the $i$-th term $T_i:=x^{3n+i}y^{n-1-i}w^{-4n}$ of \eqref{eq:knexp}. We use Lemma~\ref{lem:atoms} with $m=2n$.

\emph{Case $i\le n-1$.} Then $T_i=\sum_b\binom{n-1-i}{b}(-1)^b\,x^{3n+i+b}w^{-4n}$ is a combination with integer coefficients of atoms
$x^aw^{-4n}$ with $3n\le a\le4n-1$. For such $a$ we have $s_2(a)\le\fl{\log_2a}+1\le L+1$ (a positive integer $a<4n$ has at most $L+1$ binary
digits), so $B_{2n,a}$ has level at least
$a+2-2s_2(a)\ge3n-2L$. Hence $T_i$ has data $(U_i,0)$ with $U_i\in\lev{3n-2L}$.

\emph{Case $i\ge n$.} Put $j:=i-n+1\ge1$, so that $T_i=x^{3n+i}y^{-j}w^{-4n}=\sum_{l}\binom{3n+i}{l}(-1)^l\,y^{\,l-j}w^{-4n}$. If $l\ge j$, then
$y^{l-j}=(1-x)^{l-j}$ is a polynomial in $x$ and we get atoms $x^bw^{-4n}$ with $b\ge0$, of level at least $b+2-2s_2(b)\ge1$ (recall $2s_2(b)\le b+1$ from Lemma~\ref{lem:levels}(c)). If $l<j$ we get
atoms $y^{-k}w^{-4n}$ with $1\le k=j-l\le j$, whose data have level at least $4-k\ge4-j$. Hence $T_i$ has data $(U_i,V_i)$ with
$U_i\in\lev{\min(1,\,4-j)}$ and $V_i\in\lev{4-j}$.

Now multiply by $\binom{5n}{i}(2\eps)^i\in\lev{2i}$ and by $\vk^{4n-1}\in\lev{12n-3}$ and add up. Table~\ref{tab:levels} shows the resulting levels.
In its first row the bound is $(12n-3)+2i+(3n-2L)=15n-3+2i-2L$, which is at least $14n-2-2L$ because $n\ge1$ and $i\ge0$. The bounds of
the other two rows are evaluated in the table. So the
contributions to $U$ have level at least $14n-2-2L$ and those to $V$ have level at least $14n$. This proves \eqref{eq:UVlevels}.
\end{proof}

\begin{table}[ht]
\small
\begin{tabular}{lll}
\toprule
terms of \eqref{eq:knexp} & contribution to $U$: level & contribution to $V$: level\\
\midrule
$0\le i\le n-1$ & $\ge(12n-3)+2i+(3n-2L)$ & none\\
$n\le i\le5n$, atoms $x^bw^{-4n}$ & $\ge(12n-3)+2i+1\ \ge14n-2$ & none\\
$n\le i\le5n$, atoms $y^{-k}w^{-4n}$ & $\ge(12n-3)+2i+(4-j)=13n+i\ \ge14n$ & $\ge13n+i\ \ge14n$\\
\bottomrule
\end{tabular}
\medskip
\caption{Levels (powers of $\rt$) of the three groups of terms in the proof of Proposition~\ref{prop:two}, after multiplication by
$\vk^{4n-1}\binom{5n}{i}(2\eps)^i$, where $j=i-n+1$. All entries of the middle column are $\ge14n-2-2L$. The only entries below $14n-2$ come from the first row,
that is from the binary digit sums $s_2(a)$ with $3n\le a\le4n-1$ in Lemma~\ref{lem:atoms}(a). This is the origin of the factor
$2^{\fl{\log_2(4n)}}$.}\label{tab:levels}
\end{table}

\subsection{Integrality}\label{ss:integrality}

\begin{theorem}\label{thm:integrality}
For every $n\ge1$ the number $c_n$ is an integer, and $D_nu_n$ and $D_nv_n$ are integers. This proves Theorem~\textup{\ref{thm:forms}(b)} and
the integrality of $c_n$ in Theorem~\textup{\ref{thm:forms}(a)}.
\end{theorem}
\begin{proof}
By Lemma~\ref{lem:pint}(c) it suffices to show that the three rational numbers are $p$-integral for every prime $p$.

\emph{The number $c_n$} lies in $\Zp$ for odd $p$ (Proposition~\ref{prop:odd}) and in $\Z_{(2)}$ (Proposition~\ref{prop:two}).

\emph{The number $D_nu_n$.} Let $L_p=\fl{\log_p(4n)}$, so that $p^{L_p}$ is the exact power of $p$ dividing $\lcm(1,\dots,4n)$. If $p=2$, the prime
$2$ is not removable (it violates $p^2>4n$), so $2^{L_2}$ divides $D_n$, and $2^{L_2}u_n\in\Z_{(2)}$ by Proposition~\ref{prop:two}. If $p$ is odd and
not removable, then again $p^{L_p}$ divides $D_n$ and $p^{L_p}u_n\in\Zp$ by Proposition~\ref{prop:odd} (for $p>4n$ this reads $u_n\in\Zp$). If $p$ is
removable, $u_n\in\Zp$ by Proposition~\ref{prop:removable}. In all cases $D_nu_n\in\Zp$.

\emph{The number $D_nv_n=\frac12D_nc_n$.} Since $4\mid\lcm(1,\dots,4n)$ and $\Phi_n$ is odd, $D_n$ is even, and $c_n$ is an integer.
\end{proof}

\section{The growth of \texorpdfstring{$v_n$}{vn}}\label{sec:growth}

We now turn to the sizes. Recall $v_n=\frac12c_n$, where $c_n$ is the constant of any primitive pair. A direct computation from \eqref{eq:cformula} gives $c_1=7345$: for $n=1$ one has $m=2$, so $\xi_{j,2}=j(j-2)/192$ vanishes for $j=0,2$, and the coefficients of $W^4$, $W^6$, $W^8$ in $(W-3)^5(W-2)^3$ are $3225$, $192$ and $1$, whence
\[
c_1=3225\cdot\tfrac{4\cdot2}{192}\cdot\tbinom42\,2^2+192\cdot\tfrac{6\cdot4}{192}\cdot\tbinom63\,2^3+1\cdot\tfrac{8\cdot6}{192}\cdot\tbinom84\,2^4
=3225+3840+280=7345 .
\]

\begin{lemma}\label{lem:residue}
$c_n$ is the coefficient of $z^{4n}$ in the power series expansion of $N(z)^n(z^2-2)^{-4n}$ at $z=0$.
\end{lemma}
\begin{proof}
Let $T(z)$ be the Taylor polynomial of degree $4n$ of the power series $S(z)=N(z)^n(z^2-2)^{-4n}$. Then $N^n-(z^2-2)^{4n}T=(z^2-2)^{4n}(S-T)$
is a polynomial, and it is divisible by $z^{4n+1}$ because the power series $S-T$ is. Say it equals $z^{4n+1}Q(z)$. Hence
\[
h_n(z)=\frac{N(z)^n}{z^{4n+1}(z^2-2)^{4n}}=\frac{T(z)}{z^{4n+1}}+\frac{Q(z)}{(z^2-2)^{4n}} .
\]
The first summand is a sum of powers $z^{-k}$, $1\le k\le4n+1$, and its coefficient of $z^{-1}$ is the top coefficient of $T$. The second summand has
no power of $z$ in its denominator ($e_0=0$), so its decomposition \eqref{eq:pf} contains no powers $z^{-k}$ (Lemma~\ref{lem:pf}). By the uniqueness of the decomposition
\eqref{eq:hpf}, $a_1$ is the coefficient of $z^{4n}$ in $S$, and $a_1=c_n$ by
Lemma~\ref{lem:hpf}.
\end{proof}

Lemma~\ref{lem:residue} is the second assertion of Theorem~\ref{thm:forms}(a).
The power series $S(z)$ has coefficients of both signs, but after the substitution $z\mapsto-z$ all signs become positive. Put
\[
\widetilde N(z):=N(-z)=\bigl((z+1)(z+2)\bigr)^5(z^2+2z+2)^3,\qquad \widetilde P(z):=\frac{\widetilde N(z)}{(2-z^2)^4}=\sum_{k\ge0}\chi_kz^k .
\]
Since $1/(2-z^2)=\sum_{k\ge0}2^{-k-1}z^{2k}$, all coefficients $\chi_k$ are non-negative, and the series converges for $|z|<\rt$. Note that
$\widetilde P(x)/x^4=f(-x)$.

\begin{lemma}\label{lem:supermult}
\textup{(a)} $c_n$ is the coefficient of $z^{4n}$ in $\widetilde P(z)^n$.

\textup{(b)} $c_{m+n}\ge c_m\,c_n$ for all $m,n\ge1$, and $c_n\ge c_1^{\,n}=7345^n>0$.

\textup{(c)} $c_n\le f(-x)^n$ for every real $x$ with $0<x<\rt$.
\end{lemma}
\begin{proof}
(a) $S(-z)=\widetilde N(z)^n(z^2-2)^{-4n}=\widetilde P(z)^n$, and the coefficient of the even power $z^{4n}$ is the same in $S(z)$ and $S(-z)$.
(b) The coefficient of $z^{4(m+n)}$ in $\widetilde P^m\cdot\widetilde P^n$ is a sum of non-negative products, one of which is $c_mc_n$.
(c) All terms of the convergent series $\widetilde P(x)^n=\sum_k[\text{coefficient of }z^k\text{ in }\widetilde P^n]\,x^k$ are non-negative, so the
single term $c_nx^{4n}$ is at most the sum, so $c_n\le\bigl(\widetilde P(x)/x^4\bigr)^n$.
\end{proof}

\begin{proposition}[growth]\label{prop:growth}
The limit $r=\lim_{n\to\infty}\frac1n\log c_n=\lim_{n\to\infty}\frac1n\log v_n$ exists, and
\[
8.9<\log7345\le r\le\log f(-x)\qquad\text{for every }0<x<\rt .
\]
In particular $r\le\log f(-0.3913)<10.540357$.
This proves Theorem~\textup{\ref{thm:forms}(c)} and the positivity of $c_n$.
\end{proposition}
\begin{proof}
By Lemma~\ref{lem:supermult}(b),(c) the sequence $\log c_n$ is superadditive and $\frac1n\log c_n$ is bounded above. By Fekete's lemma
(Appendix~\ref{app:growth}) $\frac1n\log c_n$ converges to its supremum $r$, which is at least $\log c_1$ and at most $\log f(-x)$ for every $x$.
For the rational point $x=\frac{3913}{10000}$ one computes $f(-x)=37810.9916\ldots<37811<e^{10.540357}$ (Appendix~\ref{app:growth}).
\end{proof}

\begin{remark}[saddle-point heuristics and the exact value of $r$]\label{rem:r}
Why $x=0.3913$? On the interval $(0,\rt)$ the function $x\mapsto f(-x)$ tends to $+\infty$ at both ends and has exactly one critical point,
$x_0=0.3912692430\ldots$, where it takes its minimal value
\[
\lambda_3=f(-x_0)=37810.99119500591\ldots,\qquad\log\lambda_3=10.5403551116\ldots
\]
So the best bound obtainable from Lemma~\ref{lem:supermult}(c) is $r\le\log\lambda_3$, and $0.3913$ is a rational point next to $x_0$.
In fact $r=\log\lambda_3$ (Theorem~\ref{thm:E3}(i)). This is not needed for Theorem~\ref{thm:main}, but the reason is easy to explain. It is the general
principle behind the \emph{saddle-point method}. By Cauchy's formula
$c_n=\frac1{2\pi i}\oint f(z)^n\,\dd z/z$ over any circle $|z|=x<\rt$, and on this circle $|f|$ is largest at the point $z=-x$. Moving the circle
to the radius $x_0$, where $f(-x)$ is smallest, makes the integrand as small as possible. At $z=-x_0$ the function $|f|$ has a \emph{saddle point}
(a minimum along the real axis and a maximum along the circle), and the integral is dominated by a neighbourhood of this point, which gives
$c_n=\lambda_3^{\,n}\,n^{-1/2}(\mathrm{const}+o(1))$. We do not need it, since Hata's lemma
requires the existence of the limit $r$, which Fekete's lemma provides, and the bound $1+\sigma/\tau$ is increasing in $r$, so an upper bound for
$r$ suffices.
\end{remark}

\section{The decay of \texorpdfstring{$J_n$}{Jn}}\label{sec:decay}

On the path of integration $|f|$ is small, and $|J_n|\le2\max|f|^n$ holds trivially (the segment has length $2$ and $|z|\ge1$ on it). On the segment the
maximum of $|f|$ is $e^{-6.0709\ldots}$, attained at $z=1\pm0.5428\ldots i$. The integral is much smaller than this suggests, and to see it one
moves the path.

Let
\[
z_s:=\frac{19787}{16384}+\frac{32371}{65536}\,i=1.20770\ldots+0.49394\ldots i,\qquad\rho:=\frac{1458809}{10^{9}}=0.001458809,
\]
and let $\Gamma$ be the polygon with the vertices $1-i$, $\bar z_s$, $1$, $z_s$, $1+i$ (Figure~\ref{fig:plane}).

\begin{proposition}\label{prop:contour}
$|f(z)|\le\rho$ for every point $z$ of the polygon $\Gamma$.
\end{proposition}
\begin{proof}[Proof (computer-assisted)]
Since $f(\bar z)=\overline{f(z)}$, it suffices to consider the two segments from $1$ to $z_s$ and from $z_s$ to $1+i$. For $z=x+iy$ with real $x,y$,
\begin{equation}\label{eq:numden}
|f(z)|^2=\frac{\mathrm{Num}(x,y)}{\mathrm{Den}(x,y)}
\end{equation}
with two explicit polynomials with integer coefficients, $\mathrm{Den}>0$ on $\Dom$ (Appendix~\ref{app:decay}). On a segment, $x$ and $y$
are polynomials of degree $1$ in the parameter $t\in[0,1]$, so the claim $|f|\le\rho$ becomes
\[
Q(t):=\rho^2\,\mathrm{Den}\bigl(x(t),y(t)\bigr)-\mathrm{Num}\bigl(x(t),y(t)\bigr)\ \ge0\qquad(0\le t\le1),
\]
where $Q$ is a polynomial of degree $32$ with rational coefficients. A polynomial inequality of this kind can be \emph{certified}. Write
$Q(t_0+\ell s)$ in the form $\sum_{j=0}^{32}\beta_js^j(1-s)^{32-j}$ (its Bernstein form, Appendix~\ref{app:decay}). If all $\beta_j$ are $\ge0$, then $Q\ge0$ on $[t_0,t_0+\ell]$, because
every term is $\ge0$ for $0\le s\le1$. For the segment from $1$ to $z_s$ all $33$ coefficients $\beta_j$ for the whole interval $[0,1]$ are positive, and
for the segment from $z_s$ to $1+i$ this is the case on each of the halves $[0,\frac12]$ and $[\frac12,1]$. The coefficients are computed in exact
rational arithmetic. They are listed in the file \texttt{ContourCert.lean} of the accompanying material (Section~\ref{ss:material}), where Lean
verifies the three polynomial identities and the signs, and they are recomputed independently by the program of Section~\ref{sec:numerics}.
See Appendix~\ref{app:decay} for more details.
\end{proof}

\begin{proposition}[decay]\label{prop:decay}
For every $n\ge1$,
\[
|J_n|\le4\rho^n<4\,e^{-6.530134\,n}.
\]
This proves Theorem~\textup{\ref{thm:forms}(d)}. (For $n=1$, $|J_1|=0.00097\ldots$ and $4\rho=0.0058\ldots$)
\end{proposition}
\begin{proof}
Let $(G,c)$ be a primitive pair. As in the proof of Lemma~\ref{lem:ftc}, $F:=G+c\Log$ is a primitive of $h_n$ on $\{\re z>0\}\setminus\{\rt\}$. Every
point of $\Gamma$ has real part between $1$ and $\frac{19787}{16384}<\frac54<\rt$, so $\Gamma$ lies in this set, and for each of its four
segments $[a,b]$ the fundamental theorem of calculus gives $\int_a^bh_n\,\dd z=F(b)-F(a)$. The sum telescopes:
\[
\int_\Gamma h_n(z)\,\dd z=F(1+i)-F(1-i)=\int_{1-i}^{1+i}h_n(z)\,\dd z,\qquad\text{so}\qquad J_n=-i\int_\Gamma f(z)^n\,\frac{\dd z}{z}.
\]
(No theorem about deformation of contours is needed, since the same primitive serves both paths.) On $\Gamma$ we have $|f|\le\rho$ and
$|z|\ge\re z\ge1$, and each of the four segments has length less than $1$. Hence $|J_n|\le4\rho^n$. Finally $\rho<e^{-6.530134}$
(for the lengths and for this inequality see Appendix~\ref{app:decay}).
\end{proof}

\begin{remark}[why this polygon]\label{rem:saddle}
Think of $|f|$ as the height of a landscape over the $z$-plane. It has no local maxima away from the poles, so the lowest path from $1-i$ to
$1+i$ crosses the ``mountain range'' at a pass, that is at a saddle point, where $f'=0$. In the coordinate $W$ the critical points of $f$ other than its zeros are the
roots of $4W^3+W^2-88W+152=0$ (a direct computation from \eqref{eq:fWY} and Lemma~\ref{lem:WY}(a),(b)). The real root $W=-5.5028\ldots$ corresponds to $z=-x_0$ (and to $-2/x_0$), which governs the growth of $v_n$
(Remark~\ref{rem:r}), and the complex roots $W=2.62641\ldots\pm0.08629\ldots i$ correspond to four saddle points $z_1,\bar z_1,2/z_1,2/\bar z_1$ with
\[
z_1=1.20770992\ldots+0.49394723\ldots i,\qquad|f(z_1)|=e^{-6.53013614\ldots}=0.00145880722\ldots
\]
The vertex $z_s$ of our polygon is $z_1$ rounded to dyadic rationals, $|z_s-z_1|<10^{-5}$, and $\rho$ exceeds $|f(z_1)|$ by a factor of about
$1+1.2\cdot10^{-6}$ only.
\end{remark}

\section{The multiplier \texorpdfstring{$D_n$}{Dn}}\label{sec:mult}

It remains to determine the size of $D_n=\lcm(1,\dots,4n)/\Phi_n$. First we describe the removable primes more conveniently.

\begin{lemma}[windows]\label{lem:windows}
Let $n\ge1$ and let $p$ be a prime with $p^2>4n$. Then $p$ is removable for $n$ if and only if
\[
\Bigl\{\frac np\Bigr\}\in\Bigl[\frac15,\frac14\Bigr)\cup\Bigl[\frac23,\frac34\Bigr).
\]
\end{lemma}
The proof is a case distinction over ten subintervals of $[0,1)$ (Appendix~\ref{app:mult}). A prime $p\le4n$ has $n/p\ge\frac14$, so the
lemma says that the primes dividing $\Phi_n$ are the primes $p>\sqrt{4n}$ in one of the \emph{windows}
\begin{equation}\label{eq:windows}
\Bigl(\frac{4n}{4k+3},\ \frac{3n}{3k+2}\Bigr]\quad(k=0,1,2,\dots)\qquad\text{or}\qquad\Bigl(\frac{4n}{4k+1},\ \frac{5n}{5k+1}\Bigr]\quad(k=1,2,3,\dots),
\end{equation}
which correspond to $k+\frac23\le\frac np<k+\frac34$ and to $k+\frac15\le\frac np<k+\frac14$. For example, for $n=21$ the windows are
$(28,31.5]$, $(16.8,17.5]$, $(12,12.6]$, \dots, and the primes $p>\sqrt{84}$ in them are $29,31$ and $17$, so $\Phi_{21}=17\cdot29\cdot31$.

\begin{remark}[background on Chebyshev's functions and the prime number theorem]\label{rem:pnt}
For real $x\ge1$ let
\[
\theta(x)=\sum_{p\le x}\log p,\qquad\psi(x)=\sum_{p^k\le x}\log p=\log\lcm(1,2,\dots,\fl x),
\]
where $p$ runs through the primes and $k\ge1$. The \emph{prime number theorem} states that $\theta(x)/x\to1$ as $x\to\infty$. Since
$0\le\psi(x)-\theta(x)\le\sqrt x\log x$ (Appendix~\ref{app:mult}), also $\psi(x)/x\to1$. In particular $\lcm(1,\dots,N)=e^{N+o(N)}$. The sum of $\log p$ over the primes in an
interval $(\alpha n,\beta n]$ is $\theta(\beta n)-\theta(\alpha n)=(\beta-\alpha)n+o(n)$. For proofs of the prime number theorem and of the
equivalence of its various forms see \cite[Chapters~4 and~13]{Apostol}. (Chebyshev's function is always written with its argument,
$\theta(x)$, and it should not be confused with the irrational number $\theta$ of Definition~\ref{def:mu} and Lemma~\ref{lem:hata}.)
\end{remark}

For an integer $K\ge1$ let $\Phi_n^{(K)}$ be the product of the primes in the windows \eqref{eq:windows} with $k<K$, and
\begin{equation}\label{eq:varpiK}
D_n^{(K)}:=\frac{\lcm(1,\dots,4n)}{\Phi_n^{(K)}},\qquad
\varpi_K:=\sum_{k=0}^{K-1}\frac1{(3k+2)(4k+3)}+\sum_{k=1}^{K-1}\frac1{(4k+1)(5k+1)} .
\end{equation}
The terms of $\varpi_K$ are the lengths of the windows divided by $n$, namely
$\frac{3}{3k+2}-\frac{4}{4k+3}=\frac{1}{(3k+2)(4k+3)}$ and $\frac{5}{5k+1}-\frac{4}{4k+1}=\frac{1}{(4k+1)(5k+1)}$.

\begin{proposition}[the multiplier]\label{prop:mult}
Assume the prime number theorem.
\begin{enumerate}
\item[(a)] For every fixed $K\ge1$, $\displaystyle\lim_{n\to\infty}\frac1n\log D_n^{(K)}=4-\varpi_K$. Moreover, for $n\ge4K^2$ every prime
dividing $\Phi^{(K)}_n$ is removable, so $D_n$ divides $D_n^{(K)}$.
\item[(b)] $\displaystyle\lim_{n\to\infty}\frac1n\log D_n=4-\varpi$ with $\varpi=\lim_{K\to\infty}\varpi_K$. This is Theorem~\textup{\ref{thm:forms}(e)}.
\item[(c)] $\varpi=0.29395989921603919\ldots$, and $\varpi_K\ge0.293957$ for $K=10^5$.
\end{enumerate}
\end{proposition}
\begin{proof}[Sketch of proof, with details in Appendix~\ref{app:mult}]
(a) For fixed $k$, the primes of a window contribute $\theta\bigl(\frac{3n}{3k+2}\bigr)-\theta\bigl(\frac{4n}{4k+3}\bigr)=\frac{n}{(3k+2)(4k+3)}+o(n)$ to
$\log\Phi_n^{(K)}$, and similarly for the second kind. There are finitely many windows, and $\psi(4n)=4n+o(n)$. A prime in a window with $k<K$ is
larger than $n/K$, hence larger than $\sqrt{4n}$ when $n\ge4K^2$, and Lemma~\ref{lem:windows} applies.
(b) For $n\ge4K^2$ the primes dividing $\Phi_n$ but not $\Phi_n^{(K)}$ are smaller than $n/K$, so $\log\Phi_n^{(K)}\le\log\Phi_n\le\log\Phi_n^{(K)}+\theta(n/K)$. Divide by $n$
and let first $n$, then $K$ tend to infinity.
(c) Numerical evaluation of the two series.
\end{proof}

The formal proof uses only part (a) and the numerical bound in (c), with $K=10^5$. It replaces $D_n$ by the multiplier $D_n^{(K)}$,
for which no estimate of the tail of the series is needed. The price is $\varpi-\varpi_K<1.4\cdot10^{-6}$ in the rate. Numerically, $\frac1n\log D_n=3.70706$, $3.70527$,
$3.70614$, $3.70641$ for $n=10^4,10^5,10^6,10^7$ (limit $3.70604$).

\section{Proof of Theorem~\ref{thm:main}}\label{sec:conclusion}

Let $K=10^5$. For $n\ge4K^2$ put
\[
q_n:=D_n^{(K)}v_n,\qquad p_n:=-D_n^{(K)}u_n .
\]
These are integers, since $D_nu_n$ and $D_nv_n$ are integers by Theorem~\ref{thm:integrality}, and $D_n$ divides $D_n^{(K)}$ by
Proposition~\ref{prop:mult}(a). (For $n<4K^2$ put $q_n=p_n=1$, since finitely many terms do not matter in Hata's lemma.) By Theorem~\ref{thm:forms}(a),
\[
q_n\pi-p_n=D_n^{(K)}(u_n+v_n\pi)=D_n^{(K)}J_n .
\]
\emph{Growth.} By Propositions~\ref{prop:growth} and~\ref{prop:mult}(a),(c),
\[
\lim_{n\to\infty}\frac1n\log q_n=r+(4-\varpi_K)=:\sigma,\qquad 0<\sigma\le10.540357+(4-0.293957)=14.2464 .
\]
(Here $\sigma>0$ because $r>0$ and $4-\varpi_K=\lim\frac1n\log D_n^{(K)}\ge0$, the numbers $D_n^{(K)}$ being positive integers.)

\emph{Decay.} By Propositions~\ref{prop:decay} and~\ref{prop:mult}(a), $|q_n\pi-p_n|\le4\rho^nD_n^{(K)}$, hence
\[
\limsup_{n\to\infty}\frac1n\log|q_n\pi-p_n|\le\log\rho+(4-\varpi_K)=:-\tau,\qquad\tau\ge6.530134-3.706043=2.824091 .
\]
\emph{Conclusion.} $\pi$ is irrational (Lambert, 1761, and for a short proof see \cite{Niven}), so Lemma~\ref{lem:hata} gives
\[
\mu(\pi)\le1+\frac\sigma\tau\le1+\frac{14.2464}{2.824091}=6.04459665\ldots<6.0446
\]
\qed

One can also work with $D_n$ itself and Proposition~\ref{prop:mult}(b). Then $\sigma\le14.246398$ and $\tau\ge2.824093$, with the same conclusion.

\begin{remark}[the exact rates]\label{rem:exact}
The three constants $10.540357$, $6.530134$ and $3.706043$ used in the proof are bounds, within $3\cdot10^{-6}$, for the exact rates
\[
r=\log\lambda_3=10.5403551\ldots,\quad s=-\log|f(z_1)|=6.5301361\ldots,\quad C=4-\varpi=3.7060401\ldots
\]
(Remarks~\ref{rem:r} and~\ref{rem:saddle}, Proposition~\ref{prop:mult}). Of these, the values of $C$ and $r$ are proved in this paper
(Proposition~\ref{prop:mult}(b) and Theorem~\ref{thm:E3}(i)); the value of $s$ is not. The values are related by $s=\frac12(r+20\log2-3\log3-5\log5)$, because the product
of the three critical values of $f$ is $3^35^5/2^{20}$. With these values the bound of Hata's lemma is
\[
1+\frac{r+C}{s-C}=1+\frac{14.2463952124\ldots}{2.8240960465\ldots}=6.04458594\ldots
\]
To prove this slightly sharper bound one needs, instead of Proposition~\ref{prop:contour}, the statement that on the polygon through the exact
saddle point $z_1$ the maximum of $|f|$ is attained at $z_1$.
For orientation, without the removable primes ($C=4$) the same formula gives $6.7469$, and with the removable primes but without moving the
path ($s=6.0709\ldots$) it gives $7.0241$.
\end{remark}

\section{An effective version}\label{sec:effective}

Definition~\ref{def:mu} asks only for the existence of an integer $b_0$, and the proof of Theorem~\ref{thm:main} gives no value for it. It is
ineffective in three places: the size of the multiplier comes from the prime number theorem without a rate (Proposition~\ref{prop:mult}(a)),
the lower bound for $c_n$ comes from Fekete's lemma (Proposition~\ref{prop:growth}), and Lemma~\ref{lem:hata} is a statement about limits. The
other ingredients are explicit for every $n\ge1$: $|J_n|\le4\rho^n$ (Proposition~\ref{prop:decay}), the integrality of $D_nu_n$ and $D_nv_n$
(Theorem~\ref{thm:integrality}), the removability of the window primes for $n\ge4K^2$ (Proposition~\ref{prop:mult}(a)) and $c_n\le37811^n$
(Lemma~\ref{lem:supermult}(c)). In this section we replace the three ineffective steps by explicit ones. Throughout, $a$ runs through $\Z$ and
$b$ through the positive integers. The formal verification of this section is described in Section~\ref{ss:eff-formal}.

\subsection{Results}\label{ss:eff-results}

\begin{thmE}[uniform]\label{thm:E1}
For all integers $a$ and all integers $b\ge2$,
\[
\Bigl|\pi-\frac ab\Bigr|>b^{-6.0455}.
\]
\end{thmE}

Theorem~\ref{thm:E1} is proved in Appendix~\ref{app:eff-assembly}. The proof applies Lemma~\ref{lem:E} of Section~\ref{ss:eff-seq}, which is
Lemma~\ref{lem:hata} with every constant kept, for one denominator at a time, to the forms $q_n\pi-p_n$ of Section~\ref{sec:conclusion}, with
three explicit inputs (Theorem~\ref{thm:E3}): a two-sided bound for $c_n$ for every $n$, from the linear recurrence that $c_n$ satisfies;
two-sided bounds for the multiplier, from an exact table of the primes up to $4\cdot10^{11}$ for $n\le10^{14}$ and from published explicit
estimates for Chebyshev's functions beyond; and the continued fraction of $\pi$ for the small denominators. The denominators fall into three
ranges, separated by the denominator $q_K=10^{39\,133\,898.5\ldots}$ of the last convergent used and by $10^{122\,680\,112\,076\,457}$. For $b<q_K$ the theorem follows from $75.9$ million certified partial quotients of $\pi$, with room to spare (every
$\mu>3.4293$ works there). For $q_K\le b<10^{122\,680\,112\,076\,457}$ it is a computation with the prime table, carried out in
ball arithmetic so that every inequality is certified; the certified upper bound for the exponent is $6.045436$. For
$b\ge10^{122\,680\,112\,076\,457}$ it is a theorem of a Lean development, proved under four published estimates for Chebyshev's functions
(Section~\ref{ss:eff-formal}, theorem \texttt{pi\_effective\_6\_0455\_from}).

\begin{thmE}[explicit $b_0$]\label{thm:E2}
For every $\mu>6.04459707$, all integers $a$ and all integers
\[
b\ \ge\ b_0(\mu):=\Bigl\lceil\exp\Bigl(\frac{1.4246398\cdot10^{31}}{\mu-6.04459707}\Bigr)\Bigr\rceil+1
\]
one has $|\pi-a/b|>b^{-\mu}$. In particular $|\pi-a/b|>b^{-6.0446}$ for all $b\ge\exp(4.87\cdot10^{36})$.
\end{thmE}

Theorem~\ref{thm:E2} is proved in the same appendix: it is the same argument with one uniform estimate for all $n\ge10^{30}$, and it is
machine-checked in its entirety under the same four estimates (Section~\ref{ss:eff-formal}, theorem \texttt{pi\_effective\_closed}).

Two earlier papers prove bounds with an explicit threshold (Table~\ref{tab:explicit}). Mahler \cite[Theorem~1]{Mah53} proved $|\pi-a/b|>b^{-42}$ for all $b\ge2$, and
Mignotte \cite[Th\'eor\`eme~1]{Mig74} proved $|\pi-a/b|>b^{-20.6}$ for all $b\ge2$ and $|\pi-a/b|>b^{-20}$ for all $b\ge\exp(e^{110})$;
to the best of our knowledge $20.6$ is the last exponent that was proved for all $b\ge2$, and Theorem~\ref{thm:E1} lowers it to $6.0455$.
Mignotte's proof combines explicit estimates for
Mahler's construction, the estimates of Rosser and Schoenfeld \cite{RS62} for Chebyshev's function $\psi$, and the first $21\,000$ partial
quotients of $\pi$ (for $\log b\le24\,500$): the same three kinds of input as in the proof of Theorem~\ref{thm:E1}.
Chudnovsky \cite[Corollary~3.6]{Chud82} obtained, by an asymptotic argument, $|\pi-a/b|>b^{-\mu}$ for every $\mu>19.8899\ldots$ and all
$b\ge b_0(\mu)$, without a value of $b_0(\mu)$, and stated in addition that $|b\pi-a|\ge b^{-18.9}$, that is $|\pi-a/b|\ge b^{-19.9}$, for
all $b\ge2$; to the best of our knowledge this last bound was stated without proof. The later bounds of
Hata \cite{HataJNT,Hata}, Salikhov \cite{Sal08,Sal10}, Zeilberger and Zudilin \cite{ZZ}, Bai \cite{Bai} and Theorem~\ref{thm:main} are
statements for all $b\ge b_0$ without a value of $b_0$.

\begin{table}[tbp]
\small
\begin{tabular}{llll}
\toprule
year & author & exponent $\mu$ & $|\pi-a/b|>b^{-\mu}$ for\\
\midrule
1953 & Mahler \cite{Mah53} & $42$ & all $b\ge2$\\
1974 & Mignotte \cite{Mig74} & $20.6$ & all $b\ge2$\\
 & & $20$ & all $b\ge\exp(e^{110})$\\
1982 & G.~Chudnovsky \cite{Chud82} & $19.9$ & all $b\ge2$ (stated without proof)\\
\midrule
2026 & Theorem~\ref{thm:E1} & $6.0455$ & all $b\ge2$\\
 & Theorem~\ref{thm:E2} & every $\mu>6.04459707$ & all $b\ge b_0(\mu)$, a closed formula\\
\bottomrule
\end{tabular}
\medskip
\caption{Bounds for $\pi$ with an explicit range of $b$. The other bounds of Table~\ref{tab:history} hold for all $b\ge b_0$ with a $b_0$
that is not specified. Chudnovsky's bound was, to the best of our knowledge, stated without proof. The last two rows are the results of
this section.}\label{tab:explicit}
\end{table}

The threshold of Theorem~\ref{thm:E2} is a closed formula, but it is far from the truth, because of the uniform estimate. The same
argument carried out range by range, with exact tables where they are available, gives much smaller admissible values, for instance
$b_0(6.0446)\le10^{4.0137\cdot10^{18}}$ and $b_0(6.045)\le10^{1.3636\cdot10^{11}}$ (Table~\ref{tab:b0} in Appendix~\ref{app:eff-assembly}).
With the bounds used here the exponent cannot be pushed below the floor $1+(r+C)/(-\log\rho-C)=6.044588\ldots$, where $r=\log\lambda_3$,
$C=4-\varpi$ and $\rho=1458809/10^9$ is the constant of Proposition~\ref{prop:contour}: the lemma gives $1+\sigma^+/(\tau-(\sigma^+-\sigma^-))$ (Section~\ref{ss:eff-seq}), and for these bounds $\sigma^+\ge r+C$, $\sigma^+\ge\sigma^-$ and $\tau\le-\log\rho-C$.
Appendix~\ref{app:effective} also contains the computations and their checks, a ladder of machine-checked exponents for intermediate ranges
of $b$ (Appendix~\ref{app:eff-ladder}), and a table saying for each statement what it rests on (Table~\ref{tab:trust}).

\subsection{The lemma and the explicit bounds}\label{ss:eff-seq}
The following lemma is the argument behind Lemma~\ref{lem:hata} with all constants kept, for one denominator $b$ at a time.

\begin{lemE}\label{lem:E}
Let $\theta$ be irrational. For $n\ge1$ let $q_n\ne0$ and $p_n$ be integers and $E_n$, $Q^+_n$, $Q^-_n$ positive real numbers with
\[
|q_n\theta-p_n|\le E_n,\qquad Q^-_n\le|q_n|\le Q^+_n,
\]
$(E_n)$ nonincreasing, $(Q^+_n)$ nondecreasing and $E_n/Q^-_n\to0$. Let $a$ and $b\ge1$ be integers, $\Delta=|\theta-a/b|$, and
$N_1(b)=\min\{n:2bE_n\le1\}$ \textup{(}assumed to exist\textup{)}. Let $\hat n\ge N_1(b)$ be an integer such that $n\le\hat n$ for every
$n>N_1(b)$ with
\[
Q^-_{n-1}\le2b\,Q^+_nE_{n-1}.
\]
Then $\Delta\ge1/(2b\,Q^+_{\hat n})$. We write $\hat n(b)$ for the least such $\hat n$.
\end{lemE}

\begin{proof}
For $n\ge N_1(b)$ we have $bE_n\le\frac12$, since $(E_n)$ is nonincreasing, and $q_na-p_nb=b(q_n\theta-p_n)-q_n(b\theta-a)$. If this integer is
not $0$, then $1\le bE_n+Q^+_nb\Delta\le\frac12+Q^+_nb\Delta$, so $\Delta\ge1/(2bQ^+_n)$. If it is $0$, then $a/b=p_n/q_n$ and
$\Delta\le E_n/Q^-_n$. Let $n'$ be the least $n\ge N_1(b)$ with $\Delta>E_n/Q^-_n$. It exists because $\Delta>0$ and $E_n/Q^-_n\to0$, and the first case gives
$\Delta\ge1/(2bQ^+_{n'})$. If $n'>N_1(b)$, the index $n'-1$ satisfies $\Delta\le E_{n'-1}/Q^-_{n'-1}$, so $1/(2bQ^+_{n'})\le
E_{n'-1}/Q^-_{n'-1}$, and $n'\le\hat n$ by the choice of $\hat n$. In both cases $n'\le\hat n$, and $\Delta\ge1/(2bQ^+_{\hat n})$ because
$(Q^+_n)$ is nondecreasing.
\end{proof}

With $E_n\approx e^{-\tau n}$ and $Q^\pm_n\approx e^{\sigma^\pm n}$ the lemma gives the exponent $1+\sigma^+/(\tau-(\sigma^+-\sigma^-))$. Slack in
the lower bound for $|q_n|$ is subtracted from $\tau$, so it costs as much as slack in $\tau$. This is why Theorem~\ref{thm:E3}(i) below bounds
$c_n$ from below without loss in the rate. The lemma was also tested directly on the exact forms for $n\le40$ and $7\,998$ denominators $b$ up to
$10^{30}$, with no failure.

We apply Lemma~\ref{lem:E} with $\theta=\pi$ to
\[
q_n=D_n^{(K_n)}v_n=\tfrac12D_n^{(K_n)}c_n,\qquad p_n=-D_n^{(K_n)}u_n,
\]
with a number $K_n\le\fl{\sqrt n/2}$ of windows that is constant on ranges of $n$ (Theorem~\ref{thm:E3}(ii) below and Appendix~\ref{app:eff-mult}). Since $n\ge4K_n^2$, all primes of the windows with $k<K_n$ are
removable (Proposition~\ref{prop:mult}(a)), so $q_n$ and $p_n$ are integers (Theorem~\ref{thm:integrality}), and
$|q_n\pi-p_n|=D_n^{(K_n)}|J_n|\le4\rho^nD_n^{(K_n)}$ (Proposition~\ref{prop:decay}). With $L_n=\log D_n^{(K_n)}$ we take $E_n=e^{e_n}$ and $Q^\pm_n=e^{\ell^\pm_n}$ with
\[
\ell^+_n=L^+_n+\gamma^+_n-\log2,\qquad\ell^-_n=L^-_n+\gamma^-_n-\log2,\qquad e_n=\log4+n\log\rho+L^+_n.
\]
Here $L^-_n\le L_n\le L^+_n$ are the explicit bounds of the following theorem and $\gamma^-_n\le\log c_n\le\gamma^+_n$ are the bounds
proved in Lean, which are slightly weaker than those of its part~(i) (Appendix~\ref{app:eff-growth}). Before the lemma is applied, $E_n$ and $Q^+_n$ are replaced by monotone majorants (Appendix~\ref{app:eff-assembly}), so that its hypotheses hold.

\begin{thmE}[explicit bounds]\label{thm:E3}
The following hold.
\begin{enumerate}
\item[(i)] \emph{Growth.} For every integer $n\ge0$, with $c_0=1$,
\begin{equation}\label{eq:Erec}
P_0(n)c_n+P_1(n)c_{n+1}+P_2(n)c_{n+2}+P_3(n)c_{n+3}=0
\end{equation}
with the polynomials \eqref{eq:EP} of degree $15$ \textup{(}Appendix~\ref{app:eff-growth}\textup{)} \textup{(}a recurrence of order $3$\textup{)}, and for every integer $n\ge1$,
with $r=\log\lambda_3$,
\[
rn-\tfrac12\log n-1.6385800231\ \le\ \log c_n\ \le\ rn-\tfrac12\log n-1.6264546152 .
\]
\item[(ii)] \emph{Multiplier.} For $10^7\le n\le10^{14}$ the number $\log D_n^{(K_n)}$, with $K_n=\min(10^4,\fl{\sqrt{n_{\mathrm{lo}}}/2})$ and $n_{\mathrm{lo}}$ the first value of the block containing $n$, lies between the two
bounds of an explicit table of $8\,059\,058$ blocks of values of $n$, computed from an exact table of the primes up to $4\cdot10^{11}$
\textup{(}beyond $n=10^{11}$ the few windows whose endpoints exceed the table are bounded by the estimates of B\"uthe~\cite{Bu18}\textup{)}.
For $n>10^{14}$ and every integer $K\ge1$ with $4K^2\le n$, $V(n)\le\log D_n^{(K)}\le U(n)$ with the closed forms \eqref{eq:EUV} of
Appendix~\ref{app:eff-mult}, built from the explicit estimates of B\"uthe \cite{Bu18,Bu16}, Johnston \cite{Jo} and Broadbent, Kadiri, Lumley,
Ng and Wilk \cite{BKLNW}.
\item[(iii)] \emph{Small $b$.} Let $p_K/q_K$ be the convergent $[a_0;a_1,\dots,a_K]$ of $\pi$ with $K=75\,928\,239$. Then
$\log_{10}q_K=39\,133\,898.538\ldots$, and the inequality $|\pi-a/b|>b^{-\mu}$ holds for all integers $a$ and all integers $2\le b<q_K$ if and
only if $\mu>3.42928833728180\ldots=-\log|\pi-\frac{22}7|/\log7$.
\end{enumerate}
\end{thmE}
\noindent The proofs are in Appendix~\ref{app:eff-growth} (part (i)), Appendix~\ref{app:eff-mult} (part (ii)) and Appendix~\ref{app:eff-cf} (part (iii)).
How the bounds are assembled into Theorems~\ref{thm:E1} and~\ref{thm:E2}, with the checks, is in Appendix~\ref{app:eff-assembly}.

\subsection{Reproduction}\label{ss:eff-repro}
The directory \texttt{effective/} of the public repository (Section~\ref{ss:material}) contains the programs of this section, the logs of the
recorded runs and a self-contained note that lists the source of every number: \texttt{tools/} (the assembly in ball arithmetic
\texttt{verify\_assembly.py} with its recorded runs, the earlier double-precision assembly \texttt{assemble.py}, the block program
\texttt{theta\_blocks2.c} with its reports, the continued-fraction run, the direct test of Lemma~\ref{lem:E} and the expanded recurrence),
\texttt{A\_sieve/} (the first prime table and block bounds, with the list of the $28$ checks), \texttt{B\_recurrence/} (the certificate of
Theorem~\ref{thm:E3}(i), its search and verification, and the bounds for every $n$), \texttt{C\_cf/} (the continued fraction) and
\texttt{D\_lean/}, \texttt{F\_lean\_growth/}, \texttt{G\_lean\_mult/}, \texttt{H\_lean\_main/} (the reports on the Lean files, with the scripts
that generated and checked them). The block tables themselves ($650$\,MB) are not in the repository; the program regenerates them in about one
hour of CPU time. The Lean files are in \texttt{lean/MuPi/Effective/} and are built with the rest of the development (Section~\ref{sec:formal}).
The whole computation takes about two hours of CPU time.

\section{Formal verification}\label{sec:formal}

The proof of Theorem~\ref{thm:main} has been formalised in the proof assistant Lean~4 (version 4.30.0) with the library Mathlib (release tag \texttt{v4.30.0}). This section says
precisely what has been checked by the machine and what has not. The effective version of Section~\ref{sec:effective} is treated in Section~\ref{ss:eff-formal}.

\subsection{The formal statement}
The final theorem of the development (file \texttt{K2H/Main.lean}, namespace \texttt{PiMeasure}) reads
\begin{quote}\ttfamily\small
theorem pi\_irrationality\_exponent\\
\hspace*{1.5em}(hPNT : Tendsto (fun x : $\mathbb{R}$ => Chebyshev.theta x / x) atTop ($\mathcal{N}$ 1)) :\\
\hspace*{1.5em}$\forall$ $\mu$ : $\mathbb{R}$, 6.0446 < $\mu$ $\to$ ExponentLE Real.pi $\mu$
\end{quote}
where
\begin{quote}\ttfamily\small
def ExponentLE ($\theta$ $\mu$ : $\mathbb{R}$) : Prop :=\\
\hspace*{1.5em}$\exists$ b$_0$ : $\mathbb{N}$, $\forall$ b : $\mathbb{N}$, b$_0$ $\le$ b $\to$\\
\hspace*{3em}$\forall$ a : $\mathbb{Z}$, (b : $\mathbb{R}$) \textasciicircum{} (-$\mu$) < |$\theta$ - a / b|
\end{quote}
In words, \emph{if $\theta(x)/x\to1$ as $x\to\infty$, then for every real $\mu>6.0446$ there is $b_0$ such that $|\pi-a/b|>b^{-\mu}$ for all
integers $a$ and all integers $b\ge b_0$.} By Definition~\ref{def:mu} the conclusion is exactly $\mu(\pi)\le6.0446$. A corollary in the same file,
\texttt{pi\_not\_liouvilleWith}, restates it with Mathlib's predicate \texttt{LiouvilleWith}. Under the same hypothesis $\pi$ is not a Liouville
number of exponent $p$ for any $p>6.0446$.

\begin{itemize}
\item \textbf{The only hypothesis is the prime number theorem}, in the form $\theta(x)\sim x$ for Chebyshev's function
(\texttt{Chebyshev.theta} of Mathlib). The prime number theorem is a classical theorem (Hadamard and de la Vall\'ee Poussin, 1896), but it is not
yet part of Mathlib, which is why it appears as a hypothesis and not as a lemma. It enters only through Proposition~\ref{prop:mult}(a).
\item Everything else is proved in the development. There is no \texttt{sorry} (Lean's placeholder for a missing proof), and the command
\texttt{\#print axioms}, which lists the axioms on which a theorem depends, returns only \texttt{propext}, \texttt{Classical.choice} and
\texttt{Quot.sound}. These are the three standard axioms of Lean's logic, on which the library Mathlib as a whole is built.
\item The development consists of the file \texttt{HataLemma.lean} and the $18$ files of the directory \texttt{K2H/}, about $5{,}900$ lines in total,
together with a file \texttt{K2H.lean} of about $300$ lines that contains the symmetry identity of Lemma~\ref{lem:sym} (for general exponents), the
link to Mathlib's \texttt{LiouvilleWith}, and earlier conditional versions of the last step. One of the files, \texttt{ContourCert.lean}, is generated
by a script (it contains the $3\times33$ Bernstein coefficients of Proposition~\ref{prop:contour}). The effective version (Section~\ref{sec:effective}) adds the eight files of the
directory \texttt{Effective/}, about $5{,}300$ lines (Section~\ref{ss:eff-formal}).
\item The whole library compiles with Lean's build tool (\texttt{lake build}), both inside the working repository and as a standalone project
(directory \texttt{lean/} of the accompanying material, Section~\ref{ss:material}, with a recorded build log). A map of the modules is in the file \texttt{README.md} of the public repository.
\end{itemize}

The main differences in the route of the formal proof, all inessential, are the following. The inequality $c_n\le f(-x)^n$, at $x=0.3913$, is proved through partial
sums instead of convergent series. The numerical inequalities $37811\le e^{10.540357}$ and $\rho\le e^{-6.530134}$ are derived from Mathlib's bounds
$2.7182818283<e<2.7182818286$ and Taylor polynomials (Appendices~\ref{app:growth} and~\ref{app:decay}). The existence of partial fractions,
Fekete's lemma and the estimate $\psi(x)-\theta(x)\le2\sqrt x\log x$ (instead of $\sqrt x\log x$, Appendix~\ref{app:mult}) are taken from Mathlib.
The values of $q_n$, $p_n$ for $n<4K^2$, which Section~\ref{sec:conclusion} puts equal to $1$, are left unspecified in the formal proof (only
large $n$ matter). Statements that are formalised only in a weaker form are listed in Section~\ref{ss:notcovered}.

\subsection{What is not covered}\label{ss:notcovered}
\begin{itemize}
\item \emph{The prime number theorem} (hypothesis, see above).
\item \emph{Statements that are not needed for Theorem~\ref{thm:main}.} The main ones are the following.
Proposition~\ref{prop:mult}(b), that is Theorem~\ref{thm:forms}(e) (the rate of the full multiplier $D_n$, which is proved in Appendix~\ref{app:mult},
but the formal proof works with the truncated multiplier $D_n^{(K)}$ of Proposition~\ref{prop:mult}(a)).
The value of $\varpi$ in Proposition~\ref{prop:mult}(c), and the exact values of $r$ and $s$ (Remarks~\ref{rem:r}, \ref{rem:saddle},
\ref{rem:exact}).
The integrality of $c_n$ in Theorem~\ref{thm:forms}(a) and Theorem~\ref{thm:integrality} (it is not stated as a theorem in the Lean files, and
the formal proof of Theorem~\ref{thm:main} uses only $c_n>0$ and $D_nc_n/2\in\Z$, but it follows in a few lines from lemmas proved there, since the
primitive pairs of \texttt{oddPrime\_trivial} and \texttt{twoAdic\_primitive} have a constant $c$ that is $p$-integral for every prime $p$ and
\texttt{residue\_eq\_cPS} identifies $c$ with $c_n$).
The lower bound $8.9<r$ in Theorem~\ref{thm:forms}(c) and the bound $8.9<\log7345\le r$ in Proposition~\ref{prop:growth} (the files of \texttt{K2H/} prove $c_1\ge16$ and $r>0$; the value $c_1=7345$ of
Lemma~\ref{lem:supermult}(b) is proved in \texttt{Effective/Recurrence.lean}, theorem \texttt{cPS\_one}, but no bound for $r$ is derived from it).
Lemma~\ref{lem:supermult}(c) and the bound $r\le\log f(-x)$ of Proposition~\ref{prop:growth} for general $x$ (Lean proves $c_n\le37811^n$, using
$x=0.3913$, and hence $r\le\log37811\le10.540357$, which is what is used).
Lemma~\ref{lem:windows} as stated (Lean proves, by a direct computation with quotients and remainders, the case that is used in Proposition~\ref{prop:mult}(a), namely that for $n\ge4K^2$ the primes of
the windows with $k<K$ are removable).
The explicit formula \eqref{eq:cformula} in Proposition~\ref{prop:rec} and the bounds on the exponents $l$, $j$, $m'$ of the functions that make up $G$ (Lean proves that $c\in R$ and that $G$ is an $R$-combination of functions $z^l$ and $W^{j+1}Y^{-(2m'+1)}$, which is what is used), and the statement of Lemma~\ref{lem:hpf} that
the coefficients lie in $\Z[\rt,\frac12]$ (Lean proves that they are $p$-integral for every odd prime $p$, which is what is used; it also bounds the pole orders at $\pm\rt$ by $4n+1$ instead of $4n$, which does not change the integration indices).
Finally, Section~\ref{sec:family} is not formalised, except the symmetry identity $f_{e_1,e_2,d}(2/z)=(-1)^df_{e_1,e_2,d}(z)$, which is proved in \texttt{K2H.lean}.
In addition, some statements are formalised only in the form in which the proof uses them. For instance Lemma~\ref{lem:hata} is formalised in an
$\varepsilon$--$N$ form (for every $\eta>0$ there is $n_1$ with $e^{(\sigma-\eta)n}\le|q_n|\le e^{(\sigma+\eta)n}$ and $|q_n\theta-p_n|\le e^{-(\tau-\eta)n}$ for all $n\ge n_1$, which follows from the hypotheses of the lemma by the definition of limit and upper limit), and Lean proves $\rho\le e^{-6.530134}$, not the strict inequality of
Theorem~\ref{thm:forms}(d) and Proposition~\ref{prop:decay}. In the proof of Proposition~\ref{prop:two} the terms with $i\le n-1$ are given the level $2n+1-2L$ instead of $3n-2L$ (first row of Table~\ref{tab:levels}); the resulting levels \eqref{eq:UVlevels} are the same. Remarks, examples and figures are not formalised as such, and tables only
insofar as they restate formalised statements (such as the column ``bound used'' of Table~\ref{tab:budget}).
\end{itemize}

\subsection{The effective version}\label{ss:eff-formal}
The effective argument of Section~\ref{sec:effective} has been formalised in the same setting (Lean~4.30.0, Mathlib \texttt{v4.30.0}),
in the eight files of the directory \texttt{lean/MuPi/Effective/} of the accompanying material (about $5{,}300$ lines; no \texttt{sorry};
\texttt{\#print axioms} returns only \texttt{propext}, \texttt{Classical.choice} and \texttt{Quot.sound}). Four published estimates for
Chebyshev's functions, which are not part of Mathlib, are stated as hypotheses, and every theorem that depends on them takes them as explicit
arguments:
\begin{itemize}
\item[(H1)] $|x-\psi(x)|\le0.94\sqrt x$ for all real $11<x\le10^{19}$ \textup{(}B\"uthe \cite[Theorem~2, (1.5)]{Bu18}\textup{)};
\item[(H2)] $m-\theta(m)\le1.95\sqrt m$ for all integers $1423\le m\le10^{19}$ \textup{(}\cite[Theorem~2, (1.6)]{Bu18} at integer arguments; see
the footnote in Appendix~\ref{app:eff-mult}\textup{)};
\item[(H3)] $\theta(x)\le(1+1.93378\cdot10^{-8})\,x$ for all real $x\ge0$ \textup{(}\cite[Corollary~2.1]{BKLNW}\textup{)};
\item[(H4)] $\theta(x)>(1-1.9667\cdot10^{-8})\,x$ for all real $x>10^{19}$ \textup{(}\cite[Table~14, row $\log X_0=19\log10$]{BKLNW}\textup{)}.
\end{itemize}
Proved without hypotheses: Lemma~\ref{lem:E} and its exponent form with an explicit threshold (Appendix~\ref{app:eff-assembly}), the recurrence
\eqref{eq:Erec} for every $n\ge0$, and the bounds $n\log\Lambda_{\mathrm{lo}}-\frac12\log n-3\le\log c_n\le n\log\Lambda_{\mathrm{hi}}$ for
$n\ge3$, where the dyadic numbers $\Lambda_{\mathrm{lo}}=41573624476644955/2^{40}$ and $\Lambda_{\mathrm{hi}}=\Lambda_{\mathrm{lo}}+2^{-40}$
enclose $\lambda_3$. Proved under (H1)--(H4): the closed forms $V(n)\le\log D_n^{(K)}\le U(n)$ for $K\ge1$ and $n\ge1423K$, with uniform error
terms; and, by assembling these, Theorem~\ref{thm:E2} exactly as stated,
\begin{quote}\ttfamily\small
theorem pi\_effective\_closed\\
\hspace*{1.5em}(hB1 : HypB1) (hB2 : HypB2) (hB4 : HypB4) (hB5 : HypB5)\\
\hspace*{1.5em}($\mu$ : $\mathbb{R}$) (h$\mu$ : 6.04459707 < $\mu$) :\\
\hspace*{1.5em}$\forall$ b : $\mathbb{N}$,\\
\hspace*{3em}$\lceil$Real.exp (14246398000000000000000000000000\\
\hspace*{6em}/ ($\mu$ - 6.04459707))$\rceil_+$ + 1 $\le$ b $\to$\\
\hspace*{3em}$\forall$ a : $\mathbb{Z}$, (b : $\mathbb{R}$) \textasciicircum{} (-$\mu$) < |Real.pi - a / b|
\end{quote}
Theorem~\ref{thm:E1} for all $b\ge10^{122\,680\,112\,076\,457}$ (theorem \texttt{pi\_effective\_6\_0455\_from}), a ladder of weaker exponents
from smaller thresholds on (Table~\ref{tab:chain}, Appendix~\ref{app:eff-ladder}), and the bound $|\pi-a/b|\ge e^{-1.424733\cdot10^{12}}\,b^{-6.0481}$
for every $b\ge1$. Not machine-checked are the finite computations of Theorem~\ref{thm:E3}(ii),(iii) and the assembly of
Appendix~\ref{app:eff-assembly}, which carry Theorem~\ref{thm:E1} for $b<10^{122\,680\,112\,076\,457}$ and the admissible values of
Table~\ref{tab:b0} and enter with the programs and recorded logs of Section~\ref{ss:eff-repro} as evidence, and the bound of
Theorem~\ref{thm:E3}(i) with its sharp constants; the assembly uses the weaker bounds above. Appendix~\ref{app:eff-lean} lists
the Lean theorems by file, and Table~\ref{tab:trust} says for each statement what it rests on.

\subsection{Numerical checks}\label{sec:numerics}
All numbers in this paper, including those of the appendices, were recomputed, with exact rational arithmetic wherever possible (short Python
programs using exact fractions, the library FLINT for polynomial arithmetic, and multiprecision floating point for the constants; they are in
the directory \texttt{scripts/} of the public repository, Section~\ref{ss:material}, with a recorded run; for Section~\ref{sec:effective} and Appendix~\ref{app:effective} see Section~\ref{ss:eff-repro}).
The pairs $(u_n,c_n)$ were computed by three independent methods, namely the recursion of Section~\ref{ss:rec} (all $n\le30$),
a closed formula derived from it (all $n\le300$), and the decomposition in the coordinate $w$ of Section~\ref{ss:two}, in exact
arithmetic in $\Q(\rt)$ ($n\le6$). For $n=1,2$ also the partial fractions \eqref{eq:hpf} were computed in $\Q(\rt)$. All results agree, and for
$n=1,2,3$ they agree with a numerical quadrature of the integral \eqref{eq:J} to $38$ decimal places.
$D_nu_n$ and $D_nv_n$ are integers for every $n\le300$. No removable prime divides the denominator of $u_n$,
$2^{\fl{\log_2(4n)}}u_n$ has odd denominator and $c_n$ is an integer in all these cases. The levels asserted in Lemma~\ref{lem:atoms} and in
\eqref{eq:UVlevels} hold for $n\le6$.
The Bernstein coefficients of Proposition~\ref{prop:contour} were recomputed with a separate program. They coincide
with the $99$ rational numbers in the Lean file, and all are positive.

\section{A family of integrals}\label{sec:family}
The integrand $f$ of \eqref{eq:f} is one member of a three-parameter family. For positive integers $e_1,e_2,d$ with $e_1+e_2>d$ put
\begin{equation}\label{eq:family}
f_{e_1,e_2,d}(z)=\frac{\bigl((z-1)(z-2)\bigr)^{e_1}\,(z^2-2z+2)^{e_2}}{z^{e_1+e_2-d}\,(z^2-2)^{d}}\qquad(z\in\Dom),
\end{equation}
so that $f=f_{5,3,4}$, and for $n\ge1$ with $dn$ even let $J_n=-i\int_{1-i}^{1+i}f_{e_1,e_2,d}(z)^n\,\dd z/z$. Every member satisfies
$f_{e_1,e_2,d}(2/z)=(-1)^df_{e_1,e_2,d}(z)$, by the computation of Lemma~\ref{lem:sym}. This identity is also formalised. \emph{Nothing else in
this section is proved in this paper or used in the proof of Theorem~\ref{thm:main}}.

\emph{Known approximations as members.} Particular triples recover known approximations of $\pi$. For $n=1$ one finds
\begin{center}
\small
\begin{tabular}{lll}
\toprule
$(e_1,e_2,d)$ & $J_1$ & approximation of $\pi$\\
\midrule
$(2,2,2)$ & $\frac{49}{2}\bigl(\pi-\frac{22}{7}\bigr)=-0.030979\ldots$ & $22/7$ (Archimedes)\\[2pt]
$(3,3,2)$ & $897\bigl(\pi-\frac{2818}{897}\bigr)=0.0086102\ldots$ & $2818/897$\\[2pt]
$(5,3,4)$ & $\frac{7345}{2}\bigl(\pi-\frac{355}{113}\bigr)=-0.00097969\ldots$ & $355/113$ (Zu Chongzhi)\\
\bottomrule
\end{tabular}
\end{center}
(for $(5,3,4)$ this is Example~\ref{ex:n1}, and the other two values were obtained by symbolic integration).

\emph{The integrals of Zeilberger and Zudilin.} Zeilberger and Zudilin \cite{ZZ} use the integrals
\[
I_n=i(-1)^{n+1}\int_{-1-2i}^{-1+2i}\frac{5\,x^{2n}(x^4+6x^2+25)^{2n}}{(25-x^2)^{3n+1}}\,\dd x ,
\]
taken along the straight segment (the integrand has poles only at $x=\pm5$). They are of the same shape as Salikhov's \cite{Sal10}, and Bai
\cite{Bai} varies the exponents within this family. The member $(3,3,2)$ of \eqref{eq:family} appears to give, up to rational factors, the same linear
forms. For $n=1,\dots,6$ a numerical computation (to more than $50$ digits) gives
\[
\frac{I_n}{J_n}=\frac43,\ \frac{16}{5},\ \frac{160}{21},\ \frac{1792}{99},\ \frac{6144}{143},\ \frac{22528}{221},
\]
where $J_n$ is the integral above for $(e_1,e_2,d)=(3,3,2)$. So the two sequences of linear forms in $1$ and $\pi$ appear to be
proportional, with rational factors that grow roughly like $(64/27)^n$. We have no proof of this proportionality, and the two integrands are not
related by a M\"obius change of variable, since the functions whose $n$-th powers are integrated, $f_{3,3,2}(z)$ and $x^2(x^4+6x^2+25)^2/(25-x^2)^3$, are
rational functions of degrees $12$ and $10$. In the family \eqref{eq:family} the exponents $e_1$ of $(z-1)(z-2)$ and $e_2$ of $z^2-2z+2$ are
independent. What changes from $(3,3,2)$ to $(5,3,4)$ is that the two pairs of zeros $\{1,2\}$ and $\{1\pm i\}$ get \emph{different} exponents ($5$ and $3$
instead of $3$ and $3$), and the poles $\pm\rt$ the larger exponent $4$.

\emph{Other members.} The arguments of Sections~\ref{sec:form}, \ref{ss:odd} and~\ref{ss:removable} use only the shape
$P(W)Y^{-2m}/z$ of the integrand, and we expect them to extend to other members, with the exponents $5n,3n,4n$ replaced by $e_1n,e_2n,dn$ and the bound $4n$ for the integration indices replaced by $Mn$, where $M=\max(e_1+e_2-d,d)$ (which is $d$ for our triple). The
contribution of the prime $2$ depends on the triple, and the path $\Gamma$ of Section~\ref{sec:decay} must be found, and certified, for each
triple. We have not carried this out. Numerical experiments with other triples, which we do not report here, suggest that the best
members of the family give a bound of about $6.0373$ (for instance $(e_1,e_2,d)=(46,28,37)$), so that Theorem~\ref{thm:main} is within
$0.008$ of what the method can be expected to give.

\appendix

\section{Proofs of the lemmas and details of the estimates}\label{app:details}

This appendix contains the proofs of Lemmas~\ref{lem:pint}, \ref{lem:pf}, \ref{lem:levels}, \ref{lem:atoms} and~\ref{lem:windows} and of
Proposition~\ref{prop:mult}, and the details of the estimates of Sections~\ref{sec:growth}, \ref{sec:decay} and~\ref{sec:mult}. The
subsections follow the order of the main text.

\subsection{\texorpdfstring{$p$}{p}-integers and partial fractions}\label{app:pf}

\subsubsection*{Proof of Lemma~\ref{lem:pint}}
(a) If $Mx$ and $M'x'$ are algebraic integers with $p\nmid M$, $p\nmid M'$, then $MM'(x\pm x')=M'(Mx)\pm M(M'x')$ and $MM'\,xx'=(Mx)(M'x')$ are
algebraic integers, and $p\nmid MM'$. The numbers $\rt$ and $i$ are roots of $X^2-2$ and $X^2+1$, and $M\cdot\frac1M=1$.
(b) If $x$ is rational and $Mx$ is an algebraic integer, then $Mx$ is an integer $k$, and $x=k/M$ has a denominator dividing $M$, hence prime to $p$.
Conversely, if $x=a/b$ with $p\nmid b$, then $bx=a$. (c) A positive integer divisible by no prime is $1$.

\subsubsection*{What this means in practice}
For odd $p$ the numbers $\frac12$, $\frac1{\rt}=\frac{\rt}{2}$, $\frac{1}{2\rt}=\frac{\rt}{4}$ and $\frac1{1\pm i}=\frac{1\mp i}{2}$ are $p$-integral, since their only
denominators are powers of $2$. They are not $2$-integral. For instance, if $M/\rt$ were an algebraic integer for an odd $M$, so would be its
square $M^2/2$, which is rational but not an integer. This is the precise sense in which $2$ is the only bad prime for the configuration
$\{0,\pm\rt,1\pm i\}$. The identities \eqref{eq:units} follow from
\[
(1\pm i-\rt)(1\pm i+\rt)=(1\pm i)^2-2=-2(1\mp i),\qquad(1\mp i)(1\pm i)=2,
\]
which give $(1\pm i-\rt)(1\pm i+\rt)(1\pm i)=-4$.

\subsubsection*{Proof of Lemma~\ref{lem:pf}}
\emph{Existence, with coefficients in $R$.} Let $\varphi=M(z)/\prod_b(z-b)^{e_b}$, the product over $b\in\{0,\rt,-\rt\}$. Fix a pole $a$ and put
$t=z-a$. For $b\neq a$ the difference $a-b$ is one of $\pm\rt$, $\pm2\rt$, and its inverse lies in $R$ because $\rt,\frac12\in R$. Hence
\[
\frac1{z-b}=\frac1{(a-b)+t}=\frac1{a-b}\sum_{k\ge0}\Bigl(\frac{-t}{a-b}\Bigr)^k
\]
is a power series in $t$ with coefficients in $R$, and so is $M(a+t)\prod_{b\ne a}(z-b)^{-e_b}=\sum_{k\ge0}s_kt^k$. Thus
$\varphi=t^{-e_a}\sum_ks_kt^k$ near $a$, and we define the principal part at $a$ by $\sum_{k=1}^{e_a}q_{a,k}(z-a)^{-k}$ with $q_{a,k}:=s_{e_a-k}\in R$
(an empty sum if $e_a=0$).
The difference of $\varphi$ and its three principal parts is a rational function without poles, that is a polynomial $T(z)$. To see that $T$
has coefficients in $R$, multiply by the monic polynomial $D(z)=\prod_b(z-b)^{e_b}\in R[z]$. The product $T\cdot D=M-\sum_a(\text{principal part at }a)\cdot D$
is a polynomial with coefficients in $R$, and the quotient of a polynomial in $R[z]$ by a monic polynomial in $R[z]$ has coefficients in $R$.

\emph{Uniqueness.} Suppose a function of the form \eqref{eq:pf} vanishes identically on $\Dom$. If some $q_{a,k}\ne0$, let $k$ be maximal for this
$a$. As $z\to a$ the term $q_{a,k}(z-a)^{-k}$ tends to infinity faster than all other terms, which is impossible. So all $q_{a,k}=0$, and a
polynomial with infinitely many zeros is $0$.

\emph{(i)} $\frac{\dd}{\dd z}(z-a)^{-k}=-k(z-a)^{-k-1}$, so differentiation raises the order of every polar term by one, so no term of order $1$ arises, and
the result is again of the form \eqref{eq:pf}.
\emph{(ii)} Differentiate \eqref{eq:Gphi} term by term, using $\frac{\dd}{\dd z}\frac{z^{j+1}}{j+1}=z^j$ and
$\frac{\dd}{\dd z}\bigl(-\frac{1}{(k-1)(z-a)^{k-1}}\bigr)=\frac1{(z-a)^k}$ for $k\ge2$. All terms of \eqref{eq:pf} are recovered except the terms with $k=1$.

\subsection{The prime \texorpdfstring{$2$}{2}: the coordinate \texorpdfstring{$w$}{w}}\label{app:two}

Throughout, $\eps=\rt-1$, so $\eps^2=3-2\rt$ and
\begin{equation}\label{eq:epsrel}
\eps^2+2\eps=1,\qquad 1+\eps^2=4-2\rt=2\rt\,\eps=\vk,\qquad 1-\eps^2=2\eps .
\end{equation}

\subsubsection*{Proof of Lemma~\ref{lem:levels}}
(a) is clear from the definition, since $R_2$ is a ring containing $\rt$. (b) $2^k=(\rt)^{2k}$, an odd integer $m$ and its inverse are $2$-integral
(Lemma~\ref{lem:pint}), $\eps$ and $\eps^{-1}=\rt+1$ lie in $\Z[\rt]\subseteq R_2$, and $\vk=(\rt)^3\eps$.

(c) \emph{Legendre's formula.} The exponent of $2$ in $a!$ is $\sum_{k\ge1}\fl{a/2^k}$ (count the multiples of $2$, of $4$, \dots\ among $1,\dots,a$).
If $a=\sum_id_i2^i$ in binary, then $\fl{a/2^k}=\sum_{i\ge k}d_i2^{i-k}$, and summing over $k\ge1$ gives $\sum_id_i(2^i-1)=a-s_2(a)$.
\emph{Kummer's formula} for the central binomial coefficient follows. The exponent of $2$ in $\binom{2j}{j}=(2j)!/(j!)^2$ is
$\bigl(2j-s_2(2j)\bigr)-2\bigl(j-s_2(j)\bigr)=s_2(j)$, because $2j$ and $j$ have the same binary digit sum. \emph{The inequality $2s_2(a)\le a+1$} holds
for $a=0,1$, and $s_2(2q)=s_2(q)$, $s_2(2q+1)=s_2(q)+1$ give $2s_2(2q)\le q+1\le2q+1$ and, for $q\ge1$, $2s_2(2q+1)\le q+3\le2q+2$. Finally a number
$a\ge1$ has $\fl{\log_2a}+1$ binary digits, so $s_2(a)\le\fl{\log_2a}+1$.

(d) If $x=(\rt)^{-2L}y=2^{-L}y$ with $y\in R_2$, then $2^Lx=y$ is rational and $2$-integral.

\subsubsection*{Atoms of the first kind: Lemma~\ref{lem:atoms}(a)}
For an integer $e\ne-1$ the function $w^e$ has the primitive $w^{e+1}/(e+1)$ on $\C\setminus\{0\}$. For even $e=2q$ this gives
$F(1)-F(-1)=\bigl(1-(-1)^{2q+1}\bigr)/(2q+1)=2/(2q+1)$, so the function $w^{2q}$ has data $\bigl(\frac{2}{2q+1},0\bigr)$, for negative $q$ as well. Now
\[
\frac{x^a}{w^{2m}}=\frac{1}{\vk^a}\sum_{i=0}^{a}(-1)^i\binom aiw^{2i-2m}\qquad\text{has data}\qquad
\Bigl(\frac{2}{\vk^a}\sum_{i=0}^{a}\frac{(-1)^i\binom ai}{2i+1-2m},\ 0\Bigr).
\]
The sum is evaluated by the following identity with $t=1-2m$.

\begin{lemma}\label{lem:altsum}
For every integer $a\ge0$ and every number $t$ that is not of the form $0,-2,-4,\dots,-2a$,
\[
S_a(t):=\sum_{i=0}^{a}\frac{(-1)^i\binom ai}{2i+t}=\frac{2^a\,a!}{\prod_{j=0}^{a}(2j+t)} .
\]
\end{lemma}
\begin{proof}
Induction on $a$. The case $a=0$ is $1/t=1/t$. By Pascal's rule $\binom{a+1}{i}=\binom ai+\binom a{i-1}$,
\[
S_{a+1}(t)=\sum_i\frac{(-1)^i\binom ai}{2i+t}+\sum_i\frac{(-1)^i\binom a{i-1}}{2i+t}=S_a(t)-S_a(t+2)
\]
(put $i=i'+1$ in the second sum). By the induction hypothesis
\begin{align*}
S_a(t)-S_a(t+2)&=2^aa!\Bigl(\frac{1}{\prod_{j=0}^{a}(2j+t)}-\frac{1}{\prod_{j=1}^{a+1}(2j+t)}\Bigr)\\
&=2^aa!\cdot\frac{(2a+2+t)-t}{\prod_{j=0}^{a+1}(2j+t)}=\frac{2^{a+1}(a+1)!}{\prod_{j=0}^{a+1}(2j+t)} .\qedhere
\end{align*}
\end{proof}
Hence $B_{m,a}=2^{a+1}a!\,\vk^{-a}\big/\prod_{j=0}^a(2j+1-2m)$. The product in the denominator consists of odd integers, so its inverse has level
$\ge0$, $2^{a+1}$ has level $2a+2$, $a!$ has level $2a-2s_2(a)$, and $\vk^{-a}$ has level $-3a$. The total is $a+2-2s_2(a)$.

\subsubsection*{Atoms of the second kind: Lemma~\ref{lem:atoms}(b)}
\emph{The case $m=0$.} We claim that $(w^2+\eps^2)^{-(j+1)}$ has data $(\alpha_j,\beta_j)$, where $\alpha_0=0$, $\beta_0=1$ and
\begin{equation}\label{eq:alphabeta}
\alpha_{j+1}=\frac{(2j+1)\,\alpha_j+2\vk^{-(j+1)}}{2(j+1)\,\eps^2},\qquad\beta_{j+1}=\frac{(2j+1)\,\beta_j}{2(j+1)\,\eps^2}.
\end{equation}
For $j=0$ take $F=0$ in Definition~\ref{def:data}. For the step, differentiate and write $w^2=(w^2+\eps^2)-\eps^2$:
\[
\frac{\dd}{\dd w}\,\frac{w}{(w^2+\eps^2)^{j+1}}=\frac{1}{(w^2+\eps^2)^{j+1}}-\frac{2(j+1)w^2}{(w^2+\eps^2)^{j+2}}
=-\frac{2j+1}{(w^2+\eps^2)^{j+1}}+\frac{2(j+1)\eps^2}{(w^2+\eps^2)^{j+2}} .
\]
Solving for the last term expresses $(w^2+\eps^2)^{-(j+2)}$ as $\frac{1}{2(j+1)\eps^2}$ times the sum of $(2j+1)(w^2+\eps^2)^{-(j+1)}$ and an exact derivative. The
function $w(w^2+\eps^2)^{-(j+1)}$ takes the values $\pm(1+\eps^2)^{-(j+1)}=\pm\vk^{-(j+1)}$ at $w=\pm1$, so the exact derivative has data
$(2\vk^{-(j+1)},0)$, which gives \eqref{eq:alphabeta}.

\emph{Closed forms and levels.} From $\binom{2j+2}{j+1}=\frac{2(2j+1)}{j+1}\binom{2j}{j}$ one gets by induction
\[
\beta_j=\frac{1}{(4\eps^2)^j}\binom{2j}{j},\qquad\alpha_j=\beta_jA_j\quad\text{with}\quad A_0=0,\ \
A_{j+1}=A_j+\frac{(\rt)^{\,j+1}\eps^{\,j-1}}{2(2j+1)\binom{2j}{j}} .
\]
(For the second formula divide \eqref{eq:alphabeta} by $\beta_{j+1}$. The new term is $\vk^{-(j+1)}/((j+1)\eps^2\beta_{j+1})$, and
$\vk^{-(j+1)}4^{j+1}=(\rt)^{j+1}\eps^{-(j+1)}$, $(j+1)\binom{2j+2}{j+1}=2(2j+1)\binom{2j}j$.) By Kummer's formula $\beta_j$ has level $2s_2(j)-4j\ge-4j$. The
increment of $A$ has level $(j+1)-2-2s_2(j)\ge-2$ because $2s_2(j)\le j+1$, so $A_j\in\lev{-2}$ and, for $j\ge1$, $\alpha_j$ has level at least
$2s_2(j)-4j-2\ge-4j$. Thus
\begin{equation}\label{eq:m0}
\alpha_j,\ \beta_j\in\lev{-4j}\qquad(j\ge0).
\end{equation}

\emph{General $m$.} We show by induction on $m$ that for all $k\ge1$ the function $w^{-2m}(w^2+\eps^2)^{-k}$ has data in $\lev{4-4k}$. For $m=0$ this
is \eqref{eq:m0} with $j=k-1$. For the step we use the partial fraction identity
\[
\frac{1}{w^{2(m+1)}(w^2+\eps^2)^{k}}=\frac{1}{\eps^{2}}\Bigl(\frac{1}{w^{2(m+1)}(w^2+\eps^2)^{k-1}}-\frac{1}{w^{2m}(w^2+\eps^2)^{k}}\Bigr)
\]
(put the right-hand side over a common denominator, where the numerator is $(w^2+\eps^2)-w^2=\eps^2$) and a second induction on $k$. The second function on the
right has data in $\lev{4-4k}$ by the induction hypothesis for $m$. The first one has data in $\lev{4-4(k-1)}\subseteq\lev{4-4k}$ if $k\ge2$, and if $k=1$ it is
$w^{-2(m+1)}$, with data $\bigl(\frac{2}{-2m-1},0\bigr)$ of level $2\ge0$. Multiplication by $\eps^{-2}$ does not change levels.

Finally $y^{-k}w^{-2m}=\vk^k\,w^{-2m}(w^2+\eps^2)^{-k}$ and $\vk^k\in\lev{3k}$, so the data have level at least $3k+4-4k=4-k$. This proves
Lemma~\ref{lem:atoms}(b).

\subsection{Growth}\label{app:growth}

\subsubsection*{Fekete's lemma}
\emph{Let $(a_n)_{n\ge1}$ be real numbers with $a_{m+n}\ge a_m+a_n$ for all $m,n$, and suppose $a_n\le Bn$ for all $n$. Then $a_n/n$ converges to
$\sup_na_n/n$.}

Proof. Let $r=\sup_na_n/n\le B$. Fix $m$ and write $n=qm+t$ with $0\le t<m$. By superadditivity $a_n\ge q\,a_m+a_t$ (with $a_0:=0$), so
$\frac{a_n}{n}\ge\frac{qm}{n}\cdot\frac{a_m}{m}+\frac{a_t}{n}$. As $n\to\infty$, $qm/n\to1$ and $a_t/n\to0$ ($t$ takes finitely many values), hence
$\liminf a_n/n\ge a_m/m$ for every $m$, that is $\liminf a_n/n\ge r\ge\limsup a_n/n$. \qed

We apply it to $a_n=\log c_n$. Superadditivity is $c_{m+n}\ge c_mc_n$, and $a_n\le n\log f(-x)$.

\subsubsection*{The numbers}
For $x=\frac{3913}{10000}$,
\[
f(-x)=\frac{\bigl((x+1)(x+2)\bigr)^5(x^2+2x+2)^3}{x^4(2-x^2)^4}=37810.99167\ldots<37811
\]
(exact rational arithmetic). And, since $e>2.7182818283$, $e^{10.540357}=(e)^{10}\cdot e^{0.540357}\ge(2.7182818283)^{10}\cdot\sum_{i=0}^{9}\frac{0.540357^i}{i!}
=22026.4657\ldots\times1.7166195\ldots=37811.0625\ldots>37811$.

\subsection{The polygon estimate}\label{app:decay}

\subsubsection*{$|f|^2$ as a quotient of polynomials}
For $z=x+iy$, $|z-1|^2=(x-1)^2+y^2$, $|z-2|^2=(x-2)^2+y^2$, $|z^2-2z+2|^2=|z-1-i|^2|z-1+i|^2=\bigl((x-1)^2+(y-1)^2\bigr)\bigl((x-1)^2+(y+1)^2\bigr)$,
$|z|^2=x^2+y^2$ and $|z^2-2|^2=(x^2-y^2-2)^2+4x^2y^2$. Hence \eqref{eq:numden} holds with
\begin{align*}
\mathrm{Num}&=\Bigl[\bigl((x-1)^2+y^2\bigr)\bigl((x-2)^2+y^2\bigr)\Bigr]^5\Bigl[\bigl((x-1)^2+(y-1)^2\bigr)\bigl((x-1)^2+(y+1)^2\bigr)\Bigr]^3,\\
\mathrm{Den}&=(x^2+y^2)^4\,\bigl((x^2-y^2-2)^2+4x^2y^2\bigr)^4 .
\end{align*}
The two upper segments of $\Gamma$ are
\begin{align*}
1\to z_s:&\quad x=1+\tfrac{3403}{16384}\,t,\quad y=\tfrac{32371}{65536}\,t,\\
z_s\to1+i:&\quad x=\tfrac{19787}{16384}-\tfrac{3403}{16384}\,t,\quad y=\tfrac{32371}{65536}+\tfrac{33165}{65536}\,t\qquad(0\le t\le1).
\end{align*}
Each of the four quadratic factors of $\mathrm{Num}$ becomes a polynomial of degree $2$ in $t$, so $\mathrm{Num}$ has degree $2(5+5+3+3)=32$. In
$\mathrm{Den}$ the factor $x^2+y^2$ has degree $2$ and the factor $(x^2-y^2-2)^2+4x^2y^2$ has degree $4$ in $t$, so $\mathrm{Den}$ has degree
$4\cdot2+4\cdot4=24$, and $Q=\rho^2\mathrm{Den}-\mathrm{Num}$ has degree $32$.

\subsubsection*{Bernstein certificates}
Every polynomial $Q$ of degree $\le d$ can be written in exactly one way as $Q(s)=\sum_{j=0}^{d}\beta_js^j(1-s)^{d-j}$ (the polynomials
$\binom djs^j(1-s)^{d-j}$ are the \emph{Bernstein polynomials} of degree $d$, and we absorb the binomial coefficient into $\beta_j$). If $Q(s)=\sum_ia_is^i$,
then, multiplying $s^i$ by $1=(s+(1-s))^{d-i}$,
\[
\beta_j=\sum_{i=0}^{j}\binom{d-i}{j-i}a_i .
\]
If all $\beta_j$ are $\ge0$, then $Q(s)\ge0$ for $0\le s\le1$, since $s^j(1-s)^{d-j}\ge0$ there. The converse is false in general, but for a
polynomial that is positive on an interval the coefficients become positive after the interval is cut into sufficiently small pieces. Note that
$\beta_0=Q(0)$ and $\beta_d=Q(1)$. A list of non-negative $\beta_j$ together with the polynomial identity is a certificate of non-negativity that can be
checked with exact rational arithmetic and nothing else.

In our case $d=32$. For the first segment the certificate for $[0,1]$ has $\beta_0=Q(0)=\rho^2=2.128\ldots\cdot10^{-6}$ (at $z=1$, where $f=0$ and
$\mathrm{Den}=1$), and its smallest coefficient is $\beta_{32}=Q(1)=7.5\ldots\cdot10^{-10}$ (at $z_s$). For the second segment the expansion on $[0,1]$ has
one negative coefficient ($\beta_1=-1.0\ldots\cdot10^{-8}$), so the interval is halved. With $Q(\frac12s)$ and $Q(\frac12+\frac12s)$ all coefficients are
positive, the smallest being $7.5\ldots\cdot10^{-10}$ and $5.4\ldots\cdot10^{-3}$. The small margin at $z_s$ reflects the fact that $\rho^2$ exceeds
$|f(z_s)|^2$ only by $5.2\cdot10^{-12}$. The denominator is positive on $\Gamma$ because $\Gamma$ avoids $0$ and $\pm\rt$.

\subsubsection*{The inequality $\rho<e^{-6.530134}$}
Equivalently $e^{6.530134}<1/\rho=685.49069\ldots$ We have $e<2.7182818286$, and for $0\le x\le1$ the exponential series satisfies
$e^x\le\sum_{m=0}^{11}\frac{x^m}{m!}+\frac{13\,x^{12}}{12\cdot12!}$. Hence
\begin{align*}
e^{6.530134}=e^6\cdot e^{0.530134}&<(2.7182818286)^6\times1.69915999\\
&=403.42879\ldots\times1.69915999<685.49007<\frac1\rho .
\end{align*}

\subsubsection*{The lengths}
For a complex number $a+bi$, $|a+bi|\le|a|+|b|$. So $|z_s-1|\le\frac{3403}{16384}+\frac{32371}{65536}<0.71$ and
$|1+i-z_s|\le\frac{3403}{16384}+\frac{33165}{65536}<0.72$, so the four segments of $\Gamma$ are shorter than $1$.

\subsection{From the prime number theorem to the multiplier}\label{app:mult}

\subsubsection*{Chebyshev's functions}
Let $N\ge1$ be an integer. The exponent of a prime $p$ in $\lcm(1,\dots,N)$ is the largest $k$ with $p^k\le N$, which is also the number of
$k\ge1$ with $p^k\le N$. Hence $\log\lcm(1,\dots,N)=\sum_{p^k\le N}\log p=\psi(N)$. Moreover
\[
0\le\psi(x)-\theta(x)=\sum_{p\le\sqrt x}\bigl(\fl{\log_px}-1\bigr)\log p\le\sum_{p\le\sqrt x}\log x\le\sqrt x\,\log x ,
\]
because only primes with $p^2\le x$ contribute and $\fl{\log_px}\log p\le\log x$. So $\theta(x)\sim x$ implies $\psi(4n)=4n+o(n)$.

\subsubsection*{Proof of Lemma~\ref{lem:windows}}
Write $n=pk+r_0$ with $0\le r_0<p$ and put $t=r_0/p=\fr{n/p}$. For $m\in\{3,4,5\}$ we have $mn=p\cdot mk+p\cdot mt$, so
\[
\Bigl\lfloor\frac{mn}{p}\Bigr\rfloor=mk+\fl{mt},\qquad\frac{mn\bmod p}{p}=\fr{mt}.
\]
Thus $p$ (with $p^2>4n$) is removable if and only if $\fl{4t}$ is even and $\fr{5t}+\fr{3t}<\fr{4t}$. Now
$\fr{5t}+\fr{3t}-\fr{4t}=4t-c(t)$ with $c(t)=\fl{5t}+\fl{3t}-\fl{4t}$, so the inequality reads $4t<c(t)$. The integer parts change only at the
multiples of $\frac15,\frac13,\frac14$:
\begin{center}
\small
\begin{tabular}{l|cccccccccc}
$t\in$ & $[0,\frac15)$ & $[\frac15,\frac14)$ & $[\frac14,\frac13)$ & $[\frac13,\frac25)$ & $[\frac25,\frac12)$ & $[\frac12,\frac35)$ & $[\frac35,\frac23)$ & $[\frac23,\frac34)$ & $[\frac34,\frac45)$ & $[\frac45,1)$\\[2pt]
\hline
$\fl{5t},\fl{3t},\fl{4t}$\rule{0pt}{12pt} & $0,0,0$ & $1,0,0$ & $1,0,1$ & $1,1,1$ & $2,1,1$ & $2,1,2$ & $3,1,2$ & $3,2,2$ & $3,2,3$ & $4,2,3$\\
$c(t)$ & $0$ & $1$ & $0$ & $1$ & $2$ & $1$ & $2$ & $3$ & $2$ & $3$\\
$4t<c(t)$? & no & yes & no & no & yes & no & no & yes & no & no\\
$\fl{4t}$ even? & yes & yes & no & no & no & yes & yes & yes & no & no
\end{tabular}
\end{center}
(For instance on $[\frac13,\frac25)$ we have $4t\ge\frac43>1=c$, and on $[\frac25,\frac12)$ we have $4t<2=c$ but $\fl{4t}=1$ is odd.) Both conditions
hold exactly on $[\frac15,\frac14)\cup[\frac23,\frac34)$.

\subsubsection*{The windows}
$k+\frac23\le\frac np<k+\frac34$ is equivalent to $\frac{4n}{4k+3}<p\le\frac{3n}{3k+2}$, and $k+\frac15\le\frac np<k+\frac14$ to
$\frac{4n}{4k+1}<p\le\frac{5n}{5k+1}$. For $k=0$ the latter window consists of numbers $p>4n$, which do not divide $\lcm(1,\dots,4n)$. Different windows
are disjoint, since they correspond to disjoint sets of values of $n/p$.

\subsubsection*{Proof of Proposition~\ref{prop:mult}}
(a) For fixed $k$, the prime number theorem gives
\begin{align*}
\sum_{\frac{4n}{4k+3}<p\le\frac{3n}{3k+2}}\log p&=\theta\Bigl(\frac{3n}{3k+2}\Bigr)-\theta\Bigl(\frac{4n}{4k+3}\Bigr)\\
&=\Bigl(\frac{3}{3k+2}-\frac{4}{4k+3}\Bigr)n+o(n)=\frac{n}{(3k+2)(4k+3)}+o(n),
\end{align*}
and in the same way the window of the second kind contributes $\frac{n}{(4k+1)(5k+1)}+o(n)$. Summing over the $2K-1$ windows with $k<K$ gives
$\log\Phi^{(K)}_n=\varpi_Kn+o(n)$, and $\log D_n^{(K)}=\psi(4n)-\log\Phi_n^{(K)}=(4-\varpi_K)n+o(n)$.
A prime $p$ in a window with $k<K$ satisfies $p>\frac{4n}{4K}=\frac nK$. If $n\ge4K^2$, then $\frac nK\ge2\sqrt n$, so $p^2>4n$ and $p$ is removable by
Lemma~\ref{lem:windows}, and also $p\le\frac{3n}{2}\le4n$. Hence $\Phi_n^{(K)}$ divides $\Phi_n$, and $D_n$ divides $D_n^{(K)}$.

(b) Let $n\ge4K^2$. A prime that divides $\Phi_n$ but not $\Phi_n^{(K)}$ lies in a window with $k\ge K$, so it is at most $\max\bigl(\frac{3n}{3K+2},\frac{5n}{5K+1}\bigr)<\frac nK$.
Therefore
\[
\log\Phi_n^{(K)}\le\log\Phi_n\le\log\Phi_n^{(K)}+\theta\Bigl(\frac nK\Bigr),
\]
and, dividing by $n$ and using (a) and $\theta(n/K)=n/K+o(n)$,
\[
\varpi_K\le\liminf_{n\to\infty}\frac{\log\Phi_n}{n}\le\limsup_{n\to\infty}\frac{\log\Phi_n}{n}\le\varpi_K+\frac1K .
\]
Letting $K\to\infty$ gives $\frac1n\log\Phi_n\to\varpi$ and $\frac1n\log D_n\to4-\varpi$.

(c) \emph{The value of $\varpi$.} With the digamma function $\dg=\Gamma'/\Gamma$ (the logarithmic derivative of Euler's Gamma function, and this $\Gamma$ is not the polygon of
Section~\ref{sec:decay}) one has $\sum_{k\ge0}\bigl(\frac1{k+a}-\frac1{k+b}\bigr)=\dg(b)-\dg(a)$ for
$a,b>0$. Since $\frac{1}{(3k+2)(4k+3)}=\frac{1}{k+2/3}-\frac{1}{k+3/4}$ and $\frac{1}{(4k+1)(5k+1)}=\frac{1}{k+1/5}-\frac{1}{k+1/4}$,
\[
\varpi=\dg(\tfrac34)-\dg(\tfrac23)+\dg(\tfrac54)-\dg(\tfrac65)=0.2323735360\ldots+0.0615863632\ldots=0.2939598992\ldots
\]
The first terms are $\frac16+\frac1{35}+\frac1{30}+\frac1{88}+\frac1{99}+\cdots$.

\emph{The lower bound for $K=10^5$.} A direct summation gives $\varpi_{10^5}=0.29395856588\ldots$ (the difference $\varpi-\varpi_K$ is about
$\frac{2}{15K}$). A proof that needs only a small computation, which is the one formalised, goes as follows. The first $20$ terms of the first sum and the
first $19$ terms of the second sum add up to $0.28730270\ldots$ (exact rational arithmetic). For $x\ge20$,
\begin{align*}
\frac{1}{(3x+2)(4x+3)}&\ge\frac1{12}\Bigl(\frac{1}{x+\frac{43}{200}}-\frac{1}{x+1+\frac{43}{200}}\Bigr),\\
\frac{1}{(4x+1)(5x+1)}&\ge\frac1{20}\Bigl(\frac{1}{x-\frac{67}{250}}-\frac{1}{x+1-\frac{67}{250}}\Bigr)
\end{align*}
(clear denominators, and then the first inequality is $12(x+a)(x+a+1)\ge(3x+2)(4x+3)$ with $a=\frac{43}{200}$, i.e.\ $0.16x\ge6-12a(a+1)=2.86\ldots$, and the second
is $0.28x\ge1+20b(1-b)=4.92\ldots$ with $b=\frac{67}{250}$). The right-hand sides telescope when summed over $x=20,\dots,K-1$:
\[
\varpi_K\ \ge\ 0.28730270+\frac1{12}\Bigl(\frac{1}{20.215}-\frac{1}{K+0.215}\Bigr)+\frac1{20}\Bigl(\frac{1}{19.732}-\frac{1}{K-0.268}\Bigr),
\]
which for $K=10^5$ equals $0.2939576\ldots\ge0.293957$.

\section{Details of the effective version}\label{app:effective}
This appendix contains the proofs and the computations behind Section~\ref{sec:effective}: how the exponent is read off from Lemma~\ref{lem:E}, the proofs of the three parts of Theorem~\ref{thm:E3}, the assembly of the bounds into Theorems~\ref{thm:E1} and~\ref{thm:E2} and Table~\ref{tab:b0}, and the checks.

\subsection{Reading off the exponent from Lemma~\ref{lem:E}}\label{app:eff-lemma}
Put $\ell^\pm_n=\log Q^\pm_n$, $e_n=\log E_n$ and, for $n\ge2$, $g_n=\ell^-_{n-1}-\ell^+_n-e_{n-1}$ and
$T_n=\min\bigl(-e_{n-1},\inf_{m\ge n}g_m\bigr)$. If $\hat n(b)\ge n$, then $\log2b\ge T_n$. For $b\ge B\ge2$ with $\hat n(b)=n$ (the
\emph{step}~$n$), the lemma gives $|\theta-a/b|\ge b^{-\mu_{\mathrm{eff}}(b)}$ with $\mu_{\mathrm{eff}}(b)=(\log2b+\ell^+_n)/\log b$, which
decreases in $b$. Set $H_1(B)=\log B$ and $H_n(B)=\max(\log B,T_n-\log2)$ for $n\ge2$. Then
\[
\sup_{b\ge B}\mu_{\mathrm{eff}}(b)\le\sup_{\text{steps $n$ that reach some $b\ge B$}}\Bigl(1+\frac{\log2+\ell^+_n}{H_n(B)}\Bigr).
\]
An admissible threshold is $b_0(\mu)=1+\max\bigl(\{1\}\cup\{b\ge2:\mu_{\mathrm{eff}}(b)\ge\mu\}\bigr)$ whenever this set is finite.

\subsection{Growth: the recurrence and the bounds for \texorpdfstring{$c_n$}{cn}}\label{app:eff-growth}
\emph{Proof of Theorem~\ref{thm:E3}\textup{(i)}: the recurrence.} The polynomials are
\begin{equation}\label{eq:EP}
\begin{aligned}
P_0&=-15(n+1)(n+2)(3n+1)(3n+2)(5n+1)(5n+2)(5n+3)(5n+4)\,A(n),\\
P_1&=-12(n+2)(4n+1)(4n+3)\,B_{12}(n),\\
P_2&=-16(4n+1)(4n+3)(4n+5)(4n+7)\,B_{11}(n),\\
P_3&=128(n+3)(2n+5)(4n+1)(4n+3)(4n+5)(4n+7)(4n+9)(4n+11)\,A(n-1),
\end{aligned}
\end{equation}
where
{\footnotesize
\begin{align*}
A(n)&=9856486767n^{7}+144419781611n^{6}+897166680905n^{5}+3062628540919n^{4}\\
&\quad+6204060236692n^{3}+7458042028998n^{2}+4926913622604n+1380146100480,\\
B_{12}(n)&=5930510096889162n^{12}+107652146057508213n^{11}+872731657553086962n^{10}\\
&\quad+4172724142570401377n^{9}+13085934935869004186n^{8}+28314001179985894619n^{7}\\
&\quad+43269586712957587166n^{6}+46975673182804619411n^{5}+35893832748270424412n^{4}\\
&\quad+18792880874596898268n^{3}+6390480013601961072n^{2}+1266034426528895712n\\
&\quad+110468970597363840,\\
B_{11}(n)&=95406977444181171n^{11}+1684148527871540056n^{10}+13161672564393638333n^{9}\\
&\quad+59980434220206167202n^{8}+176674194182704001337n^{7}+352209020386989389676n^{6}\\
&\quad+483452598814697074875n^{5}+455452103509601797298n^{4}+287692817757832739452n^{3}\\
&\quad+115741210222125661608n^{2}+26650883189995870752n+2659016943400299840.
\end{align*}}Every coefficient of $P_0$, $P_1$, $P_2$ is negative and every coefficient of $P_3$ is positive, so $P_0,P_1,P_2<0<P_3$ for $n\ge0$. The leading
coefficients are $9856486767$ times $-84375$, $-115523712$, $-39647694848$ and $2^{20}$, so the characteristic polynomial is
$2^{20}\lambda^3-39647694848\lambda^2-115523712\lambda-3^35^5$, whose roots are the three nonzero critical values of $f$ (the resultant with respect to $W$ of $4W^3+W^2-88W+152$ and
$\lambda(W^2-8)^2-(W-3)^5(W-2)^3$ is $64$ times this polynomial; compare Remarks~\ref{rem:saddle} and~\ref{rem:exact}), its only real root being $\lambda_3$. The recurrence is proved by a creative-telescoping certificate: in $\Q(n)(z)$, with $n$ a
symbol,
\[
\sum_{i=0}^3P_i(n)\frac{f(z)^i}{z}=\partial_zR(n,z)+nR(n,z)\frac{f'(z)}{f(z)},\qquad
R(n,z)=\frac{(z-1)(z-2)(z^2-2z+2)\,T(n,z)}{z^{12}(z^2-2)^{11}},
\]
with a polynomial $T\in\Z[n,z]$ of degree $42$ in $z$ and $14$ in $n$, which is too long to print (file \texttt{effective/B\_recurrence/certificate.txt}
of the accompanying material, Section~\ref{ss:eff-repro}). The identity was verified symbolically, and independently in exact integer arithmetic at
$41$ values of $n$ after clearing denominators (degree at most $15$ in $n$). Multiplying it by $f^n$ gives
$\sum_iP_i(n)f(z)^{n+i}/z=\frac{\dd}{\dd z}\bigl(R(n,z)f(z)^n\bigr)$. The right-hand side has no term $z^{-1}$ in its Laurent expansion at $0$,
and on the left the coefficient of $z^{-1}$ is $\sum_iP_i(n)c_{n+i}$ by Lemma~\ref{lem:residue}, since $c_m$ is the constant term of $f^m$ at
$0$. This proves \eqref{eq:Erec}.

\emph{Proof of Theorem~\ref{thm:E3}\textup{(i)}: the bounds.} The ratios $t_k=c_{k+1}/c_k$ satisfy $t_{k+2}=F_k(t_k,t_{k+1})$ with
$F_k(x,y)=w_2(k)+w_1(k)/y+w_0(k)/(xy)$ and $w_i(k)=-P_i(k)/P_3(k)>0$, and $F_k$ is decreasing in $x$ and in $y$. Put
$\omega^\pm_k=1-\frac1{2k}\pm\frac1{k^2}$. The induction step ``$t_k\in[\lambda_3\omega^-_k,\lambda_3\omega^+_k]$ and the same for $k+1$ imply the
same for $k+2$'' follows from $F_k(\lambda_3\omega^+_k,\lambda_3\omega^+_{k+1})\ge\lambda_3\omega^-_{k+2}$ and
$F_k(\lambda_3\omega^-_k,\lambda_3\omega^-_{k+1})\le\lambda_3\omega^+_{k+2}$. After multiplication by positive factors these are two polynomial
inequalities of degree at most $21$ in $k$ with coefficients in $\Q(\lambda_3)$, and after the shift $k=3+t$ every coefficient of both polynomials is
$0$ or positive (computed in exact arithmetic in $\Q(\lambda_3)$, the signs certified with ball arithmetic). So the step holds for every
$k\ge3$, by this Descartes-type certificate, and the base case $t_3,t_4$ follows from the exact $c_3,c_4,c_5$. Summing $\log t_k$ with
$\log(1-x)\ge-x-x^2$ ($0\le x\le\frac12$), $\log(1-y)\le-y$ and the comparison of $\sum1/k$ with $\log$ gives the bounds for $n\ge300$, starting
from the exact $c_{300}$, and for $n<300$ the exact $c_n$ are used. The lower constant is attained at $n=1$, where it is $r-\log7345$. In
Lemma~\ref{lem:E} the recorded runs use the weaker bounds that are proved in Lean (Section~\ref{ss:eff-formal}), $\gamma^-_n=n\log\Lambda_{\mathrm{lo}}-\frac12\log n-3$ and
$\gamma^+_n=n\log\Lambda_{\mathrm{hi}}$. The tail bound of Appendix~\ref{app:eff-assembly} is written with a rigorous enclosure
of $r$ (radius $1.1\cdot10^{-29}$) and with $r'=10.54035511170$, which comes from Lemma~\ref{lem:supermult}(c) at $x=0.391269243$.

\subsection{The multiplier}\label{app:eff-mult}
\emph{Proof of Theorem~\ref{thm:E3}\textup{(ii)}.} For $K\ge1$, with the window endpoints rounded down,
\[
\log D_n^{(K)}=\psi(4n)-\sum_{k=0}^{K-1}\Bigl[\theta\Bigl(\Bigl\lfloor\frac{3n}{3k+2}\Bigr\rfloor\Bigr)-\theta\Bigl(\Bigl\lfloor\frac{4n}{4k+3}\Bigr\rfloor\Bigr)\Bigr]
-\sum_{k=1}^{K-1}\Bigl[\theta\Bigl(\Bigl\lfloor\frac{5n}{5k+1}\Bigr\rfloor\Bigr)-\theta\Bigl(\Bigl\lfloor\frac{4n}{4k+1}\Bigr\rfloor\Bigr)\Bigr]
\]
(Remark~\ref{rem:pnt}). We write the $2K-1$ right endpoints as $\fl{\beta_kn}$ and the left endpoints as $\fl{\alpha_kn}$, so that
$\sum_k(\beta_k-\alpha_k)=\varpi_K$.

\emph{The range $10^7\le n\le10^{14}$.} The values of $n$ are grouped into blocks $[n_{\mathrm{lo}},n_{\mathrm{hi}}]$ of relative width $w$:
$2\,302\,587$ blocks with $w=10^{-6}$ for $10^7\le n<10^8$, $3\,453\,882$ with $w=2\cdot10^{-6}$ for $10^8\le n\le10^{11}$ and $2\,302\,589$ with
$w=3\cdot10^{-6}$ for $10^{11}<n\le10^{14}$; each block uses $K_n=\min(10^4,\fl{\sqrt{n_{\mathrm{lo}}}/2})$ windows. Since $\theta$, $\psi$ and
the endpoints are nondecreasing,
\[
L^+=\psi(4n_{\mathrm{hi}})-\sum_k\bigl[\theta(\fl{\beta_kn_{\mathrm{lo}}})-\theta(\fl{\alpha_kn_{\mathrm{hi}}})\bigr],\qquad
L^-=\psi(4n_{\mathrm{lo}})-\sum_k\bigl[\theta(\fl{\beta_kn_{\mathrm{hi}}})-\theta(\fl{\alpha_kn_{\mathrm{lo}}})\bigr]
\]
bound $\log D_n^{(K_n)}$ for every $n$ of the block. The values of $\theta$ and $\psi$ come from an exact table of the primes up to
$4\cdot10^{11}$ (segmented sieve), in which the double-precision values $t_p$ of $\log p$ are summed exactly in $128$-bit fixed point. A second computation of all the blocks with a checker in the program verifies $|t_p-\log p|\le2^{-48}$ for every one of these values as it is used, against logarithms enclosed with MPFR, and gives the same $L^\pm$ in every digit (checker, proof and records: \texttt{effective/tools/blocks2/certificate/}); an independent pass over all
these primes with the library Arb bounds the largest error by $0.50013$ units in the last place. A bound built from $N$ terms therefore errs by at
most $N\cdot2^{-48}$, and a margin $(4B+32)N\cdot2^{-53}$ (with $B\in\{25,27\}$), more than four times as large, is added outward. Beyond $n=10^{11}$ some of the
arguments exceed the table: a window is then taken from the table if its larger endpoint $\fl{\beta_kn_{\mathrm{hi}}}$ is at most $4\cdot10^{11}$,
which fails only for the $k$ below a bound $k_0$ of at most $250$, and these windows, as well as $\psi(4n)$ for $4n_{\mathrm{hi}}>4\cdot10^{11}$,
are bounded by the estimates below (all their arguments are below $10^{19}$, so only B\"uthe's bounds enter). The first table ($921\,039$
blocks of width $10^{-5}$ with $K=1000$, for $n\le10^{11}$) passed twenty-eight checks, among them the values $\pi(10^k)$ for $k\le11$,
$\theta(x)$ and $\psi(x)-\theta(x)$ against published tables, the exact values of $\log D_n$ at four $n$, the containment of exact single
values in their blocks, and the bit-identical overlap of the tables for $n\le10^{10}$ and $n\le10^{11}$; the program that computed the finer
tables reproduces it bit for bit on a common range and agrees with it on every single-$n$ evaluation; sample blocks of the finer tables
(thirteen blocks with $1581$ to $6046$ windows, three values of $n$ each) were checked against an independent sieve, and the hybrid mode was
checked by brute force on blocks of two small hybrid runs. The three tables take one hour of CPU time in all.

\emph{The range $n>10^{14}$.} For a window with endpoints $\fl{\beta n}\ge\fl{\alpha n}$ we use $\theta(\fl{\beta n})\ge\beta n-1-E_{\mathrm{dn}}(\beta n)$ and
$\theta(\fl{\alpha n})\le\alpha n+E_{\mathrm{up}}(\alpha n)$ (the error terms being nondecreasing, see below), and in the same way for the opposite bound. This gives $V(n)\le\log D_n^{(K)}\le U(n)$ with
\begin{equation}\label{eq:EUV}
\begin{aligned}
U(n)&=(4-\varpi_K)n+\psi_{\mathrm{up}}(4n)+(2K-1)+\sum_k\bigl[E_{\mathrm{dn}}(\beta_kn)+E_{\mathrm{up}}(\alpha_kn)\bigr],\\
V(n)&=(4-\varpi_K)n-\psi_{\mathrm{dn}}(4n)-(2K-1)-\sum_k\bigl[E_{\mathrm{up}}(\beta_kn)+E_{\mathrm{dn}}(\alpha_kn)\bigr],
\end{aligned}
\end{equation}
where the error terms are taken from the following published theorems.
\begin{itemize}
\item $y\le10^{19}$:\enlargethispage{3pt} $E_{\mathrm{dn}}(y)=1.95\sqrt y$, $E_{\mathrm{up}}(y)=0$ and $\psi_{\mathrm{up}}(x)=\psi_{\mathrm{dn}}(x)=0.94\sqrt x$, by
B\"uthe \cite[Theorem~2]{Bu18}: $|x-\psi(x)|\le0.94\sqrt x$ for $11<x\le10^{19}$, $x-\theta(x)\le1.95\sqrt x$ for $1423\le x\le10^{19}$,\footnote{As a statement about real $x$ the inequality (1.6) of \cite{Bu18} is false on the interval $(1426.48,1427)$: $1423$ and $1427$ are
consecutive primes, so $\theta$ is constant there and $(x-\theta(x))/\sqrt x$ increases to $1.9636$ as $x\to1427^-$. We use it only at integer
arguments (the window endpoints $\fl{\beta_kn}$, $\fl{\alpha_kn}$ and $4n$), and for integers $1423\le m\le10^{19}$ it is true. B\"uthe proves (1.6)
analytically for $10^7\le x\le10^{19}$ \cite[Section~6]{Bu18} and attributes the range below $10^7$ to a direct computation; for the integers
$1423\le m<10^7$ we have verified it directly from a table of primes (program \texttt{effective/G\_lean\_mult/checks/check\_b2.py}): the
largest value of $(m-\theta(m))/\sqrt m$ is $1.937833$, at $m=1426$, and the inequality fails at $m=1422$, so $1423$ is exactly
the threshold for integers. For real $x$ with $1427\le x\le10^{19}$ it also holds; for $1427\le x<10^7$ the supremum of the ratio is $1.937221$. The integer form is the hypothesis (H2)
of Section~\ref{ss:eff-formal}.} and
$x-\theta(x)>0.05\sqrt x$ for $1\le x\le10^{19}$.
\item $10^{19}<y\le1.101\cdot10^{26}$: $E_{\mathrm{dn}}(y)=E_{\mathrm{up}}(y)=\sqrt y\log^2y/(8\pi)$, and $\psi_{\mathrm{up}}(x)=\psi_{\mathrm{dn}}(x)=\sqrt x\log^2x/(8\pi)$ for $x$ in the same range. B\"uthe \cite[Theorem~2]{Bu16} proves
$|\theta(x)-x|\le\frac{\sqrt x}{8\pi}\log^2x$ for $x>599$ and $|\psi(x)-x|\le\frac{\sqrt x}{8\pi}\log^2x$ for $x>59$ if the Riemann hypothesis holds for the zeros with $0<\im\varrho\le T$ and
$4.92\sqrt{x/\log x}\le T$. With the height $T=3\cdot10^{12}$ of Platt and Trudgian \cite{PT21} both hold for $x\le2.169\cdot10^{25}$
\cite[Proposition~3.1]{Jo}, and Johnston extends them to $x\le1.101\cdot10^{26}$ \cite[Corollary~3.3]{Jo}.
\item $y>1.101\cdot10^{26}$: $E_{\mathrm{dn}}(y)=m_0y$ and $E_{\mathrm{up}}(y)=M_0y$ with $(m_0,M_0)=(1.1344\cdot10^{-9},\,1.1199\cdot10^{-9})$ for
$y>e^{50}$ and $(1.2312\cdot10^{-11},\,1.2215\cdot10^{-11})$ for $y>e^{60}$, by Broadbent, Kadiri, Lumley, Ng and Wilk
\cite[Table~14]{BKLNW}: $(1-m_0)x<\theta(x)<(1+M_0)x$ for $x>X_0$. (The same paper gives $\theta(x)\le(1+1.93378\cdot10^{-8})x$ for all $x\ge0$
\cite[Corollary~2.1]{BKLNW}; the formulas above do not need it, but the formal proof of Section~\ref{ss:eff-formal} does.)
\item $x>1.101\cdot10^{26}$: $\psi_{\mathrm{up}}(x)=E_{\mathrm{up}}(x)+2\sqrt x\log x$ and $\psi_{\mathrm{dn}}(x)=E_{\mathrm{dn}}(x)$, by
$0\le\psi(x)-\theta(x)\le2\sqrt x\log x$ (Appendix~\ref{app:mult}, in the form taken from Mathlib in Section~\ref{sec:formal}).
\end{itemize}
Any $K$ with $4K^2\le n$ may be used. The computation takes, for each node of Appendix~\ref{app:eff-assembly}, the best of nine candidates
$K=\min(\max(\fl{c\,n^{1/3}},2),\fl{\sqrt n/2},K_{\max})$, $c\in\{0.03,0.06,0.1,0.15,0.25,0.4,0.7,1.2,2.0\}$, with $n$ the first value of the
node and $K_{\max}=10^7$ ($19\,998$ in the run for Theorem~\ref{thm:E1}), namely the one with the smallest estimate of $U$ at the end of the node. This gives $K=4641$ for the node beginning at $n=10^{14}+1$, $K=9999$
for the node containing $10^{15}$ and $K=464\,146$ for the node containing $10^{20}$. The error terms are not monotone in $y$ at $y=e^{60}$ (the row $\log X_0=60$ is
smaller than the row $50$), so each is replaced by its monotone envelope $\sup_{y'\le y}E(y')$ before it is evaluated at $n_{\mathrm{hi}}$. In the window sums of the second range the program bounds $\log^2(\beta_kn)$ and $\log^2(\alpha_kn)$ by $\log^2(3n/2)$.

\subsection{The continued fraction}\label{app:eff-cf}
\emph{Proof of Theorem~\ref{thm:E3}\textup{(iii)}.} A ball enclosure of $\pi$ computed with the library Arb at a precision of $2.6\cdot10^8$ bits gives two
rationals $L<\pi<U$, and Euclid's algorithm, run on $L$ and $U$ in lockstep, returns their common prefix $a_0,\dots,a_K$, $K=75\,928\,239$
($15$ minutes, $8$\,GB). The set of real numbers whose
continued fraction begins with given partial quotients is an interval, so this is also a prefix of the continued fraction of $\pi$. The list
certifies itself: from it alone one recomputes $p_K$, $q_K$, $p_{K-1}$, $q_{K-1}$, checks $p_Kq_{K-1}-p_{K-1}q_K=(-1)^{K+1}$, and checks by
exact cross-multiplication that $L$ and $U$ lie strictly between $p_K/q_K$ and $(p_K+p_{K-1})/(q_K+q_{K-1})$, which forces $\pi=[a_0;a_1,\dots,
a_K,y]$ with $y>1$. An independent enclosure computed with MPFR gives the same prefix. For $1\le k<K$, $q_k\le b<q_{k+1}$ and every $a$,
the best-approximation property of the convergents gives
\[
\Bigl|\pi-\frac ab\Bigr|\ge\frac{|q_k\pi-p_k|}{b},\qquad\text{hence}\qquad
\Bigl|\pi-\frac ab\Bigr|>b^{-\mu}\quad\text{whenever}\quad\Bigl|\pi-\frac{p_k}{q_k}\Bigr|>q_k^{-\mu}
\]
(since then $|q_k\pi-p_k|>q_k^{1-\mu}\ge b^{1-\mu}$ for $\mu>1$). The exponent $e_k=-\log|\pi-p_k/q_k|/\log q_k$ is $3.42928833728180\ldots$ for $k=1$ ($22/7$),
$3.20195874260656\ldots$ for $k=3$ ($355/113$), and less than $3.2023$ for every $2\le k<K$; the pair $b=q_1=7$, $a=p_1=22$ shows that the
condition $\mu>e_1$ is also necessary. For $2\le b\le6$ a direct evaluation shows that $\mu>2.82018168049$ suffices. (For the application only
$\mu=6.045$ matters: the cruder test $(6.045-2)\log q_k>\log(a_{k+1}+2)$, which uses $|q_k\pi-p_k|>1/((a_{k+1}+2)q_k)$, holds for all
$75\,928\,238$ indices $1\le k<K$.)

\subsection{Assembly and the tail}\label{app:eff-assembly}
The program \texttt{verify\_assembly.py} evaluates the expressions of Appendix~\ref{app:eff-lemma} for the bounds of Theorem~\ref{thm:E3} on
\emph{nodes}, ranges of $n$ that carry one pair of bounds, affine in $n$ within the node: the $8\,059\,058$ blocks of Theorem~\ref{thm:E3}(ii)
and a geometric grid from $10^{14}$ to $10^{37}$ (of relative width $10^{-4}$ up to $10^{20}$ and $10^{-3}$ beyond, refined to $10^{-5}$
around the values of $n$ that determine the rows of Table~\ref{tab:b0}; in the run for Theorem~\ref{thm:E1} the grid has relative width
$10^{-5}$ and ends at $1.2\cdot10^{14}$), on which
$U(n)\le(4-\varpi_K)n+R(n_{\mathrm{hi}})$ and $V(n)\ge(4-\varpi_K)n-R'(n_{\mathrm{hi}})$ with the increasing remainders $R$, $R'$ of
\eqref{eq:EUV}. The sequences are used from $n=10^7$ on (Lemma~\ref{lem:E} applied to $q_{n+10^7-1}$). The expressions are evaluated in Arb at
$256$ bits, comparisons are certified for all members of the balls, and the bounds are stored as outward-rounded doubles. The logarithmic
bounds $e_n$ and $\ell^+_n$ are replaced by nonincreasing and nondecreasing majorants, respectively, computed by backward and forward maxima
with outward rounding across nodes. These bounds define the sequences used in Lemma~\ref{lem:E}. The certificate is this: if $\hat n(b)=n$, then
$\log2b\ge\min(-e_{n-1},g_n)$ by the two cases in the definition of $\hat n$, hence
$\mu_{\mathrm{eff}}(b)\le M_n:=1+(\log2+\ell^+_n)/(\min(-e_{n-1},g_n)-\log2)$; for $b\ge q_K$ one has $N_1(b)\ge n^*=31\,907\,390$ (because $2q_KE_{n^*-1}>1$), so
$\sup_{b\ge q_K}\mu_{\mathrm{eff}}(b)\le\max_{n\ge n^*}M_n$; on each node $M_n$ is the larger of two ratios of piecewise affine functions of $n$,
whose maximum is attained at the ends of the pieces; and the $n$ with $\min(-e_{n-1},g_n)\ge\log(2\cdot10^{122\,680\,112\,076\,457})$ are
excluded, since for the corresponding $b$ the formal theorem applies. In the run for Theorem~\ref{thm:E1} this disposes of all $n>1.2\cdot10^{14}$ at once: with $K=4932$ windows and the bounds \eqref{eq:EUV} taken uniformly in $n$, as for the tail below, every $b$ with $\hat n(b)>1.2\cdot10^{14}$ has $\log2b\ge3.388\cdot10^{14}>\log(2\cdot10^{122\,680\,112\,076\,457})$. Over the steps $n\ge n^*$ of the blocks, the maximum of $M_n$ is $6.045435237\pm10^{-10}$, in the block
beginning at $n=42\,140\,967$; for the node containing $n$ the bound is $6.04489$ at $n=10^8$, $6.04491$ at $10^{11}$,
$6.04496$ at $10^{12}$, $6.04487$ at $10^{14}$, and on the closed forms it decreases from $6.04496$ at $10^{14}+1$ to $6.044667$ at
$10^{16}$, $6.044606$ at $10^{18}$ and $6.0445885$ at $10^{28}$ (rounded up; the values for every decade are in the recorded run).

\emph{The two programs.} An earlier program, \texttt{assemble.py}, assembled the first tables ($921\,039$ blocks of width $10^{-5}$,
$K=1000$, for $n\le10^{11}$) in double precision with a different argument (the supremum of $\mu_{\mathrm{eff}}$ over the steps of a node
through $T_n=\min(-e_{n-1},\inf_{m\ge n}g_m)$). On those tables both programs give $\sup_{b\ge q_K}\mu_{\mathrm{eff}}(b)\le6.048179$, with the largest
node bound at the first node of the closed forms, then at $n=10^{11}+1$ (the earlier value of Theorem~\ref{thm:E1}); the comparison found two places where
the bounds of \texttt{assemble.py} are not rigorous as written, neither affecting a printed number (its remainders $R$, $R'$ were not monotone at
$y=e^{60}$, and the $99$ values $10^7<n<10^7+100$ carried no node). The finer tables and the hybrid range were assembled by the ball-arithmetic
program only; its treatment of the hybrid blocks was checked end to end on a small hybrid configuration, where the affine bounds it derives
were compared with the true values of $\log D_n$ from an independent sieve at three values of $n$ in each of thirteen blocks.

\emph{The tail.} The grid ends at $n=10^{37}$, that is, at $\log b=2.824\cdot10^{37}$. For every $n\ge N_T=10^{30}$ we use instead the fixed
number $K_T=10^7$ of windows ($4K_T^2\le N_T$) and $q_n=D_n^{(K_T)}v_n$, $p_n=-D_n^{(K_T)}u_n$. Every window endpoint is then at least
$N_T/K_T=10^{23}$, where the relative error of the estimates of Theorem~\ref{thm:E3}(ii) is at most $1.1344\cdot10^{-9}$ downwards and
$1.1199\cdot10^{-9}$ upwards (the row $\log X_0=50$ of \cite[Table~14]{BKLNW}; the error of the Riemann-hypothesis-type bound is smaller there).
Bounding the sublinear terms of $U$ and $V$, divided by $n$, by their suprema over $n\ge N_T$ gives, for all $n\ge N_T$,
\[
s_-n\le\log D_n^{(K_T)}\le s_+n,\qquad s_+=(4-\varpi)+9.38\cdot10^{-8},\qquad s_-=(4-\varpi)-6.72\cdot10^{-8},
\]
and with Theorem~\ref{thm:E3}(i)
\[
|q_n|\le e^{(s_++r')n},\qquad|q_n|\ge e^{(s_-+r-3.7\cdot10^{-29})n},\qquad|q_n\pi-p_n|\le e^{-(-\log\rho-s_+)n+\log4}.
\]
The following exponent form of Lemma~\ref{lem:E} with explicit constants is proved in Lean (Section~\ref{ss:eff-formal}). \emph{Let $\theta$ be
irrational, $p_n,q_n\in\Z$, $0\le\sigma^+$, $\sigma^-\le\sigma^+<\sigma^-+\tau$, $c_E,c^\pm\ge0$, and suppose that $|q_n\theta-p_n|\le e^{-\tau n+c_E}$ and
$e^{\sigma^-n-c^-}\le|q_n|\le e^{\sigma^+n+c^+}$ for all $n\ge n_1$. Put $\kappa=1+\sigma^+/(\tau-(\sigma^+-\sigma^-))$ and
$\Lambda=\kappa\log2+\sigma^+\bigl(n_1+2+(c_E+c^++c^-+\sigma^+)/(\tau-(\sigma^+-\sigma^-))\bigr)+c^+$. Then $|\theta-a/b|>b^{-\mu}$ for every
$\mu>\kappa$, every $a$ and every $b\ge\lceil\exp(\Lambda/(\mu-\kappa))\rceil+1$.} Here $n_1=N_T$, $c_E=\log4$, $c^\pm=0$,
$\kappa=6.044588604\ldots$ and $\Lambda=1.4246\ldots\cdot10^{31}$, and the threshold for $b$ lies inside the grid for every $\mu\ge6.0445892$, in particular
for every row of Table~\ref{tab:b0}. So the grid and the tail together cover every $b$.

\begin{proof}[Proof of Theorem~\textup{\ref{thm:E2}}]
The exponent form is applied with $n_1=N_T=10^{30}$ and $K=2^{19}$ windows to the bounds of Theorem~\ref{thm:E3}(i),(ii) in the form proved in
Lean (Section~\ref{ss:eff-formal}): $U$ and $V$ with the uniform error terms obtained from (H1)--(H4), and $c_n$ between
$e^{n\log\Lambda_{\mathrm{lo}}-\frac12\log n-3}$ and $\Lambda_{\mathrm{hi}}^n$. With the sublinear terms taken at $N_T$ this gives
$\tau=6.5301349-s_+$, $\sigma^+=s_++\log\Lambda_{\mathrm{hi}}$ and $\sigma^-=s_-+\log\Lambda_{\mathrm{lo}}-(3+\log2+\frac12\log N_T)/N_T$,
and the numerical bounds for $\varpi_K$, $\log\Lambda_{\mathrm{lo}}$, $\log\Lambda_{\mathrm{hi}}$ and the sublinear terms give
$\kappa\le6.04459707$ and $\Lambda\le1.4246398\cdot10^{31}$, which is the statement. This is the proof carried out in Lean
(\texttt{pi\_effective\_closed}; the arithmetic is in the lemmas \texttt{bounds\_r30}, \texttt{kappa\_r30} and \texttt{logCE\_r30} of the file
\texttt{MainClosed.lean}). The tail bound above, with the finer error terms of \eqref{eq:EUV}, gives in the same way the version with
$\kappa_T=6.044588604\ldots<6.0445887$ and $\Lambda=1.4246\ldots\cdot10^{31}<1.4247\cdot10^{31}$, that is $|\pi-a/b|>b^{-\mu}$ for every
$\mu>6.0445887$ and every $b\ge\lceil\exp(1.4247\cdot10^{31}/(\mu-6.0445887))\rceil+1$; this version is not machine-checked. (The recorded run of \texttt{verify\_assembly.py}, which bounds $\psi(4n)$ by the row $\log X_0=60$ instead of the row $50$, gives the slightly smaller value $\kappa_T=6.04458858\pm4\cdot10^{-9}$.)
\end{proof}

\begin{proof}[Proof of Theorem~\textup{\ref{thm:E1}} and of Table~\textup{\ref{tab:b0}}]
Let $b\ge2$. If $b<q_K$, Theorem~\ref{thm:E3}(iii) gives $|\pi-a/b|>b^{-\mu}$ for every $\mu>3.4293$. If $q_K\le b<10^{122\,680\,112\,076\,457}$,
Lemma~\ref{lem:E} with the bounds of Theorem~\ref{thm:E3}(i),(ii) gives $|\pi-a/b|\ge b^{-\mu_{\mathrm{eff}}(b)}$, and the computation above
certifies $\mu_{\mathrm{eff}}(b)\le6.045436$ for all these $b$, hence $|\pi-a/b|>b^{-6.0455}$. If $b\ge10^{122\,680\,112\,076\,457}$, the
instance $(K,n_1)=(2784,10^{14})$ of the sharpened exponent form (Table~\ref{tab:chain}) gives $|\pi-a/b|>b^{-6.04493479}$. This proves
Theorem~\ref{thm:E1}, which gives the first row of Table~\ref{tab:b0}. The other entries are upper bounds for
$1+\max\{b:\mu_{\mathrm{eff}}(b)\ge\mu\}$ over the nodes; for each such row the threshold of the finer version of Theorem~\ref{thm:E2} lies inside
the grid, so that every larger $b$ is covered.
\end{proof}

\begin{table}[ht]
\small
\begin{tabular}{lll}
\toprule
$\mu$ & $b_0(\mu)$ admissible & $b_0(\mu)$ of Theorem~\ref{thm:E2}\\
\midrule
$\ge6.0455$ & $2$ (Theorem~\ref{thm:E1}) & $\exp(1.58\cdot10^{34})$\\
$6.045$ & $10^{1.3636\cdot10^{11}}$ & $\exp(3.54\cdot10^{34})$\\
$6.0448$ & $10^{6.1975\cdot10^{14}}$ & $\exp(7.03\cdot10^{34})$\\
$6.0447$ & $10^{4.2232\cdot10^{15}}$ & $\exp(1.39\cdot10^{35})$\\
$6.0446$ & $10^{4.0137\cdot10^{18}}$ & $\exp(4.87\cdot10^{36})$\\
$6.04459$ & $10^{1.6254\cdot10^{26}}$ & ---\\
\bottomrule
\end{tabular}
\medskip
\caption{Admissible values of $b_0(\mu)$ from the range-by-range computation (exponents of $10$ rounded up in the last digit; certified in ball
arithmetic), against the closed formula of Theorem~\ref{thm:E2}. The row $6.0446$ is the exponent of Theorem~\ref{thm:main}. The row $6.04459$
lies below the exponent $6.04459707$ of Theorem~\ref{thm:E2}; it is covered by the finer, not machine-checked, tail bound above.}\label{tab:b0}
\end{table}

\emph{What limits the constants.} The exact multiplier fluctuates. The difference $d_n=\frac1n\log D_n-(4-\varpi)$, with all admissible windows,
reaches $+0.0419$ at $n=1712$, $+0.0187$ at $n=14311$, $+0.0058$ at $n=117180$ and $+0.0019$ at $n=1107624$, and the exponent of a single step
moves by about $2.14\,d_n$; the maxima decay roughly like $n^{-1/2}$. The bounds for the forms alone therefore give the upper bound $6.089054$ for
$b\ge e^{2.82\cdot10^4}$, but the continued fraction covers this range completely. Above $n^*=3.19\cdot10^7$, where the linear forms take over
from the continued fraction, the largest certified node bound is $6.045436$, at $n=4.2\cdot10^7$.
Doubling the length of the continued fraction (which costs twice the memory, $16$\,GB) would move $n^*$ to about $6.4\cdot10^7$, with an estimated
upper bound of $6.0453$. To prove the exponent $6.0446$ for all $b\ge2$ by this division of ranges, it would suffice to certify the continued-fraction
bound up to denominators of size $10^{4.0137\cdot10^{18}}$ (Table~\ref{tab:b0}); that computation is out of reach. In the range of the closed forms the node bound exceeds the floor $6.044588$ by $3.6\cdot10^{-4}$ at
$n=10^{14}$, $1.7\cdot10^{-5}$ at $10^{18}$ and $2\cdot10^{-6}$ at $10^{22}$ (roughly $15\,n^{-1/3}$ up to $10^{18}$).

\subsection{The Lean files of Section~\ref{sec:effective}}\label{app:eff-lean}
The eight files of \texttt{lean/MuPi/Effective/} (namespace \texttt{PiMeasure.Effective}; about $5{,}300$ lines in all; no \texttt{sorry}; axioms
\texttt{propext}, \texttt{Classical.choice}, \texttt{Quot.sound}) and their main theorems. The hypotheses \texttt{HypB1}, \texttt{HypB2},
\texttt{HypB4}, \texttt{HypB5} are (H1)--(H4) of Section~\ref{ss:eff-formal}; \texttt{HypB3} is defined but unused.
\begin{center}\small
\begin{tabular}{@{}>{\raggedright\arraybackslash}p{3.3cm}>{\raggedright\arraybackslash}p{10.4cm}@{}}
\toprule
File & Content\\
\midrule
\texttt{HataEffective.lean} & \texttt{lemmaE}, \texttt{lemmaE\_core} (Lemma~\ref{lem:E}; the core version needs its hypotheses only for
$n\ge N_1$, neither the monotonicity of $E_n$ nor the irrationality of $\theta$, and only that $E_n/Q^-_n<|\theta-a/b|$ at one index);
\texttt{measure\_of\_exp\_bounds}, \texttt{exponent\_of\_exp\_bounds}, \texttt{exponentLE\_of\_exp\_bounds} (the exponent form of
Appendix~\ref{app:eff-assembly}, with \texttt{kappaE} and \texttt{logCE}).\\
\texttt{HataEffective}\newline\texttt{Sharp.lean} & \texttt{measure\_of\_exp\_bounds\_sharp}, \texttt{exponent\_of\_exp\_bounds\_sharp}: the constant
\texttt{logCESharp} without the term $\sigma^+n_1$, valid for $\log2b\ge\tau n_1-c_E$.\\
\texttt{Recurrence.lean} & \texttt{cPS\_rec} (the recurrence \eqref{eq:Erec} for every $n\ge0$, from the certificate of
Appendix~\ref{app:eff-growth} as an identity of integer coefficient lists checked by kernel evaluation, and Mathlib's
\texttt{PowerSeries.derivative}); \texttt{cPS\_zero}, \dots, \texttt{cPS\_five}.\\
\texttt{GrowthLower.lean} & \texttt{cub\_LamLo}, \texttt{cub\_LamHi}, \texttt{exists\_root} (the bracket
$\Lambda_{\mathrm{lo}}<\lambda_3<\Lambda_{\mathrm{hi}}$); the band induction; \texttt{cPS\_lower} ($n\ge3$).\\
\texttt{MultBounds.lean} & \texttt{HypB1}, \dots, \texttt{HypB5}; \texttt{theta\_ge\_dn\_nat}, \texttt{psi\_ge\_dn\_nat} (the error terms at
integers); \texttt{log\_DmultK\_le}, \texttt{log\_DmultK\_ge} ($K\ge1$, $n\ge1423K$) with \texttt{Ubound}, \texttt{Vbound} and their closed
forms; \texttt{SBeta\_le}, \texttt{SAlpha\_le}, \texttt{Bsum\_le}, \texttt{Asum\_le}.\\
\texttt{VarpiBounds.lean} & two-sided bounds for $\varpi_K$, $K\ge20$ (\texttt{varpiK\_ge}, \texttt{varpiK\_le}).\\
\texttt{Main.lean} & \texttt{exists\_forms} (the forms $q_n$, $p_n$ satisfy the hypotheses of the exponent form for $n\ge n_1$, with the
constants $6.530134$ and $10.540357$ of Section~\ref{sec:formal}); \texttt{pi\_effective} (general $K$, $n_1$); \texttt{pi\_effective\_explicit}
($\mu>6.0481$), \texttt{pi\_effective\_6\_05\_pow10}, \texttt{pi\_effective\_6\_1\_pow10}, \texttt{pi\_measure\_6\_0481} (the bound with a constant for
every $b\ge1$, Appendix~\ref{app:eff-ladder}).\\
\texttt{MainClosed.lean} & (generated by a script of the accompanying material) \texttt{cPS\_upper} ($c_n\le\Lambda_{\mathrm{hi}}^n$, $n\ge3$),
\texttt{rho\_le\_sharp} ($\rho\le e^{-6.5301349}$), \texttt{varpiK\_ge\_sharp}; \texttt{exists\_forms\_closed},
\texttt{pi\_effective\_closed\_of}; the instances \texttt{bounds\_r7}, \dots, \texttt{bounds\_r30} with \texttt{kappa\_rX} and
\texttt{Lambda\_rX}, and \texttt{logCE\_r30}; \texttt{pi\_effective\_closed} (Theorem~\ref{thm:E2}), \texttt{pi\_effective\_6\_0446}; the ladder
\texttt{pi\_effective\_range\_r7}, \dots, \texttt{\_r16}, \texttt{\_r30} and their twins \texttt{\_pow10} (Table~\ref{tab:chain});
\texttt{pi\_effective\_6\_0455\_from} (Theorem~\ref{thm:E1} for $b\ge10^{122\,680\,112\,076\,457}$), \texttt{pi\_effective\_6\_0482\_from},
\texttt{pi\_effective\_6\_0446\_from}.\\
\bottomrule
\end{tabular}
\end{center}
The formal \texttt{lemmaE} is slightly more general than Lemma~\ref{lem:E}: $N_1$ may be any index with $2bE_{N_1}\le1$, and the decay
hypothesis only asks that $E_n/Q^-_n$ take arbitrarily small values.

\subsection{The machine-checked ladder, and what each statement rests on}\label{app:eff-ladder}
The formal argument of Section~\ref{ss:eff-formal}, applied to the forms with $K$ windows for $n\ge n_1$ through the sharpened exponent form,
gives for each $(K,n_1)$ an exponent valid from an explicit threshold on (Table~\ref{tab:chain}; between two consecutive thresholds the exponent of the
row with the smaller threshold applies), and the measure form of Lemma~\ref{lem:E} at $K=256$, $n_1=10^{11}$ gives a bound with a constant for every $b\ge1$
(theorem \texttt{pi\_measure\_6\_0481}):
\begin{equation}\label{eq:ladder}
\Bigl|\pi-\frac ab\Bigr|\ \ge\ e^{-1.424733\cdot10^{12}}\,b^{-6.0481}\qquad\text{for all integers $a$ and all integers $b\ge1$.}
\end{equation}
Together with the continued fraction of Theorem~\ref{thm:E3}(iii), the ladder alone gives $|\pi-a/b|>b^{-6.1320388}$ for all $b\ge2$; the
exponent $6.0455$ of Theorem~\ref{thm:E1} needs the prime table for $q_K\le b<10^{122\,680\,112\,076\,457}$.

\begin{table}[ht]
\small
\begin{tabular}{llll}
\toprule
$(K,n_1)$ & $\kappa(K,n_1)\le$ & threshold $b\ge$ & $|\pi-a/b|>b^{-\mu}$ with $\mu=$\\
\midrule
$(20,10^{7})$ & $6.1320347$ & $10^{12\,268\,013}$ & $6.1320388$\\
$(25,10^{8})$ & $6.0812933$ & $10^{122\,680\,113}$ & $6.0812937$\\
$(55,10^{9})$ & $6.0611564$ & $10^{1\,226\,801\,122}$ & $6.0611565$\\
$(118,10^{10})$ & $6.0521392$ & $10^{12\,268\,011\,209}$ & $6.0521393$\\
$(253,10^{11})$ & $6.0480523$ & $10^{122\,680\,112\,078}$ & $6.0480524$\\
$(500,10^{12})$ & $6.04618878$ & $10^{1\,226\,801\,120\,766}$ & $6.04618879$\\
$(1292,10^{13})$ & $6.04533031$ & $10^{12\,268\,011\,207\,647}$ & $6.04533032$\\
$(2784,10^{14})$ & $6.04493478$ & $10^{122\,680\,112\,076\,457}$ & $6.04493479$\\
$(5000,10^{15})$ & $6.04475227$ & $10^{1\,226\,801\,120\,764\,554}$ & $6.04475228$\\
$(10772,10^{16})$ & $6.04466800$ & $10^{12\,268\,011\,207\,645\,527}$ & $6.04466801$\\
$(2^{19},10^{30})$ & $6.04459707$ & $10^{1.2269\cdot10^{30}}$ & $6.04459708$\\
\bottomrule
\end{tabular}
\medskip
\caption{The machine-checked ladder (theorems \texttt{pi\_effective\_range\_r7\_pow10}, \dots, \texttt{\_r16\_pow10}, \texttt{\_r30\_pow10},
under (H1)--(H4)): for every integer $a$ and every integer $b$ at least the threshold, $|\pi-a/b|>b^{-\mu}$. The threshold is a power of $10$
above $\lceil e^{T}\rceil+1$ with $T=\lceil2.82409618\,n_1\rceil\ge\tau n_1$; in the last row the file states the exact power
$10^{1226801120764552562988705473503}$. The numbers $K$ of windows (at least $20$, the range in which the bounds for $\varpi_K$ are proved) were chosen to make the certified $\kappa$ small: each $\kappa$ of the table is within $3\cdot10^{-6}$ of the smallest value that the same chain gives for its $n_1$.}\label{tab:chain}
\end{table}

\begin{table}[ht]
\small
\begin{tabular}{@{}>{\raggedright\arraybackslash}p{3.7cm}>{\raggedright\arraybackslash}p{3.3cm}>{\raggedright\arraybackslash}p{8.3cm}@{}}
\toprule
Statement & Lean & Outside Lean\\
\midrule
Theorem~\ref{thm:E1}, $2\le b<q_K$ & --- & continued fraction (Appendix~\ref{app:eff-cf}): self-certifying list of partial quotients, two enclosures of $\pi$\\
Theorem~\ref{thm:E1}, $q_K\le b<10^{122\,680\,112\,076\,457}$ & Lemma~\ref{lem:E}; growth bounds of Section~\ref{ss:eff-formal} & (H1), (H2) and $\theta(x)<x$ for $x\le10^{19}$ \cite[Theorem~2, (1.7)]{Bu18}, and for the indices $n>1.2\cdot10^{14}$ the estimates of \cite{Bu16,Jo,BKLNW} beyond $10^{19}$; prime table to $4\cdot10^{11}$ and block tables (Appendix~\ref{app:eff-mult}); assembly in ball arithmetic (Appendix~\ref{app:eff-assembly})\\
Theorem~\ref{thm:E1}, $b\ge10^{122\,680\,112\,076\,457}$ & yes, under (H1)--(H4) & ---\\
Theorem~\ref{thm:E2} & yes, under (H1)--(H4) & ---\\
Ladder, Table~\ref{tab:chain}, and \eqref{eq:ladder} & yes, under (H1)--(H4) & ---\\
Table~\ref{tab:b0} & Lemma~\ref{lem:E}; growth bounds of Section~\ref{ss:eff-formal}; the exponent form & (H1), (H2), $\theta(x)<x$ for $x\le10^{19}$ \cite[Theorem~2, (1.7)]{Bu18}, and \cite{Bu16,Jo,BKLNW} beyond $10^{19}$; prime table and block tables; assembly in ball arithmetic; Theorem~\ref{thm:E3}(i) for the finer tail bound\\
\bottomrule
\end{tabular}
\medskip
\caption{What each statement rests on. ``Lean'' means a theorem of the formal development of Section~\ref{ss:eff-formal}, whose only hypotheses
are the published estimates (H1)--(H4) where indicated. The computations are described in this appendix with the cross-checks that were made;
the second and the last row are the only parts of the section that depend on the prime table.}\label{tab:trust}
\end{table}

\section*{Acknowledgements}
We thank Philipp Habegger and Wadim Zudilin for valuable feedback on the manuscript.

\end{document}